\documentclass[12pt]{article}
\usepackage{amsfonts,amssymb,latexsym}

\newtheorem{theorem}{Theorem}[section]
\newtheorem{lemma}[theorem]{Lemma}
\newtheorem{corollary}[theorem]{Corollary}
\newtheorem{proposition}[theorem]{Proposition}
\newtheorem{definition}[theorem]{Definition}
\newenvironment{proof}
{\par\addvspace{0.3cm}\noindent{\rm Proof. }}
{\nopagebreak\mbox{}\hfill $\Box$\par\addvspace{0.25cm}}

\newcommand{\reseteqn}{\setcounter{equation}{0}}

\newcommand{\C}{{\mathbb C}}
\newcommand{\Co}{\dot{\C}}
\newcommand{\Z}{{\mathbb Z}}
\newcommand{\D}{{\mathbb D}}
\newcommand{\R}{{\mathbb R}}
\newcommand{\sP}{{\mathbb P}}
\newcommand{\cP}{{\cal P}}
\newcommand{\ep}{\varepsilon}
\newcommand{\ka}{\varkappa}

\newcommand{\La}{\Lambda}
\newcommand{\Ga}{\Gamma}
\newcommand{\ga}{\gamma}
\newcommand{\bqn}{\begin{eqnarray}}
\newcommand{\eqn}{\end{eqnarray}}
\newcommand{\ben}{\begin{equation} }
\newcommand{\een}{\end{equation} }
\newcommand{\ba}{\begin{array}}
\newcommand{\ea}{\end{array}}
\newcommand{\nn}{\nonumber}
\newcommand{\iv}{^{-1}}
\newcommand{\iy}{\infty}
\newcommand{\ovl}{\overline}

\renewcommand{\Re}{{\rm Re\,}}

\newcommand{\im}{{\rm im\,}}
\newcommand{\id}{{\rm id}}
\newcommand{\diag}{{\rm diag}\,}
\newcommand{\wt}[1]{\widetilde{#1}}
\newcommand{\wh}[1]{\widehat{#1}}
\newcommand{\mdg}{{\rm min}}
\newcommand{\abs}[1]{\vert#1\vert}
\newcommand{\tz}{\tilde{z}}
\newcommand{\tD}{\wt{D}}
\newcommand{\tS}{\wt{S}}

\begin{document}

\date{ }
\title{Factorization of piecewise constant matrix functions
and systems of linear differential equations}

\author{Torsten Ehrhardt
\thanks{tehrhard@mathematik.tu-chemnitz.de.
Research at MSRI is supported in part by NSF grant DMS-9701755.
Research has also been supported in part by DAAD Grant 213/402/537/5.}\\
Fakult\"{a}t f\"{u}r Mathematik\\
Technische Universit\"{a}t Chemnitz\\
09107 Chemnitz, Germany
\and
Ilya M. Spitkovsky
\thanks{ilya@math.wm.edu. Research is supported in part by NSF
Grant DMS-9800704.}\\
Department of Mathematics\\
College of William and Mary\\
Williamsburg, VA 23187, USA}

\maketitle

\begin{abstract}
Let $G$ be a piecewise constant $n\times n$ matrix function which
is defined on a smooth closed curve $\Ga$ in the complex sphere
and which has $m$ jumps. We consider the problem of determining
the partial indices of the factorization of the matrix function
$G$ in the space $L^p(\Ga)$. We show that this problem can be
reduced to a certain problem for systems of linear differential
equations.

Studying this related problem, we obtain some results for the partial
indices for general $n$ and $m$.  A complete answer is given for
$n=2$, $m=4$ and for $n=m=3$.  One has to distinguish several cases.
In some of these cases, the partial indices can be determined
explicitly.  In the remaining cases, one is led to two possibilities
for the partial indices.  The problem of deciding
which is the correct possibility is equivalent to the description of
the monodromy of $n$-th order linear Fuchsian differential equations
with $m$ singular points.
\end{abstract}

\section{Introduction}

Let $\Ga$ be a smooth closed positively oriented curve in the
complex sphere $\Co=\C\cup\{\iy\}$ which divides $\Co$ into two
domains $D_+$ and $D_-$. Assume without loss of generality that
$0\in D_+$ and $\iy\in D_-$.
Let $L^p(\Ga)$ be the Banach space of all
Lebesgue measurable functions on $\Ga$ for which
\bqn
\|f\|_{L^p(\Ga)} &=&
\Big(\int_{\Ga} |f(\tau)|^p\, |d\tau|\Big)^{1/p}
\quad<\quad\iy.
\eqn
For $1<p<\iy$, consider the singular integral operator $S_{\Ga}$
defined on $L^p(\Ga)$,
\bqn
(S_{\Ga}f)(t) &=&
\frac{1}{\pi i}\int_{\Ga} \frac{f(\tau)}{\tau-t}\,
d\tau,\qquad t\in\Ga.
\eqn
The operator $S_{\Ga}$ is bounded on  $L^p(\Ga)$ and
satisfies $S_{\Ga}^2=I$ (see \cite{GKru92}). Introduce the projections
\ben
P_{\Ga} \quad=\quad (I+S_{\Ga})/2,\qquad
Q_{\Ga} \quad=\quad (I-S_{\Ga})/2,
\een
and the Banach spaces
\ben
L_+^p(\Ga) \quad=\quad P_{\Ga}(L^p(\Ga)),\qquad
L_-^p(\Ga) \quad=\quad Q_{\Ga}(L^p(\Ga))\dotplus\C,
\een
i.e., the image of $P_{\Ga}$ and the image of $Q_{\Ga}$ plus the
set of constant functions.

It is well known (\cite{Gol}, see also \cite{LiSp}) that functions $f\in
L^p_+(\Ga)$ can be identified with functions $\hat{f}$ which are analytic in
$D_+$ and for which there exists an expanding sequence of domains $D_k$ with
rectifiable boundaries $\Ga_k$ such that $D_k\cup\Ga_k\subset D_+$, $\bigcup_k
D_k=D_+$, and $\sup_k\int_{\Ga_k}\abs{f(t)}^p\,\abs{dt}<\infty$. In fact,
$D_k$ can always be chosen as the image of $\{\,z\,\colon\,\abs{z}<r_k\,\}$
under a conformal mapping of the unit disk $\D=\{\,z\,\colon\,\abs{z}<1\,\}$
onto $D_+$, where $r_k$ is an arbitrary increasing sequence of positive
numbers converging to $1$. The function $f$ is given by the non-tangent limits
of the function $\hat{f}$ at the boundary $\Ga$ of $D_+$. Conversely, the
function $\hat{f}$ represents the analytic extension of $f$ into $D_+$. A
similar result holds for  $L^p_-(\Ga)$.

Let $L^\iy(\Ga)$ denote the set of all Lebesgue measurable
and essentially bounded functions on $\Ga$.
A {\em factorization} of an $n\times n$ matrix function
$G\in L^\iy(\Ga)^{n\times n}$ in the space $L^p(\Ga)$ is a
representation of the form
\bqn \label{fac}
G(t) &=& G_+(t) \La(t) G_-(t), \qquad t\in\Ga,
\eqn
where $\La(t)=\diag(t^{\ka_1},\dots,t^{\ka_n})$ is a diagonal
matrix with $\ka_1,\dots,\ka_n\in\Z$,
and the factors $G_+$ and $G_-$ satisfy the following conditions:
\begin{itemize}
\item[(i)]
$G_+\in L^p_+(\Ga)^{n\times n}$,
$G_+\iv\in L^q_+(\Ga)^{n\times n}$,
\item[(ii)]
$G_-\in L^q_-(\Ga)^{n\times n}$,
$G_-\iv\in L^p_-(\Ga)^{n\times n}$.
\end{itemize}
Here $1/p+1/q=1$.
The integers $\ka_1,\dots,\ka_n$ are called the {\em partial
indices} of the factorization.

Due to (i)-(ii), the operator $G_-\iv Q_\Ga G_+\iv$ is defined on the
(dense in $L^p(\Gamma)^n$) set of rational vector functions with poles
off $\Gamma$ and maps this set into $L^1(\Gamma)^n$. If this mapping
is bounded in the $L^p$-norm, then it can be extended by continuity to a
bounded operator on the whole space $L^p(\Gamma)^n$. In this case,
the representation (\ref{fac}) is called a {\em $\Phi$-factorization}
in the space $L^p(\Ga)$.

A detailed discussion of the $\Phi$-factorization, its properties
and the history of the subject can be found in \cite{CG,LiSp};
the latter monograph touches also on the (not necessarily $\Phi$-)
factorization in $L^p(\Gamma)$. The proofs of all the subsequent
results stated in this section can be found in \cite{CG,LiSp}.
The next theorem explains why the notion of the $\Phi$-factorization
is natural and important.

\begin{theorem}\label{t1.1}
The matrix function $G\in L^\iy(\Ga)^{n\times n}$ admits a
$\Phi$-factorization in the space $L^p(\Ga)$ if and only if the
singular integral operator $A=P_\Ga+GQ_\Ga$ is a
Fredholm operator on the space $L^p(\Ga)^n$.
In this case, we have
$$
\dim\ker A \quad=\quad \sum_{\ka_j>0} \ka_j,
\qquad\quad
\dim\ker A^* \quad=\quad -\sum_{\ka_j<0} \ka_j.
$$
\end{theorem}

Hence the partial indices $\ka_1,\dots,\ka_n$ contain important
information about the operator $A$. In particular, $A$ is
invertible if and only if the function $G$ admits a
so-called {\em canonical $\Phi$-factorization}, i.e.,
a $\Phi$-factorization
where all partial indices are zero. If $A$ is a Fredholm operator,
then its index, ${\rm ind\,}A:=\dim\ker A-\dim\ker A^*$, is
equal to
\bqn
\ka &=& \sum _{j=1}^n \ka_j,
\eqn
the {\em total index} of the factorization.

The notion of $\Phi$-factorization plays also an important
role in the study of more general singular integral operators.
One can show that the operator
$G_1P_\Ga+G_2Q_\Ga$  with $G_1,G_2\in L^\iy(\Ga)^{n\times n}$
is a Fredholm operator if and only if the functions $G_1$ and $G_2$
are invertible in $L^\iy(\Ga)^{n\times n}$ and the function
$G_1\iv G_2$ admits a $\Phi$-factorization.

The factorization of a matrix function is not unique. However, the
partial indices are always uniquely determined up to the change of
order. For this reason, we will assume that \ben
\ka_1\;\geq\;\ka_2\;\geq\;\dots \;\geq\;\ka_n. \een The relation
between factors $G_+$ and $G_-$ corresponding to different
factorizations of a function $G$ is described in the following
theorem. The characterization involves certain rational block
triangular matrix functions whose size is determined by the
multiple occurrence of same values
for the partial indices. 
Let $I_l$ denote the identity matrix of size $l\times l$,
and let $G\C^{l\times l}$ stand for the group of all
invertible $l\times l$ matrices.

\begin{theorem}\label{t1.2}
Assume that $G\in L^\iy(\Ga)^{n\times n}$ admits two
factorizations in $L^p(\Ga)$,
$$
G(t) \;\;=\;\; G_+(t)\La(t)G_-(t) \;\;=\;\;
\wt{G}_+(t)\La(t)\wt{G}_-(t),\qquad t\in\Ga,
$$
where
$\La(t)=\diag(
t^{\ovl{\ka}_1} I_{l_1},
t^{\ovl{\ka}_2} I_{l_2},\dots,
t^{\ovl{\ka}_k} I_{l_k})$
with $\ovl{\ka}_1,\dots,\ovl{\ka}_k\in\Z$,
$\ovl{\ka} _1>\ovl{\ka}_2>\dots>\ovl{\ka}_k$,
$l_1,\dots,l_k\in\{1,2,\dots\}$, $l_1+\dots+l_k=n$.
Then
\ben\label{f1.8}
G_+(t)\;\;=\;\; \wt{G}_+(t)V(t), \qquad\quad
\wt{G}_-(t)\;\;=\;\; U(t)G_-(t),
\een
where $U(t)$ and $V(t)$ are matrix functions of the form
$$
V(t)\;=\;\left(\ba{cccc}
A_{11} & V_{12}(t) & \dots & V_{1k}(t) \\
0 & A_{22} & & V_{2k}(t) \\
\vdots & & \ddots & \vdots \\
0 & 0 & \dots & A_{kk} \ea\right),
\quad
U(t)\;=\;\left(\ba{cccc}
A_{11} & U_{12}(t) & \dots & U_{1k}(t) \\
0 & A_{22} & & U_{2k}(t) \\
\vdots & & \ddots & \vdots \\
0 & 0 & \dots & A_{kk} \ea\right),
$$
with $A_{jj}\in G\C^{l_j\times l_j}$,
$V_{ij}(t) = \sum_{r=0}^{\ovl{\ka}_i-\ovl{\ka}_j}
A_{ij}^{(r)} t^r$,
$U_{ij}(t) = t^{\ovl{\ka}_j-\ovl{\ka}_i}V_{ij}(t)$
and $A_{ij}^{(r)}\in\C^{l_i\times l_j}$.
\end{theorem}

In fact, also the converse of the above theorem is true.
Suppose we are given a factorization ($\Phi$-factorization)
$G=G_+\La G_-$, and we determine $\wt{G}_+$ and
$\wt{G}_-$ by (\ref{f1.8}), then also
$G=\wt{G}_+\La \wt{G}_-$ is a factorization
($\Phi$-factorization).
Using this theorem and its converse, we can conclude
the following statement.

\begin{corollary}\label{c1.2a}
Suppose $G\in L^\iy(\Ga)^{n\times n}$ admits a
$\Phi$-factorization and $G=G_+\La G_-$ is a
factorization. Then $G=G_+\La G_-$ is
such a $\Phi$-factorization.
\end{corollary}

Let $PC(\Ga)^{n\times n}$ stand for the set of all
{\em piecewise continuous $n\times n$ matrix functions},
i.e., functions for which the one-sided limits $G(\tau+0)$
and $G(\tau-0)$ exist at each $\tau\in\Ga$. Here
$G(\ga(\theta_0)\pm 0):=\lim_{\theta\to\theta_0\pm 0}
G(\ga(\theta))$, where the periodic function $\ga:\R\to\Ga$ is
a parameterization of the positively oriented curve $\Ga$.
The following result solves the problem of the
existence of a $\Phi$-factorization for functions in
$PC(\Ga)^{n\times n}$.

\begin{theorem}\label{t1.3}
The matrix function $G\in PC(\Ga)^{n\times n}$ admits a
$\Phi$-factorization in the space $L^p(\Ga)$ if and only
if the following is satisfied:
\begin{itemize}
\item[(i)]
The matrices $G(\tau+0)$ and $G(\tau-0)$ are invertible for
each $\tau\in\Ga$.
\item[(ii)]
For each $j=1,\dots,n$ and each $\tau\in\Ga$, we have
\ben\label{f1.8x}
 \frac{1}{2\pi}\arg\lambda_j(\tau)+\frac{1}{p}
\;\notin\;\Z,
\een
where $\lambda_1(\tau),\dots,\lambda_n(\tau)$ are the eigenvalues
of the matrices $G(\tau-0)G(\tau+0)\iv$.
\end{itemize}
\end{theorem}

If the assumptions of this theorem are satisfied, then one can
define the numbers $\zeta_j(\tau)$ as the (unique) values of
$-(2\pi)\iv \arg\lambda_j(\tau)$ which lie in the interval
$J_p:=(1/p-1,1/p)$.

It is also possible to evaluate the total index $\ka$ related
to the above factorization. For simplicity assume that $G$ has
only finitely many jumps at points $a_1,\dots,a_m\in\Ga$, which
lie in this order on the positively oriented curve $\Ga$. Then
\bqn\label{f1.9}
\ka &=&
\sum_{k=1}^m \left[\frac{1}{2\pi}\arg\det G(\tau)
\right]_{\tau=a_k+0}^{a_{k+1}-0}\,+\,
\sum_{k=1}^m \sum_{j=1}^n \zeta_j(a_k).
\eqn
Here $a_{m+1}=a_1$, and $[\dots]$ denotes the total increment
along the subarcs on which $G$ is continuous.
Note that, in general, the value of $\ka$ depends on the space $L^p(\Ga)$
because so do the numbers $\zeta_j(\tau)$. Only in the case
where all $\lambda_j(\tau)$ are positive real numbers
(hence $\zeta_j(\tau)=0$), the value of $\ka$ does not depend on
the underlying space $L^p(\Ga)$. 

The calculation of the partial indices $\ka_1,\dots,\ka_n$ in the
case $n>1$ is a more delicate problem and can be solved only in
special situations. It is the purpose of this paper to consider
the case of {\em piecewise constant matrix functions}, i.e.,
functions which have only a finite number of jump discontinuities
and which are constant along the subarcs joining these points. We
will show that the problem of determining the partial indices can
be reduced to a certain problem for systems of linear differential
equations. Elaborating on this related problem, we will obtain
some information about the partial indices.

Preparing the following sections, we define certain matrices.
Let $G\in PC(\Ga)^{n\times n}$ be a piecewise constant matrix
function with $m$ jumps at the points $a_1,\dots,a_m\in\Ga$,
which are arranged on $\Ga$ in this order.
Assume that $G$ admits a
$\Phi$-factorization in the space $L^p(\Ga)$.
Introduce the matrices
\bqn
M_k &=& G(a_k-0)G(a_k+0)\iv,\qquad\quad
k=1,\dots,m.
\eqn
Because $G$ is constant along the subarcs joining
$a_k$ and $a_{k+1}$ and takes the values $G(a_k+0)=G(a_{k+1}-0)$
there, these matrices satisfy the relation
\bqn
M_1M_2\cdots M_m &=& I.
\eqn
Moreover, let $E_1,\dots,E_m$ be matrices such that
\bqn
M_k &\sim& \exp(-2\pi iE_k),\qquad k=1,\dots,m,\quad
\eqn
and the real parts of the eigenvalues of $E_k$ are contained
in the interval $J_p$. Here ``$\sim$'' stands for similarity
of matrices. Because $G$ admits a $\Phi$-factorization and
$J_p$ has length one, the matrices $E_k$ exist and are uniquely
determined up to similarity.
In fact, the numbers $\zeta_1(a_k),\dots,\zeta_n(a_k)$
defined above are equal to
the real parts of the eigenvalues of $E_k$. Formula
(\ref{f1.9}) for the total index can also be simplified:
\bqn\label{f1.13}
\ka &=& \sum_{k=1}^m {\rm trace\,} E_k.
\eqn
Note that the matrices $E_1,\dots, E_m$ depend on the
underlying space $L^p(\Ga)$.
Roughly speaking, they
contain the information about the ``proper'' choice
of the logarithm of the eigenvalues of $M_1,\dots,M_m$.
For sake of further reference,
the $m$-tuples of $n\times n$ matrices
\ben
[M_1,\dots,M_m]\quad \mbox{and} \quad
[E_1,\dots,E_m]
\een
will be called the ``{\em data}'' associated with the piecewise
constant matrix function $G$ with respect to the space $L^p(\Ga)$.

Let us recall what is already known about the factorization of
piecewise constant $n\times n$ matrix functions with $m$ jump discontinuities.
Any piecewise constant scalar function ($n=1$) can be easily factored
explicitly. Namely, $G(t)=G_+(t)t^{\ka}G_-(t)$, where
\ben\label{f1.16}
G_+(t)\;\;=\;\; c\prod_{k=1}^m (t-a_k)^{-\ep_k},\qquad
G_-(t)\;\;=\;\; \prod_{k=1}^m \left(1-\frac{a_k}{t}\right)^{\ep_k},
\een
$\ep_k=E_k$, and $\ka$ is given by (\ref{f1.13}). The branches of the
analytic functions in (\ref{f1.16}) are chosen in such a way that
$G_\pm$ is analytic on $D_\pm$, and
$c$ is a suitable nonzero constant. Hence, the case $n=1$ is solvable for any
$m$. Another easy case is $m=2$ for any $n$. Indeed, multiplying $G$ by the
constant matrix $G(a_1-0)^{-1}$ on the left, we may suppose without loss of
generality that $G$ assumes the values $I$ and $M(=M_2)$ only. These
values obviously commute with each other. In other words, $G$ is a so-called
{\it functionally commuting} matrix function, and its factorization can
be constructed explicitly (see \cite[Section 4.4]{LiSp}).

The simplest non-trivial case ($n=2$, $m=3$) was first studied by Zverovich
and Khvoschinskaya \cite{ZK} in a slightly different setting. In the setting
considered here it was treated by Tashbaev and one of the authors in
\cite{SpTa}. Using the appropriate modification of the results from \cite{ZK},
they gave explicit formulas for the factors and the partial indices. Therein
they had to distinguish several cases. Apart from trivial cases (where $G$ can
be reduced to a functionally commuting matrix function), the factors were
constructed by using the hypergeometric function. We will not restate here the
complete results of \cite{SpTa}, however, the results on the partial indices
will be obtained again in a corollary in Section \ref{s6} in the setting of
the related problem.

Let us describe the content of this paper in more detail. In Section 2,
we recall some basic facts from the theory of systems of linear differential
equations and from the theory of scalar linear differential equations
(of higher order).

In Section 3, we introduce a class of $n\times n$ systems
$Y'(\tz)=A(z)Y(\tz)$ of linear differential equations, which we
call systems of standard form. Their solutions $Y(\tz)$ are
characterized by the monodromy, which is given by
$[M_1,\dots,M_m]$, and by the behavior at the singularities
$a_1,\dots,a_m\in\C$, which is given by $[E_1,\dots,E_m]$. In
addition, we allow an apparent singularity at infinity, where the
behavior is described by certain integers $\ka_1,\dots,\ka_n$,
called the indices. Then we prove an ``existence'' and a
``uniqueness'' theorem for such systems. Moreover, we show that
these systems of linear differential equations provide us with a
solution to our factorization problem. Regarding the calculation
of the partial indices, we are led to the following (well-posed)
question: What are the indices of the systems which are
characterized by given data $[M_1,\dots,M_m]$ and
$[E_1,\dots,E_m]$ and singularities $a_1,\dots,a_m$? In fact, the
indices for the systems correspond to the partial indices of the
factorization.

In Section 4, we reformulate the conditions regarding the behavior at the
singularities at $a_1,\dots,a_m$ and at infinity in terms of the matrix
function $A(z)$ rather than in terms of the solution $Y(\tz)$. This is a simple
step, but all further considerations are based on it. We note that the above
question can now be expressed as follows: What is the monodromy of a system
given by $A(z)$, where $A(z)$ has a certain form?

In Section 5, we obtain two ``general'' results. Firstly, if the
monodromy is irreducible, then the indices (which are supposed to
be ordered decreasingly) satisfy $\ka_k-\ka_{k+1}\le m-2$ for each
$k=1,\dots, n-1$. Secondly, we show that this condition is sharp.
For this we produce a subclass of systems for which
$\ka_k-\ka_{k+1}= m-2$ for each $k=1,\dots, n-1$. It turns out
that the monodromy of such systems coincides with the monodromy of
certain scalar Fuchsian linear differential equations of $n$-th
order and with $m$ singularities. We remark that Bolibruch
\cite{Bo4} already obtained the first result in the setting of
vector bundles (and by using a different method) and that there is
also a strong connection of the second result with his work.

In Section 6, we consider the case $n=2$. One is led naturally to
distinguishing three cases relating to the ``reducibility type''
of $[M_1,\dots,M_m]$, i.e., the structure of the invariant
subspaces of these matrices. For general $m$, we can determine the
indices in some cases explicitly. In the remaining cases, we
obtain the estimate $\ka_1-\ka_2\le m-2$.  These remaining cases
include irreducible ones (for which this result is already clear
from Section 5) as well as certain reducible ones. The essential
point in the proof is a theorem which characterizes the data
(hence, in particular, the monodromy) of a certain class of
$2\times 2$ triangular systems.

With these results we are immediately in a position to give an explicit answer
for the case $n=2$, $m=3$, hence restating the results of \cite{SpTa}.
However, the results obtained here are slightly more general because
the assumptions on the matrices $E_1,\dots,E_m$ are less restrictive.

The complete answer is given for the case $n=2$ and $m=4$ by
resorting in addition to the second result of Section 5. However,
the answer is explicit only in some cases. In fact, in the other
cases, we are led to two possibilities for the indices. The
problem of deciding which possibility is the correct one is
equivalent to the description of the monodromy of second order
linear Fuchsian differential equations with four singular points.

In Section 7, we treat the case $n=m=3$. It is technically much more
complicated than the previous case as now nine reducibility types for
$[M_1,M_2,M_3]$ occur. We also have to prove theorems which characterize the
data (monodromy) of triangular and block-triangular systems. Finally, in the
main theorem we arrive at several cases, where, similar as before, in some
cases the indices can be given explicitly whereas in the other cases they
depend on the description of the monodromy of third order linear Fuchsian
differential equations with three singular points.

\reseteqn
\section{Systems of linear differential equations}
\label{s2}

A system of linear differential equations in the complex domain
can be written as \bqn\label{f2.1} \frac{dy}{dz} &=& A(z) y, \eqn
where $A(z)$ is a given $n\times n$ matrix function. We assume
that $A(z)$ is analytic on the complex plane $\C$ with the
exception of a finite number of points, which are called the {\em
singularities} of the system. The point $a=\iy$ is by definition a
singularity if the function $z^{-2}A(z^{-1})$ is not analytic at
$z=0$, which is justified by a change of variables $z\mapsto 1/z$.

Let $a_1,\dots,a_m\in\Co$ be the singularities of the
system, and let $\tS$ be the {\em universal covering surface}
of $S:=\Co\setminus\{a_1,\dots,a_m\}$.
Then the solutions $y$ of (\ref{f2.1}) are analytic functions
defined on $\tS$ and taking values in $\C^n$.
Associated with $\tS$ there is a {\em covering map}
$\rho:\tS\to S$.
Points on $\tS$ will usually be denoted by $\tz$ and points on $S$ by
$z=\rho(\tz)$. The set of all solutions of (\ref{f2.1}) forms
an $n$-dimensional linear space.

Taking $n$ linearly independent solutions
$y_1(\tz),\dots,y_n(\tz)$, one can consider the analytic
$n\times n$ matrix function $Y(\tz)=(y_1(\tz),\dots,y_n(\tz))$.
It satisfies the matrix equation
\bqn\label{f2.4}
Y'(\tz) &=& A(z)Y(\tz),\qquad \tz\in\tS,
\eqn
and we have $\det Y(\tz)\neq 0$ for all $\tz\in\tS$.
In fact, for any solution of (\ref{f2.4}) with
$\det Y(\tz_0)\neq0$ for some $\tz_0\in\tS$ we have
$\det Y(\tz)\neq 0$ for all $\tz\in\tS$.
We will consider only such solutions (i.e., analytic functions
$Y:\tS\to G\C^{n\times n}$) because
only they contain the full information.
Any two solutions $Y_1(\tz)$ and $Y_2(\tz)$
are related to each other by $Y_2(\tz)=Y_1(\tz)C$ with
$C\in G\C^{n\times n}$.

Associated with the universal covering surface $\tS$,
there is a group
$\Delta$ of deck transformations. A deck transformation
is an analytic bijection $\sigma:\tS\to\tS$ which satisfies
$\rho\circ\sigma=\rho$.

Given a solution $Y(\tz)$ of (\ref{f2.4}) and $\sigma\in\Delta$,
then $Y(\sigma(\tz))$ is also a solution of (\ref{f2.4}).
Hence there exists a unique matrix
$\chi(\sigma)\in G\C^{n\times n}$ such that
\bqn\label{f2.5}
Y(\tz) &=& Y(\sigma(\tz))\chi(\sigma).
\eqn
In fact, the mapping
$\chi:\Delta\to G\C^{n\times n}$ is a representation
(group homomorphism), i.e.,
$\chi(\sigma\tau)=\chi(\sigma)\chi(\tau)$.
This mapping is called the {\em monodromy representation}
of $Y(\tz)$.

Now let $\wh{Y}(\tz)$ be another solution of (\ref{f2.4}),
and let $\wh{\chi}:\Delta\to G\C^{n\times n}$ be the monodromy
representation of $\wh{Y}(\tz)$. If $Y$ and $\wh{Y}$ are
related by $\wh{Y}(\tz)=Y(\tz)C$ with
$C\in G\C^{n\times n}$, then
\bqn\label{f2.6}
\wh{\chi}(\sigma) &=& C\iv\chi(\sigma)C
\eqn
for all $\sigma\in\Delta$. This means that to the system
(\ref{f2.4}) there corresponds a class of mutually conjugate
representations. This class is called the
{\em monodromy} of the system (\ref{f2.4}).

It is interesting to note that the matrix function
$A(z)$ can be reconstructed from a solution $Y(\tz)$.
Assume that $Y:\tS\to G\C^{n\times n}$ is an analytic function
which satisfies (\ref{f2.5}) with some representation
$\chi:\Delta\to G\C^{n\times n}$. Then one can define
\bqn\label{f2.5x}
A(z) &:=&
Y'(\tz)Y(\tz)\iv,\qquad \tz\in\tS.
\eqn
This definition is correct since the right hand side turns out to
be a single-valued analytic function (i.e., it is invariant under
any deck transformation).

We will occasionally have to modify a solution of
a system
$Y'(\tz)=A(z)Y(\tz)$ with a (single-valued) analytic
matrix function
$V:S\to G\C^{n\times n}$ as follows: $\wh{Y}(\tz)=V(z)Y(\tz)$.
Then formula (\ref{f2.5x}) tells us that $\wh{Y}(\tz)$ is a
solution of the system $\wh{Y}'(\tz)=\wh{A}(z)\wh{Y}(\tz)$ with
\bqn\label{f2.6x}
\wh{A}(z) &=& V(z)A(z)V\iv(z)+V'(z)V\iv(z),\qquad z\in S.
\eqn
Moreover, the modified system has the same monodromy as the
original one. Conversely, suppose that $A(z)$ and $\wh{A}(z)$
are related to each other by (\ref{f2.6x}) with some
single-valued matrix function $V(z)$.
Then the solutions of the corresponding systems are connected
with each other by $\wh{Y}(\tz)=V(z)Y(\tz)C$ with
$C\in G\C^{n\times n}$.

Because the group $\Delta$ is finitely generated,
by elements $\sigma_1,\dots,\sigma_m$ say,
a complete characterization of the monodromy is given by
the matrices $\chi(\sigma_1),\dots,\chi(\sigma_m)$.
Although there is no distinguished choice for the generators,
we will choose them in a specific way.

Assume that $a_1,\dots,a_m\in\C$ lie in this order on a smooth
closed positively oriented curve $\Ga\subset\C$. Let $\tz_0\in\tS$
be a point such that $\rho(\tz_0)$ lies inside of $\Ga$. For
$k=1,\dots,m$, determine the (unique) point $\tz_k\in\tS$
satisfying $\rho(\tz_k)=\rho(\tz_0)$ and the following condition:
there exists a simple path $\wt{\gamma}_k:[0,1]\to\tS$ with
$\wt{\gamma}_k(0)=\tz_0$ and $\wt{\gamma}_k(1)=\tz_k$ such that
the curve $\rho(\wt{\gamma}_k)$ intersects $\Ga$ exactly twice,
first on the subarc joining $a_{k-1}$ and $a_k$, and second on the
subarc joining $a_{k}$ and $a_{k+1}$. (That means, the curve
$\rho(\wt{\gamma}_k)$ goes once around $a_k$ in positive
direction.) Having determined $\tz_k$, there exists a unique deck
transformation $\sigma_k\in\Delta$ for which
$\sigma_k(\tz_0)=\tz_k$. These deck transformations
$\sigma_1,\dots,\sigma_m\in\Delta$ satisfy the condition
\bqn\label{f2.8} \sigma_1\circ\sigma_2\circ\dots\circ\sigma_m &=&
\id. \eqn The construction of $\sigma_1,\dots,\sigma_m$ depends
primarily on the points $a_1,\dots,a_m$ and the curve $\Ga$. It
also depends on $\tz_0$, but for a different point $\tz_0$, we get
deck transformations $\sigma_k'=\tau\iv\circ\sigma_k\circ\tau$
with some $\tau\in\Delta$.

The monodromy representation of a solution of
(\ref{f2.4}) can now be described by the
$m$-tuple $[\chi(\sigma_1),\dots,\chi(\sigma_m)]$,
which (because of (\ref{f2.8})) satisfies the relation
\bqn
\chi(\sigma_1)\chi(\sigma_2)\cdots\chi(\sigma_m) &=& I.
\eqn
As we are interested in monodromy, let us introduce the following
equivalence relation:
\bqn\label{f2.equi}
[M_1,\dots,M_m] &\sim& [M_1',\dots,M_m']
\eqn
if and only if there exists a $C\in G\C^{n\times n}$ such that
$M_k'=C\iv M_k C$ for all $k=1,\dots,m$. The equivalence classes
are denoted by $[M_1,\dots,M_m]_{\sim}$ and will be called
the {\em simultaneous similarity class}.
Hence the monodromy of the system (\ref{f2.4}) is given by
class $[\chi(\sigma_1),\dots,\chi(\sigma_m)]_{\sim}$.

In order to study the behavior of the solutions of the system
(\ref{f2.4}) in a neighborhood of the singularities, we
introduce the punctured neighborhoods
$D_a=\{z\in\C:0<|z-a|<\ep\}$ for $a\in\C$ and
$D_\iy=\{z\in\C:|z|>1/\ep\}$ for $a=\iy$, where $\ep$ is
sufficiently small. Let $\tD_a$ be the universal
covering surface of $D_a$. The associated group of
deck transformations is singly generated by a transformation
$\sigma_a$, which sends a point $\tz\in\tD_a$ into its correspondent
after going once around the point $a$ in positive direction.

Let us restrict the solutions of the system (\ref{f2.4})
onto $\tD_{a_k}$ and assume for simplicity that $a_k\neq\iy$.
Then the function $Y:\tD_{a_k}\to G\C^{n\times n}$
has the property $Y(\tz)=Y(\sigma_{a_k}(\tz))M_k$, where
$M_k$ is similar to $\chi(\sigma_k)$.
Let $E_k$ be {\em any} matrix satisfying
$M_k = \exp(-2\pi iE_k)$, and introduce
\bqn
(\tz-a_k)^{E_k} &:=& \exp(E_k\ln(\tz-a_k)),\qquad
\tz\in\tD_{a_k},
\eqn
noting that the logarithm is well defined on $\tD_{a_k}$.
Then $Y(\tz)$ can be written as
\bqn\label{f2.14}
Y(\tz) &=& Z_{k}(z)(\tz-{a_k})^{E_k},\qquad
\tz\in\tD_{a_k},
\eqn
where $Z_k:D_{a_k}\to G\C^{n\times n}$ is a (single-valued)
analytic function, which may have (along with
its inverse) an isolated singularity at $a_k$.
The description of the local behavior
of the solution is now reduced to
characterizing the matrix $E_k$ and
the behavior of $Z_k(z)$ and its inverse.

When speaking of a restriction of a solution, which is defined
on $\tS$, onto $\tD_{a_k}$, we have to be a little bit careful.
Namely, if $\rho:\tS\to S$ is the covering map, then
the preimage $\rho\iv(D_{a_k})$ can be identified
with $\tD_{a_k}\times \Omega$ where $\Omega$ is an index set,
which is countably infinite for $m\ge 3$.
For a solution $Y:\tS\to G\C^{n\times n}$, the
relation (\ref{f2.14}) should therefore correctly be written as
\bqn\label{f2.15}
Y(\tz_*) &=& Z_k(z)(\tz-a_k)^{E_k}C_{\omega},\qquad
\tz_*=(\tz,\omega)\in\tD_{a_k}\times \Omega,
\eqn
where $C_\omega\in G\C^{n\times n}$ depends on the connected
component $\Omega$ of the preimage $\rho\iv(D_{a_k})$.

A singularity $a_k$ of a system is called {\em apparent}
if the solutions have no branching behavior
at $a_k$, i.e., if $\chi(\sigma_k)=I$.
In this case it is not necessary to take this singularity
into account when defining the universal covering surface.

The following fundamental theorem says that there exist systems
with prescribed singularities and prescribed monodromy.
A first proof was given by Plemelj \cite{Pl}, who employed
the theory of singular integral equations.
Another proof \cite{Ro,AnBo} can be given by using the
theory of vector bundles \cite{Fo,St}.

\begin{theorem}\label{t.mono}
For any distinct points $a_1,\dots,a_m\in\Co$ and
matrices $M_1,\dots,M_m\in\C^{n\times n}$ satisfying
$M_1\cdots M_m=I$, there exists a system
with singularities only at $a_1,\dots,a_m$ and with
monodromy given by $[M_1,\dots,M_m]_\sim$.
\end{theorem}

This theorem (as it is stated) does not say anything
about the type of the singularities.
A singularity $a_k$ of a system is called {\em Fuchsian}
if $A(z)$ has only a simple pole at $z=a_k$.
A singularity $a_k$ is called {\em regular}
if the function $Z_k(z)$ in (\ref{f2.15})
has not an essential singularity at $z=a_k$.
A Fuchsian singularity is always regular.

The result of Plemelj is actually stronger than the statement of the previous
theorem. One can construct a system which has only Fuchsian singularities with
the possible exception of one singularity, which is (in general) merely
regular.

The question whether there exists a system with prescribed
singularities and monodromy such that all singularities are
Fuchsian is known as the Riemann--Hilbert problem or
Hilbert's 21st problem. Surprisingly, Bolibruch
\cite{Bo1,AnBo,Bo4}
showed that the answer may be negative.
Notice, however, that Plemelj's results still implies that
the answer is always positive if one allows one additional
apparent singularity.

A generalized version of the Riemann--Hilbert problem, where
one asks for systems with prescribed Fuchsian singularities,
prescribed monodromy and prescribed local behavior
(in a certain specified sense)
was also studied by Bolibruch \cite{Bo3}.

The Riemann--Hilbert problem and its modifications have an
intimate relation to the problem which will be proposed in
Section \ref{s3}. In fact, all these problems can be
formulated by using the language of vector bundles.
We will not elaborate on this direction, but merely confine
ourselves to some comments later on.

Finally, let us recall some facts about (scalar) linear
differential equations of higher order \cite{Ha,In,CL}.
A linear differential equation of $n$-th order
can be written as
\bqn\label{f2.16}
y^{(n)}+q_1(z)y^{(n-1)}+\dots+q_n(z)y&=&0,
\eqn
where $y^{(k)}$ stands for the $k$-th derivative of the
complex valued function $y$.
We assume that $q_1(z),\dots,q(z)$ are analytic functions
on $\C$ with the exception of a finite number of points, which
represent the {\em singularities} of the equation.
The point $a=\iy$ is a singularity if
so is the point $z=0$ for the equation resulting from a
change of variables $z\mapsto1/(z-a)$. (This definition
does not depend on the choice of $a\in\C$.)

The solutions of the above equation with singularities
$a_1,\dots,a_m\in\Co$ are analytic functions
$y:\tS\to\C$ defined on the corresponding universal
covering surface. The set of all solutions
forms an $n$-dimensional linear space.
Given $n$ linearly independent solutions
$y_1(\tz),\dots,y_n(\tz)$, the
monodromy representation $\chi:\Delta\to G\C^{n\times n}$
is defined by
\bqn
(y_1(\tz),\dots,y_n(\tz)) &=&
(y_1(\sigma(\tz)),\dots,y_n(\sigma(\tz)))\chi(\sigma),
\eqn
where $\chi(\sigma)\in G\C^{n\times n}$ is a uniquely
determined matrix and $\sigma\in\Delta$.

A singularity $a\in\C$ of equation (\ref{f2.16}) is called {\em Fuchsian}
if one can write $q_k(z)=r_k(z)/(z-a)^k$ where $r_k(z)$ is
analytic at $z=a$ for each $k=1,\dots,n$.
The {\em local exponents} $\rho_1,\dots,\rho_n$
of this singularity are by definition the roots of the
{\em indical equation}
\bqn\label{f2.ind}
\rho(\rho-1)\cdots(\rho-n+1)+
\sum_{k=1}^n r_k(a)\rho(\rho-1)\cdots(\rho-n+1+k)
&=&0.
\eqn
The singularity at infinity is Fuchsian if
$q_k(z)=O(1/z^k)$ as $z\to\iy$ for each $k=1,\dots,n$,
and the local exponents are defined similarly.

Now suppose that equation (\ref{f2.16}) has exactly $m$
Fuchsian singularities $a_1,\dots,a_m$
with local exponents
$\rho_1^{(k)},\dots,\rho_n^{(k)}$ for $a_k$.
Then the well known Fuchs' relation says that
\bqn\label{f2.fuchs}
\sum_{k=1}^m (\rho_1^{(k)}+\dots+\rho_n^{(k)})
&=& (m-2)\frac{n(n-1)}{2}.
\eqn
A couple of more facts about linear differential equations
will be stated in Section \ref{s5} when they are needed.

\reseteqn
\section{Systems of standard form and
their relation to the factorization problem}
\label{s3}

In this section, we introduce a class of systems
which are characterized by their singularities, their
monodromy and their local behavior near the singularities.
In addition, we allow an apparent singularity at infinity
the behavior of which is described by certain integers.

We prove that such systems exist for given singularities,
monodromy and local data. We also discuss the ``uniqueness'',
i.e., we examine the relation between two systems having the same
singularities and the same data. Finally, we show that the factors appearing
in the factorization of a piecewise constant matrix function can be
constructed by using the solutions of such systems.
The problem of determining the partial indices reduces to a
corresponding problem for systems of linear differential equations.

Let $a_1,\dots,a_m\in\C$ be distinct points which lie
in this order on a smooth closed positively oriented curve
$\Ga\subset\C$. Denote by $\tS$ the universal covering
surface of $S:=\Co\setminus\{a_1,\dots,a_m\}$.
Define deck transformations
$\sigma_1,\dots,\sigma_m$ for $\tS$ in the way
explained in the previous section.

A matrix $E$ is called {\em non-resonant} if
the difference of any two eigenvalues of $E$
is not a nonzero integer.
A collection $[M_1,\dots,M_m]$ and $[E_1,\dots,E_m]$
of $m$-tuples of $n\times n$ matrices
will be called {\em admissible data} if the following
conditions are satisfied:
\begin{itemize}
\item[(a)]
$M_1M_2\cdots M_m=I$.
\item[(b)]
$M_k \sim \exp(-2\pi iE_k)$ for each $k=1,\dots,m$.
\item[(c)]
The matrices $E_1,\dots,E_m$ are non-resonant.
\end{itemize}

\begin{definition}\label{d3.1}
We say that the system
$Y'(\tz) = A(z) Y(\tz)$
is of standard form with respect to admissible data
$[M_1,\dots,M_m]$ and $[E_1,\dots,E_m]$, with
singularities $a_1,\dots,a_m\in\C$ and with indices
$\ka_1,\dots,\ka_n\in\Z$, where $\ka_1\ge\ka_2\ge\dots\ge\ka_n$,
if the following conditions are satisfied:
\begin{itemize}
\item[(i)]
The system has singularities only at $a_1,\dots,a_m$ and at
infinity, where the singularity at infinity is apparent.
\item[(ii)]
The monodromy of the system is given by
$[\chi(\sigma_1),\dots,\chi(\sigma_m)]_\sim =
[M_1,\dots,M_m]_\sim$.
\item[(iii)]
For each $k=1,\dots,m$, the solutions of the system
defined on $\tD_{a_k}$ can be written as
$$ Y(\tz)\quad=\quad
Z_k(z)(\tz-a_k)^{E_k}C,
\qquad\quad\tz\in\tD_{a_k},$$
where the function $Z_k:D_{a_k}\cup\{a_k\}\to G\C^{n\times n}$
is analytic and $C\in G\C^{n\times n}$.
\item[(iv)]
The solutions of the system defined on $D_\iy$
can be written as
$$ Y(z)\quad=\quad
\diag(z^{\ka_1},\dots,z^{\ka_n})Z_\iy(z)C,
\qquad\quad z\in D_{\iy},$$
where the function $Z_\iy:D_{\iy}\cup\{\iy\}\to G\C^{n\times n}$
is analytic and $C\in G\C^{n\times n}$.
\end{itemize}
\end{definition}

The necessity of the assumption (a) and (b) on the data
is obvious. Assumption (c) is not necessary for this definition,
but it will be important later on.
Note that the solutions $Y(\tz)$ are
analytic and invertible functions on
$\tS^*=\tS\setminus\rho\iv(\iy)$.

A simple observation is in order.
The data $[M_1,\dots,M_m]$ and $[E_1,\dots,E_m]$
of a given system of standard form can be replaced by data
$[CM_1C\iv,\dots,CM_mC\iv]$ and
$[C_1E_1C_1\iv,\dots,C_mE_mC_m\iv]$, where
$C,C_1,\dots,C_m\in G\C^{n\times n}$.
Conversely, to a given system of standard form one can
uniquely associate the monodromy data
$[M_1,\dots,M_m]_\sim$ and the local data
$[(E_1)_\sim,\dots,(E_m)_\sim]$.
Indeed, in order to see that the matrices $E_k$ are unique up
to similarity, one can use (iii) and consider the residue of
$A(z)=Y'(\tz)Y\iv(\tz)$ at $z=a_k$.

Finally, also the indices of a system of standard form
are unique. Here one can use (iv) and consider the residue
of $Y\iv(z)Y'(z)$ of any solution defined on $D_\iy$.

Before discussing the existence and uniqueness of
systems of standard form for given singularities
and data, we will illustrate
how these systems provide us with the solution of the
original factorization problem for piecewise constant
matrix functions.

In the following theorem we assume that the curve $\Ga$
divides the complex sphere $\Co$ into two domains
$D_+$ and $D_-$ such that $0\in D_+$ and $\iy\in D_-$.
As usual, we identify functions which are analytic on
the domains $D_+$ or $D_-$ with their boundary values,
considered as functions on $\Ga$.

\begin{theorem}\label{t3.2}
Let $G\in PC(\Ga)^{n\times n}$ be a piecewise constant function
with jumps only at $a_1,\dots,a_m$. Assume that $G$ admits a
$\Phi$-factorization in the space $L^p(\Ga)$, $1<p<\iy$, and let
$[M_1,\dots,M_m]$ and $[E_1,\dots,E_m]$ be the data
associated with $G$ with respect to $L^p(\Ga)$.

Assume that there exists a system of standard form with
these data, with singularities $a_1,\dots,a_m$ and
with indices $\ka_1,\dots,\ka_n$. Let $Y_1$ and $Y_2$
be solutions of this system defined on $D_+$ and
$D_-\setminus\{\iy\}$, respectively. Then there exist
$C_1,C_2\in G\C^{n\times n}$ such that the representation
\bqn
G(t) &=& G_+(t)\La(t)G_-(t),
\qquad\quad t\in\Ga,\nn
\eqn
is a $\Phi$-factorization of $G$ in the space $L^p(\Ga)$,
where
$\La(t)=\diag(t^{\ka_1},\dots,t^{\ka_n})$ and
$$
G_+(z)\;\;=\;\;C_1\iv Y_1\iv(z),\quad z\in D_+,\qquad
G_-(z)\;\;=\;\;\La\iv(z)Y_2(z)C_2,\quad
z\in D_-\setminus\{\iy\}.
$$
\end{theorem}

\begin{proof}
Note first that $Y_1$ and $Y_2$ exist because
$D_+$ and $D_-$ are simply connected and infinity is an
apparent singularity of the system.
Let $\wt{G}_+(z)=Y_1\iv(z)$, $z\in D_+$, and
let $\wt{G}_-(z)=\La\iv(z)Y_2(z)$, $z\in D_-\setminus\{\iy\}$.
{}From the condition (iv) of Definition \ref{d3.1}, it follows
that $\wt{G}_-$ is also analytic and invertible at $\iy$.
Hence $G_+$ and $G_-$ are analytic and invertible on all of
$D_+$ and $D_-$, respectively.
It follows from (iii)
that in a neighborhood of the singularity $a_k$,
the entries of the functions $\wt{G}_+\iv$ and $\wt{G}_-$
behave as $(z-a_k)^{\lambda}(\ln(z-a_k))^{\sigma}$
and the entries of the functions $\wt{G}_+$ and $\wt{G}_-\iv$
behave as $(z-a_k)^{-\lambda}(\ln(z-a_k))^{\sigma}$,
where $\lambda$ is taken from the eigenvalues of $E_k$ and
$\sigma\in\{0,1,\dots\}$ refers to their multiplicity.
By assumption we have $-1/q=1/p-1<\Re\lambda<1/p$.
This implies that
$\wt{G}_+\in L^p_+(\Ga)^{n\times n}$,
$\wt{G}_-\iv\in L^q_-(\Ga)^{n\times n}$,
$\wt{G}_+\iv\in L^p_+(\Ga)^{n\times n}$ and
$\wt{G}_-\in L^q_-(\Ga)^{n\times n}$.

The functions $Y_1$ and $Y_2$ can be extended analytically
onto the boundary $\Ga\setminus\{a_1,\dots,a_m\}$.
The values of the functions defined by $\wt{G}_+(t)=Y_1\iv(t)$
and $\wt{G}_-(t)=\La\iv(t)Y_2(t)$ with
$t\in\Ga\setminus\{a_1,\dots,a_m\}$ coincide with
the boundary values of the above functions $\wt{G}_+$
and $\wt{G}_-$, which are defined on $D_+$ and $D_-$.

We introduce the function
$\wt{G}(t)=\wt{G}_+(t)\La(t)\wt{G}_-(t)$, $t\in\Ga$,
and observe that this representation is a factorization of
$\wt{G}$ in the space $L^p(\Ga)$. On the other hand, we have
$\wt{G}(t)=Y_1\iv(t)Y_2(t)$, and because
$Y_1$ and $Y_2$ are solutions of the system, it follows
that $\wt{G}$ is a piecewise constant function on $\Ga$
with jumps only at $a_1,\dots,a_m$.

Now let us express the monodromy of the system in terms of
$\wt{G}$. We start from the solution $Y_1(z)$ on $D_+$
and continue analytically by going once around the point
$a_k$ (see the construction of the deck
transformations $\sigma_k$). When we first cross $\Ga$, the
proper analytic continuation is the function
$Y_2(z)\wt{G}(a_k-0)\iv$ on $D_-$. After the second crossing of
$\Ga$, we obtain the function $Y_1(z)\wt{G}(a_k+0)
\wt{G}(a_k-0)\iv$ on $D_+$. Hence, up to simultaneous similarity,
we have $\chi(\sigma_k)=\wt{G}(a_k-0)\wt{G}(a_k+0)\iv$.
On the other hand, by condition (ii), the monodromy is
given by $[M_1,\dots,M_m]$. Hence $\chi(\sigma_k)=
C_1M_kC_1\iv=C_1G(a_k-0)G(a_k+0)\iv C_1\iv$
with some $C_1\in G\C^{n\times n}$.
This implies that
$\wt{G}(a_k-0)\wt{G}(a_k+0)\iv=C_1G(a_k-0) G(a_k+0)\iv C_1\iv$,
hence
$\wt{G}(a_k+0)\iv C_1G(a_k+0)=\wt{G}(a_k-0)\iv C_1G(a_k-0)$
for each $k=1,\dots,m$. Because both $G$ and $\wt{G}$ are
piecewise constant, there is a $C_2\in G\C^{n\times n}$
such that $C_2=\wt{G}(t)\iv C_1G(t)$.
Hence $G(t)=C_1\wt{G}(t)C_2$, which implies that
$G(t)=G_+(t)\La(t)G_-(t)$ is a factorization in the space
$L^p(\Ga)$.
By Corollary \ref{c1.2a}, this is even a
$\Phi$-factorization.
\end{proof}

Note that the data associated with a piecewise constant matrix
function is necessarily admissible data.
Actually, assumption (c) is even weaker than the condition
that the real parts of the eigenvalues of the matrices
$E_k$ lie in the open interval $J_p$.

We need the following factorization-type result about
analytic matrix functions defined on a punctured
disk $D_a$. It is the analytic version of the
Birkhoff-Grothendieck theorem \cite{Bi}.
A quite simple proof has also been given by
Leiterer \cite{Lei} (see also \cite{AnBo}).

\begin{lemma}\label{l3.3}
For $a\in\C$, let $\Phi:D_a\to G\C^{n\times n}$ be an analytic
function. Then there exist $\ka_1,\dots,\ka_n\in\Z$
such that $\Phi$ can be written in the form
\bqn
\Phi(z) &=& U(z)\diag((z-a)^{\ka_1},\dots,(z-a)^{\ka_n})
V(z),\qquad z\in D_a,\nn
\eqn
where $U:\Co\setminus\{a\}\to G\C^{n\times n}$ and
$V:D_a\cup\{a\}\to G\C^{n\times n}$ are analytic functions.
\end{lemma}

Using a well known modification idea \cite{AnBo},
we prove the existence of systems of standard form
with prescribed data and singularities, which
have certain (a priori not known) indices.

\begin{theorem}\label{t3.3}
Let the $n\times n$ matrices
$[M_1,\dots,M_m]$ and $[E_1,\dots,E_m]$ be
admissible data, and let $a_1,\dots,a_m\in\C$
be distinct points.
Then there exists a system of standard form
with these data and singularities and with certain indices
$\ka_1,\dots,\ka_n\in\Z$.
\end{theorem}
\begin{proof}
{}From Theorem \ref{t.mono} it follows that there exists a system with
singularities only at $a_1,\dots,a_m$ whose monodromy is
given by $[M_1,\dots,M_m]_{\sim}$. Let $Y_0$ be a solution of
this system, i.e., $Y_0:\tS^*\to G\C^{n\times n}$ is
analytic and has the given monodromy representation.

We are going to modify $Y_0$ step by step in order to obtain an
analytic matrix function $Y:\tS^*\to G\C^{n\times n}$ which also
satisfies (iii) and (iv) of Definition \ref{d3.1}. Let us consider
$Y_0$ near the point $a_1$. Because the monodromy is given as
above, the solution changes by a matrix similar to $M_1$ after
going once around $a_1$ in positive direction. Hence, restricted
to $\tD_{a_1}$, the function $Y_0$ can be written in the form
$Y_0(\tz)=\Phi_1(z)(\tz-a_1)^{E_1}C$, where $\Phi_1:D_{a_1}\to
G\C^{n\times n}$ is a certain analytic function. We apply Lemma
\ref{l3.3} with $\Phi=\Phi_1$ and $a=a_1$. We obtain that
$\Phi_1(z)=\wt{U}_1(z)V_1(z)$, $z\in D_{a_1}$, where $\wt{U}_1$ is
analytic and invertible on $\C\setminus\{a_1\}$ and $V_1$ is
analytic and invertible on $D_{a_1}\cup\{a_1\}$. Hence
$\wt{U}_1\iv(z) Y_0(\tz) = V_1(z)(\tz-a_1)^{E_1}C$ for
$\tz\in\tD_{a_1}$. We define the function
$Y_1(\tz):=\wt{U}_1\iv(z) Y_0(\tz)$ for $\tz\in\tS^*$. This
function satisfies the conditions (i), (ii), and also the
condition (iii) at the point $a_1$.

The next steps consist in doing the same for $a_2,\dots,a_m$. The
essential point is that, at the $k$-step, we define the function
$Y_k(\tz)=\wt{U}_k\iv(z)Y_{k-1}(\tz)$, where $\wt{U}_k(z)$ is
analytic and invertible on $\C\setminus\{a_k\}$. This modification
does not produce additional singularities (except at infinity), it
does not change the monodromy, and it changes the local behavior
only at the $a_k$ (in the desired way), but not at the other
singularities. So we obtain a function $Y_m(\tz)$ which satisfies
the conditions (i)--(iii).

Finally, we consider the point infinity. Restricting $Y_m(\tz)$ to
$D_{\iy}$ we can write $Y_m(z)=\Phi(1/z)C$, where $\Phi(1/z)$ is
analytic and invertible on $D_{0}$. We apply Lemma \ref{l3.3} with
$a=0$ and factor $\Phi(z)=V(z)\La(z)U(z)$, $z\in D_{0}$, where
$V(z)$ is analytic and invertible on $\Co\setminus\{0\}$, $\La(z)$
is the diagonal matrix, and $U(z)$ is analytic and invertible on
$D_0\cup\{0\}$. It follows that $V\iv(1/z)
Y_m(z)=\La(1/z)U(1/z)C$. We define $Y(\tz)=V\iv(1/z) Y_m(\tz)$,
$\tz\in\tS^*$. Because $V(1/z)$ is analytic and invertible on
$\C$, the function $Y(\tz)$ satisfies also the required properties
(i)--(ii). Since $Y(z)=\La(1/z)U(1/z)C$ on $D_{\iy}$, this
function also satisfies condition (iv) after an appropriate
permutation of the entries of $\La(1/z)$.

Finally, we define the corresponding matrix function
$A(z)$ by $A(z)=Y'(\tz)Y\iv(\tz)$.
\end{proof}

Next we consider the question about the uniqueness of systems of
standard form with prescribed data and singularities. First of
all, the indices are uniquely determined. The systems themselves
are not unique, but they are related to each other in a simple
way. Because of the relation to factorization (described in
Theorem \ref{t3.2}) it should not be surprising that the
description is analogous to that given in Theorem \ref{t1.2}.

\begin{theorem}\label{t3.5}
Let $Y_1'(\tz)=A_1(z)Y_1(\tz)$ and $Y_2'(\tz)=A_2(z)Y_2(\tz)$
be two systems of standard form with the same
admissible data
and singularities $a_1,\dots,a_m\in\C$.
Then the indices corresponding to both systems coincide.
Moreover, if the indices are given by
$$\ovl{\ka}_1,\dots,\ovl{\ka}_1,\ovl{\ka}_2,\dots\ovl{\ka}_2,
\ovl{\ka}_3,\dots,\ovl{\ka}_{k-1},
\ovl{\ka}_k,\dots,\ovl{\ka}_k,$$
where $\ovl{\ka}_j$ occurs exactly $l_j$ times,
$\ovl{\ka}_1>\ovl{\ka}_2>\dots>\ovl{\ka}_k$,
$l_1,\dots,l_k\in\{1,2,\dots\}$, $l_1+\dots+l_k=n$,
then both systems are related to each other by
\bqn\label{f.equ1}
A_2(z) &=& V'(z)V\iv(z)+V(z)A_1(z)V\iv(z),\qquad
z\in S\setminus\{\iy\},
\eqn
or, equivalently, by
\bqn\label{f.equ2}
Y_2(\tz) &=& V(z)Y_1(\tz)C,
\qquad \tz\in\tS^*,
\eqn
where $C\in G\C^{n\times n}$ and
$V:\C\to  G\C^{n\times n}$ is a rational
matrix function of the form
\bqn\label{f.Vmat}
V(z) &=& \left(\ba{cccc}
V_{11} & V_{12}(z) & \dots & V_{1k}(z) \\
0 & V_{22} & & V_{2k}(z) \\
\vdots & & \ddots & \vdots \\
0 & 0 & \dots & V_{kk} \ea\right)
\eqn
with $V_{jj}\in G\C^{l_j\times l_j}$,
$V_{ij}(z) = \sum_{r=0}^{\ovl{\ka}_i-\ovl{\ka}_j}
V_{ij}^{(r)} z^r$,
where $V_{ij}^{(r)}\in\C^{l_i\times l_j}$
for $i<j$.
\end{theorem}
\begin{proof}
Because the systems
have the same singularities and the same monodromy, there
exist an analytic function $V:S\setminus\{\iy\}\to
G\C^{n\times n}$ and a $C\in G\C^{n\times n}$ such that
$Y_2(\tz)=V(z)Y_1(\tz)C$. We are first going to show
that $V(z)$ is analytic and invertible on all of $\C$.

Let us restrict the solutions $Y_1$ and $Y_2$ onto
a neighborhood $\tD_{a_k}$ of the point $a_k$.
{}From (iii) of Definition \ref{d3.1} it follows
that there are $C_1,C_2\in G\C^{n\times n}$ such that
$$
Y_1(\tz)\;\;=\;\; Z_k^{(1)}(z)(\tz-a)^{E_k}C_1,\qquad
Y_2(\tz)\;\;=\;\; Z_k^{(2)}(z)(\tz-a)^{E_k}C_2,\qquad
\tz\in \tD_{a_k},
$$
where $Z_k ^{(1)}$ and $Z_k^{(2)}$ are analytic and invertible
at $D_{a_k}\cup\{a_k\}$. Hence
\bqn\label{f3.x}
Z_k^{(2)}(z)(\tz-a)^{E_k}\wt{C} &=&
V(z)Z_k^{(1)}(z)(\tz-a)^{E_k},
\eqn
where $\wt{C}=C_2C\iv C_1\iv$.
By going once around the point $a_k$ in positive
direction it follows that
$$
Z_k^{(2)}(z)(\tz-a)^{E_k}N_k \wt{C} \quad=\quad
V(z)Z_k^{(1)}(z)(\tz-a)^{E_k}N_k,
$$
where $N_k:=\exp(2\pi iE_k)$.
Combining both equations it follows that
$N_k=\wt{C}\iv N_k\wt{C}$, i.e.,
$N_k$ and $\wt{C}$ commute with each other.
Now we use the fact that the matrix $E_k$ is non-resonant.
Because $N_k=\exp(2\pi iE_k)$, this implies that
$E_k$ can be written as a polynomial in the matrix $M_k$.
Hence also $E_k$ and $\wt{C}$ commute.
Using this in connection with (\ref{f3.x}) we obtain
\bqn
Z_k^{(2)}(z)\wt{C}(\tz-a)^{E_k} &=&
V(z)Z_k^{(1)}(z)(\tz-a)^{E_k}.\nn
\eqn
Hence $V(z)=Z_k^{(2)}(z)\wt{C}(Z_k^{(1)}(z))\iv$, which implies
that $V(z)$ is analytic and invertible at $a_k$.

So far we have shown that $V(z)$ and $V\iv(z)$ are entire analytic
functions. We analyze their behavior at infinity by restricting
the  solutions $Y_1$ and $Y_2$ onto $D_\iy$. From (iv) of
Definition \ref{d3.1}, it follows that $$
Y_1(z)\;\;=\;\;\La_1(z)Z_{\iy}^{(1)}C_1,\qquad
Y_2(z)\;\;=\;\;\La_2(z)Z_{\iy}^{(2)}C_2,\qquad z\in D_\iy, $$
where $C_1,C_2\in G\C^{n\times n}$ and $Z_{\iy}^{(1)}$ and
$Z_{\iy}^{(2)}$ are analytic and invertible on a neighborhood of
infinity. Here we assume that $\La_1$ and $\La_2$ are of the form
$$ \La_s(z)\;\;=\;\; \diag(z^{\ovl{\ka}_1}I_{l_1^{(s)}},
z^{\ovl{\ka}_2}I_{l_2^{(s)}},\dots,
z^{\ovl{\ka}_k}I_{l_k^{(s)}}),\qquad s=1,2, $$ where
$\ovl{\ka}_1,\dots,\ovl{\ka}_k\in\Z$,
$\ovl{\ka}_1>\ovl{\ka}_2>\dots>\ovl{\ka}_k$,
$l_1^{(s)},\dots,l_k^{(s)}\in\{0,1,\dots\}$,
$l_1^{(s)}+\dots+l_k^{(s)}=n$. Here $I_{l}$ stands for the
identity matrix of size $l\times l$. Now the above equations give
$$ \La_2(z)Z_{\iy}^{(2)}C_2\quad=\quad
V(z)\La_1(z)Z_{\iy}^{(1)}C_1C,\qquad z\in D_\iy. $$ This means
that $$ V(z)\;\;=\;\; \La_2(z)Z_\iy(z)\La_1\iv(z),\qquad
V\iv(z)\;\;=\;\; \La_1(z)Z_\iy\iv(z)\La_2\iv(z), $$ where
$Z_\iy:=Z_\iy^{(2)}C_2C\iv C_1\iv (Z_\iy^{(1)})\iv$ is analytic
and invertible in a neighborhood of infinity. Now we introduce
block partitions of the matrices $V(z)$ and $V\iv(z)$: $$
V(z)\;\;=\;\;(V_{ij}(z))_{i,j=1}^{k},\qquad
V\iv(z)\;\;=\;\;(U_{ij}(z))_{i,j=1}^{k}, $$ where $V_{ij}(z)$ is
of size $l_i^{(2)}\times l_j^{(1)}$ and $U_{ij}(z)$ is of size
$l_i^{(1)}\times l_j^{(2)}$. Using the structure of $\La_1(z)$ and
$\La_2(z)$, it follows that
$\|V_{ij}(z)\|=O(z^{\ovl{\ka}_i-\ovl{\ka}_j})$ and
$\|U_{ij}(z)\|=O(z^{\ovl{\ka}_i-\ovl{\ka}_j})$ as $z\to\iy$.
Because $V_{ij}(z)$ and $U_{ij}(z)$ are entire analytic functions,
we obtain that $V_{ij}=0$ and $U_{ij}=0$ for $i>j$, and that
$V_{ij}(z)$ and $U_{ij}(z)$ are matrix polynomials of degree at
most $\ovl{\ka}_i-\ovl{\ka}_j$ if $i\le j$. In particular, the
diagonal blocks are constant matrices. The block triangular
structure of $V(z)$ and $V\iv(z)$ implies that
$V_{ii}U_{ii}=I_{l_i^{(2)}}$ and $U_{ii}V_{ii}=I_{l_i^{(1)}}$.
Hence $V_{ii}$ and $U_{ii}$ are invertible square matrices, i.e.,
$l_i^{(1)}=l_i^{(2)}$. Hence the indices corresponding to both
systems are the same. Thereby, we have also proved that $V(z)$ is
of the desired form.
\end{proof}

It is easy to see that also the converse of the above theorem
holds. This suggests the following definition of an equivalence
relation for systems of standard form. Two systems
$Y_1'(\tz)=A_1(z)Y_1(\tz)$ and $Y_2'(\tz)=A_2(z)Y_2(\tz)$ of
standard form will be called {\em equivalent} if their indices are
the same and if there exists a matrix $V(z)$ of the form
(\ref{f.Vmat}) and the properties stated there such that
(\ref{f.equ1}) holds (or, equivalently, (\ref{f.equ2}) holds).

The next corollary resumes the results about
existence a uniqueness in a very brief way.
In fact, the uniqueness results will be employed
implicitly in Section \ref{s6}.

\begin{corollary}\label{c3.6}
Let $a_1,\dots,a_m\in\C$ be distinct points.
Then Definition
\ref{d3.1} establishes a one-to-one correspondence
between the equivalence classes of admissible data,
$$
[M_1,\dots,M_m]_\sim\quad\mbox{ and }\quad
[(E_1)_\sim,\dots,(E_m)_\sim],
$$
and the equivalence classes of systems of standard form
with above singularities.
\end{corollary}

The previous considerations show that to each class of admissible
data and given singularities one can associate a class of systems
of standard form, and hence certain (uniquely determined) indices.
For the problem of determining the partial indices for a piecewise
constant matrix function, one has now to solve the following
question: what are the indices associated to given admissible data
and singularities? As we will see, this question can be answered
explicitly only in some cases. In general, it leads to monodromy
problems for which no explicit solution is known.

In order to indicate the importance of the non-resonance
assumption for the matrices $E_1,\dots,E_m$, we show that
Theorem \ref{t3.5} breaks down if this assumption is dropped.
Consider the following matrix functions, which are the
solutions of certain $2\times 2$ systems with
three singularities at $a_1,a_2,a_3$,
\bqn
Y_1(z) &=& \left(\ba{cc}
(z-a_1)(z-a_2)(z-a_3) & 0 \\ 0 & 1 \ea\right),
\\[.5ex]
Y_2(z) &=& \left(\ba{cc}
(z-a_1)(z-a_2) & 0 \\ 0 & (z-a_3) \ea\right),
\\[.5ex]
Y_3(z) &=& \left(\ba{cc} (z-a_1)(z-a_2) & 0 \\ 1 & (z-a_3)
\ea\right). \eqn The corresponding systems have trivial monodromy
($M_1=M_2=M_3=I$), the same local behavior (with
$E_1=E_2=E_3=\diag(1,0)$), but their indices are $(3,0)$ for $Y_1$
and $(2,1)$ for $Y_2$ and $Y_3$. Moreover, one can show (which
requires a little effort) that $Y_2$ and $Y_3$ are not equivalent
in the above sense, although they have the same indices.
Obviously, Corollary \ref{c3.6} also breaks down, and the question
``what are the indices for given data?'' does not make sense. One
might think of replacing this question by asking for all possible
indices for given data, but we will not pursue this direction in
this paper.

Finally, let us mention the connection to vector bundles
(we refer to \cite[Sect.~5.1]{AnBo} for details). For given data and
singularities one constructs in a certain way a vector bundle
over $\Co$. Due to the fact that the matrices $E_k$ are
non-resonant, this construction is essentially unique.
Any such vector bundle is characterized by certain integers
$\ka_1,\dots,\ka_n$, which are called the {\em splitting type}
of the vector bundle.
Using the construction of the vector bundle and the
splitting type, one can easily construct the corresponding
system of standard form, which has indices
$\ka_1,\dots,\ka_n$.
So the above question takes the following form: what is
the splitting type of the vector bundles constructed for given
data and singularities?

\reseteqn
\section{Equivalent characterizations of systems of standard form}

The definition of systems of standard form has been given
entirely in terms of the solutions $Y(\tz)$ of the system.
In this section, we show by using standard arguments that
the conditions (i), (iii) and (iv) of Definition \ref{d3.1}
(i.e., except the monodromy condition)
can be replaced by equivalent characterizations given entirely
in terms of $A(z)$.

We prepare with the following well known lemma.
Note that the assumption that $E$ is non-resonant is
essential for the validity of the lemma.

\begin{lemma}\label{l4.1}
Let $A(z)$ be an analytic matrix function on $D_a$, $a\in\C$,
and assume that $A(z)$ has only a simple pole at $z=a$
whose residue is a non-resonant matrix $E$,
\bqn
A(z) &=& \frac{E}{z-a}\;+\;
\mbox{``analytic term at $z=a$''},\qquad
z\in D_a.\nn
\eqn
Then the solutions $Y:\tD_a\to G\C^{n\times n}$
of $Y'(\tz)=A(z)Y(\tz)$ can be written as
\bqn
Y(\tz) &=& Z(z)(\tz-a)^{E}C,\qquad
\tz\in\tD_a,\nn
\eqn
where the function $Z:D_a\cup\{a\}\to G\C^{n\times n}$
is analytic and $C\in G\C^{n\times n}$.
\end{lemma}

\begin{proposition}\label{p4.2}
The system  $Y'(\tz)=A(z)Y(\tz)$ is of
standard form with (admissible) data
$[M_1,\dots,M_m]$ and $[E_1,\dots,E_m]$,
with singularities $a_1,\dots,a_m\in\C$,
and with indices $\ka_1,\dots,\ka_n\in\Z$
if and only if the following conditions are satisfied:
\begin{itemize}
\item[(i)]
The matrix function $A(z)$ is analytic on
$\C\setminus\{a_1,\dots,a_m\}$.
\item[(ii)]
The monodromy of the system is given by
$[\chi(\sigma_1),\dots,\chi(\sigma_m)]_\sim
=[M_1,\dots,M_m]_\sim$.
\item[(iii)]
For each $k=1,\dots,m$, the function $A(z)$ has at most a simple
pole at $z=a_k$ whose residue is a matrix similar to $E_k$:
\bqn
A(z) &=& \frac{S_kE_kS_k\iv}{z-a_k}\;+\;
\mbox{``analytic term at $z=a_k$''},\qquad
z\in D_{a_k}.\nn
\eqn
\item[(iv)]
At infinity, the function $A(z)$ is of the form
\bqn
A(z) &=& \frac{\diag(\ka_1,\dots,\ka_n)}{z}
\;+\;
\La(z)R(z)\La\iv(z),
\qquad z\in D_\iy,\nn
\eqn
where $\La(z)=\diag(z^{\ka_1},\dots,z^{\ka_n})$ and
$R:D_\iy\to\C^{n\times n}$ is an analytic function satisfying
$\|R(z)\|=O(1/z^2)$ as $z\to\iy$.
\end{itemize}
\end{proposition}
\begin{proof}
Assume (i)--(iv) of Definition \ref{d3.1} holds.
Then also conditions (i) and (ii) of this proposition
are fulfilled.
As to (iii), remark that
$Y(\tz)=Z_k(z)(\tz-a_k)^{E_k}C$ implies
\bqn
A(z) &=& Y'(\tz)Y\iv(\tz)
\quad=\quad
Z_k'(z)Z_k\iv(z)+Z_k(z)\frac{E_k}{z-a_k}Z_k\iv(z)\nn\\
&=& \frac{Z_k(a_k)E_kZ_k\iv(a_k)}{z-a_k}\;+\;
\mbox{``analytic term at $z=a_k$''}\nn
\eqn
because $Z_k(z)$ is analytic and invertible at $z=a_k$.
At infinity we have $Y(z)=\La(z)Z_\iy(z)C$.
Hence, similarly,
\bqn
A(z) &=& Y'(\tz)Y\iv(\tz)
\quad=\quad
\La'(z)\La\iv(z)+
\La(z)Z_\iy'(z)Z_\iy\iv(z)\La\iv(z)\nn\\
&=& \frac{\diag(\ka_1,\dots,\ka_n)}{z}
+\La(z) R(z)\La\iv(z),\nn
\eqn
where $R(z)=Z_\iy'(z)Z_\iy\iv(z)$. Since $Z_\iy(z)$
is analytic and invertible,
$\|R(z)\|=O(1/z^2)$ as $z\to\iy$.
This implies (iv).

In order to prove the converse, we remark first
that (i) and (ii) of this proposition imply (i) and (ii) of
Definition \ref{d3.1} with the exception of the apparentness of
the singularity at infinity, which will follow from (iv).

Condition (iii) of this proposition in connection with
Lemma \ref{l4.1} immediately implies (iii) of Definition
\ref{d3.1}.
Note that the mere similarity of the residue to $E_k$
does not cause problems because
$(\tz-a_k)^{S_kE_kS_k\iv}=S_k(\tz-a_k)^{E_k}S_k\iv$.

Finally, let us consider the singularity at infinity. Because of
the asymptotics of $R(z)$ at infinity, there exists an analytic
solution $Z_\iy:D_\iy\cup\{\iy\}\to G\C^{n\times n}$
of the system
$Z_\iy'(z)=R(z)Z_\iy(z)$ near infinity. Now define
$Y_\iy(z)=\La(z)Z_\iy(z)$, $z\in D_\iy$. Then,
\bqn
Y_\iy'(z)Y_\iy\iv(z) &=& \La'(z)\La\iv(z)+
\La(z)Z_\iy'(z)Z_\iy\iv(z)\La\iv(z)\nn\\
&=& \frac{\diag(\ka_1,\dots,\ka_n)}{z}
+\La(z) R(z)\La\iv(z)
\quad=\quad A(z).\nn
\eqn
Hence $Y_\iy(z)$ is a solution of our original system considered
on the domain $D_\iy$. However, the relation of such a
``restricted'' solution to any solution of the system is given by
$Y(z)=Y_\iy(z)C$ with some $C$.
This is exactly condition (iv) of Definition \ref{d3.1}.
Obviously, this also shows that infinity is an apparent
singularity.
\end{proof}

Now we present the desired description in terms
of the matrix function $A(z)$.
For $k\in\Z$, $k\ge0$, let $\cP_k$ stand for the set of
all polynomials of degree less than or equal to $k$.
For $k\in\Z$, $k<0$, let $\cP_k$ stand for the set containing
only the function identically equal to zero.

\begin{theorem}\label{t4.3}
The system $Y'(\tz)=A(z)Y(\tz)$ is of standard form
with (admissible) data $[M_1,\dots,M_m]$ and $[E_1,\dots,E_m]$,
with singularities $a_1,\dots,a_m\in\C$
and
indices $\ka_1,\dots,\ka_n \in \Z$
if and only if the following conditions are satisfied:
\begin{itemize}
\item[(i)]
$A(z)=[A_{ij}(z)]_{i,j=1}^n$, where
\bqn\label{f4.1x}
A_{ij}(z) &=& \frac{\ka_i\delta_{ij}z^{m-1}+p_{ij}(z)}
{(z-a_1)(z-a_2)\cdots(z-a_m)},
\eqn
where $p_{ij}\in\cP_{m-2+\ka_i-\ka_j}$ and $\delta_{ij}$ is the
Kronecker symbol.
\item[(ii)]
For each $k=1,\dots,m$, the residue of $A(z)$ at $z=a_k$
is similar to the matrix $E_k$.
\item[(iii)]
The monodromy of the system is given by
$[\chi(\sigma_1),\dots,\chi(\sigma_m)]_\sim
=[M_1,\dots,M_m]_\sim$.
\end{itemize}
\end{theorem}
\begin{proof}
In regard to the previous proposition, we have to show that the
characterization given by (\ref{f4.1x}) is equivalent to the
property that $A(z)$ is analytic on
$\C\setminus\{a_1,\dots,a_m\}$, has simple poles at
$a_1,\dots,a_m$ and behaves at infinity as described in (iv) of
Proposition \ref{p4.2}.

The fact that $A(z)$ is analytic on $\C\setminus\{a_1,\dots,a_m\}$
and has simple poles at $a_1,\dots,a_m$ is equivalent to writing
the entries of $A(z)$ as
$A_{ij}(z)=q_{ij}(z)(z-a_1)\iv (z-a_2)\iv\cdots(z-a_m)\iv$,
where $q_{ij}$ is an entire analytic function.
The characterization of $q_{ij}$ as a certain polynomial
is now equivalent to (iv) of Proposition \ref{p4.2}.
\end{proof}

Considering the residues of $A(z)$
at $a_1,\dots,a_m$ and at infinity and taking traces,
we obtain the following simple relation,
\bqn\label{f4.1}
\ka_1+\dots+\ka_n &=& {\rm trace\,} E_1+\dots+
{\rm trace\,} E_m.
\eqn
Hence the sum of the indices is determined by the
the eigenvalues of $E_k$. This relation
is the counterpart to formula (\ref{f1.13})
for the total index of the factorization of a
piecewise constant matrix function.

\reseteqn
\section{Some general results in the irreducible case}
\label{s5}

In this section we present some necessary conditions for the indices of
systems of standard form with irreducible monodromy data $[M_1,\dots,M_m]$. By
establishing a relationship to scalar linear Fuchsian differential equations
of $n$-th order, we show that, in general, these conditions cannot be improved.

A subspace $X\subseteq\C^n$ is called an {\em invariant subspace} of a
collection of $n\times n$ matrices $\{M_\omega\}_{\omega\in\Omega}$ if the
image $M_\omega X$ is contained in $X$ for each $\omega\in\Omega$. We say that
a matrix function $F(z)$ has the invariant subspace $X$ if $F(z)X\subseteq X$
for each $z$.

A collection of matrices $\{M_\omega\}_{\omega\in\Omega}$ is
called {\em reducible} if there exists a non-trivial invariant
subspace $X$ (i.e., $X\neq\{0\}$ and $X\neq\C^n$). Otherwise, it
is called {\em irreducible}. The monodromy of a system of linear
differential equations is called reducible (resp., irreducible) if
the collection $\{\chi(\sigma)\}_{\sigma\in\Delta}$ (or,
equivalently, the collection
$[\chi(\sigma_1),\dots,\chi(\sigma_m)]$) is reducible (resp.,
irreducible) for some (hence each) monodromy representation $\chi$
of the system.

The following auxiliary result relates invariant subspaces of a matrix
function $A(z)$ with invariant subspaces for the solution and the monodromy of
the corresponding system.

\begin{lemma}\label{l5.1}
Let $A(z)$ be an analytic matrix function which has invariant subspaces
$\{X_r\}_{r\in R}$. Then there exists a solution of the system
$Y'(\tz)=A(z)Y(\tz)$ such that this solution and its corresponding monodromy
representation also have the invariant subspaces $\{X_r\}_{r\in R}$.
\end{lemma}
\begin{proof}
Let $P_r\in\C^{n\times n}$ be projections such that the image of $P_r$ is the
subspace $X_r$. Then $(I-P_r)A(z)P_r=0$, and we obtain
\bqn
(I-P_r)Y'(\tz)P_r &=&
(I-P_r)A(z)(I-P_r)Y(\tz)P_r.\nn
\eqn
Choose any point $\tz_0\in\tS$ and consider the
solution $Y(\tz)$ for which $Y(\tz_0)=I$. Hence
$(I-P_r)Y(\tz_0)P_r=0$. Because the matrix function
$(I-P_r)Y(\tz)P_r$ is a solution of
the above first order system, it follows that
$(I-P_r)Y(\tz)P_r=0$ for all $\tz\in \tS$.
This means that $X_r$ is an invariant subspace of $Y(\tz)$.
Hence $X_r$ is also an invariant subspace of the
inverse $Y\iv(\tz)$.
Let $\chi$ be the monodromy representation of $Y(\tz)$ and
$\sigma\in\Delta$. Then $\chi(\sigma)=Y\iv(\sigma(\tz))Y(\tz)$.
Hence $X_r$ is an invariant subspace of $\chi(\sigma)$.
\end{proof}

If we choose another solution of the system, $\wh{Y}(\tz)=Y(\tz)C$, then the
corresponding monodromy representation has the invariant subspaces $\{C\iv
X_r\}_{r\in R}$ (see formula (\ref{f2.6})).

Now we state the conditions on the indices in the case of
irreducible monodromy. These conditions (formulated in the setting
of vectors bundles and with a different proof) have already been
obtained by Bolibruch \cite{Bo2}. Recall that the indices are
ordered decreasingly, $\ka_1\ge\ka_2\ge\dots\ge\ka_n$.

\begin{theorem}\label{t5.2}
Let $Y'(\tz)=A(z)Y(\tz)$ be a system of standard
form with $m$ singularities $a_1,\dots,a_m\in\C$.
If the monodromy of this system is irreducible, then
the indices $\ka_1,\dots,\ka_n\in\Z$
satisfy, for each $k=1,\dots,n-1$, the condition
\bqn\label{f5.c}
\ka_{k}-\ka_{k+1} &\le& m-2.
\eqn
\end{theorem}
\begin{proof}
Suppose that this condition is not fulfilled for some $k$,
i.e., $\ka_{k}-\ka_{k+1} >m-2$. Because the indices are
ordered decreasingly, it follows that
$\ka_{j}-\ka_{i}>m-2$ for each $j=1,\dots,k$ and each
$i=k+1,\dots,n$. Theorem \ref{t4.3} implies that for those
$i$ and $j$, the entries $A_{ij}(z)$ of the matrix $A(z)$ vanish
identically. This means that the matrix $A(z)$ is of
block triangular form, i.e., it possesses the non-trivial invariant
subspace $X=\C^{k}\oplus\{0\}^{n-k}$. By the previous lemma,
the monodromy representation of some solution has also this
invariant subspace. Hence the monodromy of the system
is reducible, which contradicts the assumption.
\end{proof}

We want to emphasize that the above condition holds regardless of the location
of the singularities and the values of $[E_1,\dots,E_m]$. Reinterpreted in
terms of the factorization problem, this implies the independence of this
condition from the underlying space $L^p(\Ga)$.

It is also easy to see that this condition (even combined with the knowledge
of the total index $\ka$) is in general not sufficient to determine the values
of the partial indices uniquely. In this connection, we are going to show that
(\ref{f5.c}) cannot be improved.

For this purpose we introduce systems of standard form for which
the above relation holds with equality, i.e.,
$\ka_{k}-\ka_{k+1}=m-2$ for each $k=1,\dots,n-1$. We single out a
subclass of such systems which have -- generically -- irreducible
monodromy.

Suppose we are given indices $\ka_1,\dots,\ka_n\in\Z$ satisfying
$\ka_{k}-\ka_{k+1}=m-2$. It follows from Theorem \ref{t4.3} that systems
of standard form $Y'(\tz)=A(z)Y(\tz)$ with singularities $a_1,\dots,a_m$
and above indices have a matrix function $A(z)$ which is of the form

\ben\label{f5.2}
\frac{1}{p} \left(\ba{ccccc} \ka_1z^{m-1}+p_{11} & p_{12} &
p_{13} & \dots & p_{1n}\\ \alpha_1 & \ka_2z^{m-1}+p_{22} & p_{23} & & p_{2n} \\
0 & \alpha_2 & \ka_3z^{m-1}+p_{33} & \ddots & \vdots \\ \vdots & & \ddots &
\ddots & p_{n-1,n} \\ 0 & 0 & \dots & \alpha_{n-1} & \ka_nz^{m-1}+p_{nn}
\ea\right),
\een
where $p(z)=(z-a_1)\cdots(z-a_m)$, $\alpha_k\in\C$ and
$p_{ij}\in\cP_{(m-2)(1+j-i)}$. The condition on the polynomials
$p_{ij}$ follows from $\ka_i-\ka_j=(m-2)(j-i)$. Conversely, any system with
matrix $A(z)$ given by (\ref{f5.2}) such that the residues of $A(z)$ at
$a_1,\dots,a_m$ are non-resonant matrices is a system of standard form with
above indices and singularities and with certain admissible data.

If $\alpha_k=0$ for some $k=1,\dots,n-1$, then the monodromy of
the above system is reducible because $A(z)$ is of block
triangular form. Hence, as we are interested in irreducible
monodromy, we will focus on systems with
$\alpha_1\cdots\alpha_{n-1}\neq0$. A preliminary characterization
of the data of such systems is given next.

For $S\in\C^{n\times n}$, let $\mdg(S)$ stand
for the degree of the minimal polynomial of $S$,
i.e., the polynomial $p$ of smallest degree for which
$p(S)=0$.

\begin{proposition}\label{p5.3}
Assume that $Y'(\tz)=A(z)Y(\tz)$ is a system of standard form with
admissible data $[M_1,\dots,M_m]$ and $[E_1,\dots,E_m]$ where
$A(z)$ is given by (\ref{f5.2}) with
$\alpha_1\cdots\alpha_{n-1}\neq0$. Then $\mdg(E_k)=\mdg(M_k)=n$
for each $k=1,\dots,m$.
\end{proposition}
\begin{proof}
Because $\alpha_1\cdots\alpha_{n-1}\neq0$ the residue $\wh{E}_k$
of $A(z)$ at $z=a_k$ has nonzero entries on the diagonal below the
main diagonal. All entries below this diagonal are zero. By
considering the powers of $E_k$, it follows that
$\mdg(\wh{E}_K)=n$. By Proposition \ref{p4.2}(iii) we have
$\wh{E}_k\sim E_k$, which implies $\mdg(E_k)=n$. Because
$M_k\sim\exp(-2\pi iE_k)$ and $E_k$ is non-resonant, we obtain
$\mdg(M_k)=n$ by considering the Jordan normal forms.
\end{proof}

The next result relates the systems singled out above to scalar
linear Fuchsian differential equations. Therein, the local
exponents of the linear Fuchsian differential equations at the
points $a_1,\dots,a_m$ correspond to the eigenvalues of the
matrices $E_1,\dots,E_m$ of the given data.

This result has a strong relationship to a result obtained by
Bolibruch \cite{Bo2}. He examined the minimal number of additional apparent
singularities which a scalar Fuchsian differential equation must have in order
to realize given irreducible monodromy. The answer he obtained
involves the notion of the {\em maximal Fuchsian weight},
which is defined in terms of the splitting type
for a class of vector bundles.

\begin{theorem}\label{t5.4}
Let $[M_1,\dots,M_m]$ and $[E_1,\dots,E_m]$ be admissible data
of $n\times n$ matrices, and let $a_1,\dots,a_m\in\C$ be distinct points.
For $k=1,\dots,m$, denote by $\ep_k^{(1)},\dots,\ep_k^{(n)}$ the
eigenvalues of the matrices $E_k$.
Then the following two statements are equivalent:
\begin{itemize}
\item[(i)]
$[M_1,\dots,M_m]$ and $[E_1,\dots,E_m]$ is the data associated
to some system with $A(z)$ given by (\ref{f5.2})
and $\alpha_1\cdots\alpha_{n-1}\neq0$.
\item[(ii)]
$[M_1,\dots,M_m]$ is the monodromy of some $n$-th order linear
differential equation with Fuchsian singularities at
$a_1,\dots,a_m$ and local exponents given by
$\{\ep_k^{(1)},\dots,\ep_k^{(n)}\}$ for $z=a_k$,
$k=1,\dots,m-1$, and
$\{\ep_m^{(1)}-\ka_n,\dots,\ep_m^{(n)}-\ka_n\}$ for $z=a_m$.
\end{itemize}
\end{theorem}
\begin{proof}
(ii)$\Rightarrow$(i): Suppose we are given an $n$-th order linear
differential equation (\ref{f2.16}) with the above properties. Let
$y_1(\tz),\dots,y_n(\tz)$ be $n$ linear independent solutions. We
introduce \bqn W(\tz) &=& \left(\ba{cccc} y_1^{(n-1)}(\tz) &
y_2^{(n-1)}(\tz) & \dots & y_n^{(n-1)}(\tz) \\ \vdots & \vdots & &
\vdots \\ y_1'(\tz) & y_2'(\tz) & \dots & y_n'(\tz) \\ y_1(\tz)  &
y_2(\tz) & \dots & y_n(\tz) \ea\right),\qquad \tz\in\tS^*, \eqn
Then $W(\tz)$ is analytic on $\tS^*$, and it is well known that
the Wronskian $\det W(\tz)$ does not vanish on all of $\tS^*$. In
particular, $W(\tz)$ is a solution of the $n\times n$ system
$W'(\tz)=A_0(z)W(\tz)$, where \bqn\label{f5.4new} A_0(z) &=&
\left(\ba{cccc} -q_1(z) & -q_2(z) & \dots & -q_n(z) \\ 1 & 0 &
\dots & 0\\ \vdots & \ddots & \ddots & \vdots \\ 0 & \dots & 1 & 0
\ea\right),\qquad z\in S\setminus\{\iy\}, \eqn and
$q_1(z),\dots,q_n(z)$ are the coefficients of the scalar
differential equation. Obviously, the system for $A_0(z)$ has the
same singularities $a_1,\dots,a_m$ (and possibly an apparent
singularity at infinity) and the same monodromy as the scalar
differential equation.

Below we will replace the system $W'(\tz)=A_0(z)W(\tz)$ by a modified system
$Y'(\tz)=A(z)Y(\tz)$ which will have the desired properties. However, in
order to analyze the behavior at the singular points $a_k$, we will first
consider differently modified systems $Y_k(\tz)=A_k(z)Y_k(\tz)$ near these
singularities.

For $k=1,\dots,m$, we introduce \bqn P_k(z) &=&
\diag\Big((z-a_k)^{n-1},(z-a_k)^{n-2},\dots, (z-a_k),1\Big), \eqn
and write $q_j(z)=r_{jk}(z)/(z-a_k)^j$ for each $j=1,\dots,n$.
Because the singularities of the scalar equation are assumed to be
Fuchsian, the functions $r_{jk}(z)$ are analytic at $z=a_k$.

For $k=1,\dots,m-1$, we define $Y_k(\tz)=P_k(z)W(\tz)$. Then
$Y_k(\tz)$ is a solution of the system $Y_k'(\tz)=A_k(z)Y_k(\tz)$ with
$A_k(z)=P_k(z)A_0(z)P_k\iv(z)+P_k'(z)P_k\iv(z)$. From this we obtain that
$A_k(z)=\wh{E}_k/(z-a_k)+$``analytic term at $z=a_k$'' with
\bqn
\wh{E}_k &=&
\left(\ba{cccc}
-r_{1k}(a_k) & -r_{2k}(a_k) & \dots & -r_{nk}(a_k) \\
1 & 0 & \dots & 0\\
\vdots & \ddots & \ddots & \vdots \\
0 & \dots & 1 & 0
\ea\right)\;+\;
\diag(n-1,\dots,1,0).
\eqn
A straightforward calculation shows that the characteristic
equation $\det(\wh{E}_k-\rho I)=0$ for the matrix $\wh{E}_k$ coincides with
the indical equation (\ref{f2.ind}) for the local exponents of
the scalar equation with respect to the singularity $a_k$.
Hence the eigenvalues of $\wh{E}_k$ are just $\ep_k^{(1)},\dots,\ep_k^{(n)}$.

For $k=m$ we make a modification. We define
$Y_m(\tz)=(z-a_m)^{\ka_n}P_m(z)W(\tz)$, which is a solution of
$Y_m'(\tz)=A_m(z)Y_m(\tz)$.  Again
$A_m(z)=\wh{E}_m/(z-a_m)+$``analytic term at $z=a_m$'', but now $\wh{E}_m$
contains an additional term $\ka_n I$. Consequently, $\det(\wh{E}_m-\ka_n
I-\rho I)=0$ coincides with the indical equation for the singularity $a_m$.
Because the assumption on the local exponents is also different, we
obtain nevertheless that the eigenvalues of $\wh{E}_m$ are
$\ep_m^{(1)},\dots,\ep_m^{(n)}$.

In fact, the converse is also true. The scalar differential equation is
Fuchsian at $a_k$ if and only if the system $Y_k'(\tz)=A_k(z)Y_k(\tz)$
obtained in the above way is Fuchsian at $a_k$.

Next we are going to analyze the singularity at infinity. We make
a change of variables $z\mapsto\xi=(z-a_k)\iv$ and introduce
$\hat{y}_1(\xi),\dots,\hat{y}_n(\xi)$ by
$\hat{y}_j((z-a_m)\iv)=y_j(z)$ for  $j=1,\dots,n$. These new
functions are solutions of a modified scalar linear differential
equation. The point $\xi=0$ corresponds to the point $z=\iy$. By
assumption they are not singular points. Hence \bqn \wh{W}(\xi)
&=& \left(\ba{cccc} \hat{y}_1^{(n-1)}(\xi) &
\hat{y}_2^{(n-1)}(\xi) & \dots & \hat{y}_n^{(n-1)}(\xi) \\ \vdots
& \vdots & & \vdots \\ \hat{y}_1'(\xi) & \hat{y}_2'(\xi) & \dots &
\hat{y}_n'(\xi) \\ \hat{y}_1(\xi)  & \hat{y}_2(\xi) & \dots &
\hat{y}_n(\xi) \ea\right),\qquad \xi\in D_{0}, \eqn is analytic in
a neighborhood of $\xi=0$ and $\det\wh{W}(0)\neq0$. We remark that
the identity \bqn\label{f5.8} W(z) &=&
P_m\iv(z)SP_m\iv(z)\wh{W}((z-a_m)\iv), \qquad z\in D_\iy, \eqn
holds, which can be verified by a direct calculation. Here $S$ is
a certain invertible upper triangular $n\times n$ matrix. This
matrix does not depend on any parameters and can (in principle) be
evaluated explicitly.

Now we are prepared to define the system which will satisfy the
conditions of (i). Let
\bqn\label{f5.9}
U(z) &=& (z-a_m)^{\ka_n}P_1(z)\cdots P_{m-1}(z) S\iv P_m(z),
\eqn
and define $Y(\tz)=U(z)W(\tz)$. Because the function $U(z)$ is analytic and
invertible on $S\setminus\{\iy\}$, it follows from the afore-mentioned
properties of the system $W'(\tz)=A_0(z)W(\tz)$
that the new system $Y'(\tz)=A(z)Y(\tz)$ has
singularities only at $a_1,\dots,a_m$ (and possibly an apparent singularity at
infinity) and has the same monodromy as the scalar differential equation.

Next we analyze the local behavior of the new system at infinity.
Combining (\ref{f5.8}) and (\ref{f5.9}), it follows that \bqn Y(z)
&=& (z-a_m)^{\ka_n}P_1(z)\cdots
P_{m-1}(z)P_m\iv(z)\wh{W}((z-a_m)\iv) \nn\\[.5ex] &=&
\diag\Big(z^{\ka_n+(m-2)(n-1)},z^{\ka_n+(m-2)(n-2)},\dots,
z^{\ka_n+m-2},z^{\ka_n}\Big)Z_\iy(z)\qquad\label{f5.10} \eqn for
$z\in D_\iy$, where $Z_\iy$ is analytic and invertible on $D_\iy$.
Hence the behavior at infinity is as required for systems of
standard form with indices $\ka_1,\dots,\ka_n$. (Recall that
$\ka_k-\ka_{k+1}=m-2$ is assumed.)

In order to analyze the behavior at the singularity $a_k$, we make
the connection with the systems $Y_k'(\tz)=A_k(z)Y_k(\tz)$. First
write $U(z)=(z-a_m)^{\ka_n}\wh{U}(z)P_1(z)\cdots P_m(z)$, where
$\wh{U}(z)=P_1(z)\cdots P_{m-1}(z)S\iv P_{m-1}\iv(z)\cdots
P_1\iv(z)$ can be shown to be analytic and invertible on all of
$\C$. In fact, one uses that $S$ is an upper triangular matrix and
that $P_1(z)\cdot P_m(z)$ is a diagonal matrix function with
particular entries. Now, recalling how $Y_k(\tz)$ and $Y(\tz)$
were obtained from $W(\tz)$, it follows that
$Y(\tz)=U_k(z)Y_k(\tz)$ where  $U_k(z)$ is a certain matrix
function analytic and invertible on $D_{a_k}$. Hence the local
properties of $Y(\tz)$ and $Y_k(\tz)$ are essentially the same. In
particular, $A(z)$ has a simple pole at $z=a_k$ as so has
$A_k(z)$, and the residue of $A(z)$ at $z=a_k$ is similar to
$\wh{E}_k$.

Using the characterizations given in Proposition \ref{p4.2} and Theorem
\ref{t4.3}, it follows that $Y'(\tz)=A(z)Y(\tz)$ is a system of standard form
with ``local'' data $[\wh{E}_1,\dots,\wh{E}_m]$ and the properties mentioned
above. In particular, we obtain that $A(z)$ is of the form (\ref{f5.2}).

The assertion that $\alpha_1\cdots\alpha_{n-1}\neq0$ is a consequence of the
special structure of $A_0(z)$ given in (\ref{f5.4new}) and the fact that
$U(z)$ is an upper triangular matrix function. Observe that
\bqn\label{f5.11}
A(z) &=& U(z)A_0(z)U\iv(z)+U'(z)U\iv(z).
\eqn

Finally, we need to show that the ``local'' data is given also by
$[E_1,\dots,E_m]$, i.e., $E_k\sim\wh{E}_k$. In fact, Proposition \ref{p5.3}
implies that $\mdg(M_k)=\mdg(\wh{E}_k)=n$. Because by assumption
$M_k\sim\exp(-2\pi iE_k)$, we obtain $\mdg(E_k)=n$. Again by assumption,
the eigenvalues of $E_k$ are $\ep_k^{(1)},\dots,\ep_k^{(n)}$, hence they are
the same as those of $\wh{E}_k$. Now we can conclude $E_k\sim\wh{E}_k$.

\vspace{1ex}
(i)$\Rightarrow$(ii):
The essence of the above construction was the passage from the matrix
function $A_0(z)$ to the matrix function $A(z)$ by means of the transformation
(\ref{f5.11}) with $U(z)$. The fact that the entries below the first row of
$A_0(z)$ are of a special form is reason for the connection of $A_0(z)$ to the
scalar equation.

Now consider the singularities $a_1,\dots,a_m$ and the indices
$\ka_1,\dots,\ka_n$ as fixed and the functions $q_1(z),\dots,q_n(z)$ as
arbitrary. It follows from the above argumentation that the matrices $A(z)$
obtained by (\ref{f5.11}) are necessarily of the form (\ref{f5.2}). Moreover,
a straightforward computation (using the special structure of $U(z)$ and
$A_0(z)$) shows that the entries below the first row of $A(z)$ are also of a
special form, i.e., they only depend on $a_1,\dots,a_m$ and $\ka_n$ but not on
$q_1(z),\dots,q_n(z)$. Moreover, one obtains that
$\alpha_1\cdots\alpha_{n-1}\neq0$. In this sense, these entries of $A(z)$ have
to be considered as ``fixed''. In principle, they can be computed.

We want to make the converse transformation with $U\iv(z)$, i.e., to pass from
$A(z)$ to $A_0(z)$ by means of
\bqn\label{f5.12}
A_0(z) &=& U\iv(z)A(z)U(z)-U\iv(z)U'(z).
\eqn
Obviously, the matrix $A_0(z)$ will be only of the desired form
(\ref{f5.4new}) if the entries below the first row of $A(z)$ are chosen
``properly''. The crucial point is that if these entries of $A(z)$
are given indeed in the appropriate way, then the matrix $A_0(z)$ will be of
the form (\ref{f5.4new}) with certain $q_1(z),\dots,q_n(z)$. This can be
shown again by a straightforward computation.

Hence we cannot start with just any matrix $A(z)$ of the form
(\ref{f5.2}). Fortunately, the following result is true. Any
system $Y'(\tz)=A(z)Y(\tz)$ with $A(z)$ given by (\ref{f5.2}) and
$\alpha_1\cdots\alpha_{n-1}\neq0$ is equivalent to some system
with another matrix $\wh{A}(z)$ which is of the same form and
whose entries below the first row are prescribed arbitrarily.

In order to prove this statement we need only recall Theorem
\ref{t3.5}. We have to prove that there exists a matrix $V(z)$ of
the form (\ref{f.Vmat}) and the properties stated there such that
a transformation with $V(z)$ takes the matrix $A(z)$ to the matrix
$\wh{A}(z)$. Note that the entries of $V(z)$ are polynomials of a
certain degree. We can split this transformation into several
steps. First, let $V(z)$ be a diagonal matrix with suitable
diagonal entries such that $A(z)$ is taken into a matrix with
entries $\alpha_1,\dots,\alpha_{n-1}$. At the $k$-step,
$k=1,\dots,n$, we then choose the matrix $V(z)$ to have ones on
the main diagonal, suitable polynomials on the $k$-th diagonal
above the main diagonal and zero entries everywhere else. In fact,
a straightforward computation (analyzing $A\mapsto VAV\iv+V'V\iv$)
shows that these polynomials can be chosen in such a way that the
entries of the $(k-1)$-th diagonal above the main diagonal of
$A(z)$ (except the entry in the first row) can be modified in any
desired way. At the $k$-th step, all entries of $A(z)$ below this
diagonal remain unaltered. Hence after the $n$-th step we arrive
at a matrix $\wh{A}(z)$ of the same form and with desired entries
below the first row.

After this, we may assume that the matrix $A(z)$ of the system described
in (i) is of such a form that the transformation (\ref{f5.12}) gives a matrix
$A_0(z)$ which is of the form (\ref{f5.4new}) with certain
$q_1(z),\dots,q_n(z)$. Hence this new system gives raise to an $n$-th order
linear differential equation. We have to show that this equation has the
properties stated in (ii).

In fact, the argumentation given above can be reversed. Using the
properties of $U(z)$ (and employing the systems
$Y_k'(\tz)=A_k(z)Y_k(\tz)$ at an intermediate step) it follows
that $q_1(z),\dots,q_n(z)$ are analytic on $S\setminus\{\iy\}$ and
that $q_j(z)$ has at most a $j$-th order pole at $z=a_k$. Hence
the scalar equation is a Fuchsian equation and it follows also
from above that the local exponents are equal to the eigenvalues
of the residue of $A(z)$, which is a matrix similar to $E_k$.
Obviously, the monodromy is the same as for the given system.

The assertion that the scalar equation has no singularity at
infinity is only slightly more difficult. We know that the
solution $Y(\tz)$ of the system can be written as (\ref{f5.10}).
Introducing a modified scalar equation by change of variables as
above, we see that (\ref{f5.8}) holds with $W(\tz)$ being the
solution for the system with $A_0(z)$, or, equivalently, the
matrix (\ref{f5.3}) corresponding to the solution of the scalar
equation. Combining this with (\ref{f5.9}) and the fact that
$Y(\tz)=U(z)W(\tz)$, it is possible to conclude that $\wh{W}(\xi)$
is analytic and invertible near $\xi=0$. Hence $\xi=0$ is no
singular point for the modified scalar equation, and neither is
$z=\iy$ a singular point for the desired scalar equation.
\end{proof}

We remark in this connection that the identity
\bqn\label{f5.3}
{\rm trace}\,E_1+\dots+{\rm trace}\,E_m \;\;=\;\;
\ka_1+\dots+\ka_n &=& \frac{n(n-1)}{2}(m-2)+n\ka_n,\qquad
\eqn
which follows from (\ref{f4.1}) and
$\ka_{k}-\ka_{k+1}=m-2$, corresponds to
\bqn\label{f5.4}
\sum_{j=1}^n (\ep_1^{(j)}+\ep_2^{(j)}+\dots
+\ep_m^{(j)}-\ka_n) &=& \frac{n(n-1)}{2}(m-2),
\eqn
which is Fuchs' relation (\ref{f2.fuchs}) for the above
linear differential equation.

The previous theorem provides us with an implicit description of a class of
data which has indices satisfying $\ka_{k}-\ka_{k+1}=m-2$. (It does not
describe all such data because of the assumption
$\alpha_1\cdots\alpha_{n-1}\neq0$.) Namely, we may choose indices satisfying
$\ka_{k}-\ka_{k+1}=m-2$ and prescribe the eigenvalues of $E_1,\dots,E_m$ under
the restrictions (\ref{f5.3}) or (\ref{f5.4}). Then the corresponding linear
differential equations, which contain $(m-2)n(n+1)/2-(m-1)n+1$ free
parameters, give rise to a class of monodromy $[M_1,\dots,M_m]$. Although in
general there is no explicit description for this monodromy, one
knows that this data has the above indices.

The monodromy of linear Fuchsian differential equations is generically
irreducible. In fact, if the numbers $\ep_k^{(j)}$ are chosen such that
\ben
\sum_{k=1}^m \sum_{r=1}^R \ep_k^{(j_{r,k})} \;\notin\;\Z
\een
for all possibilities of $1\le R<n$ and $1\le j_{1,k}< j_{2,k}<\dots<j_{R,k}\le
n$, then the monodromy of the corresponding differential equation is always
irreducible. Hence the condition (\ref{f5.2}) cannot be improved in general.
We note, however, that the (non-generic) case of reducible monodromy is also
covered by the previous theorem.

The analogue of the condition $\mdg(M_k)=n$ stated in Proposition \ref{p5.3}
for scalar linear differential equations is well known. If the difference of
any two local exponents of a Fuchsian singularity is not a nonzero integer,
then the corresponding monodromy matrix $M_k$ satisfies $\mdg(M_k)=n$.

\reseteqn
\section{Some results for the {\boldmath $2\times 2$} matrix case}
\label{s6}

In this section we consider the $2\times 2$ matrix case with $m$
singularities. For general $m$ we obtain some information about
the indices corresponding to given data. For $m=3$, this
information provides us with a complete explicit answer. For
$m=4$, the answer can be given either explicitly or in terms of
the monodromy of second order linear Fuchsian differential
equations with four singularities.

In order to formulate our results we have to make some
preparations related to the ``reducibility type'' of the data.
Let $[M_1,\dots,M_m]$ and $[E_1,\dots,E_m]$ be admissible data
of $2\times 2$ matrices. Then the following (exclusive) cases
can occur:
\begin{itemize}
\item[(A)]
$[M_1,\dots,M_m]$ is irreducible;
\item[(B)]
$[M_1,\dots,M_m]$ possesses exactly one non-trivial invariant
subspace;
\item[(C)]
$[M_1,\dots,M_m]$ possesses at least two non-trivial invariant
subspaces.
\end{itemize}

In all these cases we define
\bqn
\ka &=& {\rm trace\,}E_1+\dots+{\rm trace\,}E_m.
\eqn

In case (B), the $m$-tuple $[M_1,\dots,M_m]$
is equivalent to
\ben\label{f6.1}
\left[
\left(
\ba{cc} \mu_1^{(1)} & \alpha_1 \\ 0 & \mu_1^{(2)} \ea
\right),
\dots,
\left(
\ba{cc} \mu_m^{(1)} & \alpha_m \\ 0 & \mu_m^{(2)} \ea
\right)
\right],
\een
where
$\mu_1^{(1)}\cdots\mu_m^{(1)}=\mu_1^{(2)}\cdots\mu_m^{(2)}=1$,
and $\alpha_1,\dots,\alpha_m\in\C$ satisfy the linear equation
\bqn\label{f6.2}
\alpha_1\mu_2^{(2)}\cdots\mu_m^{(2)}\;+\;
\mu_1^{(1)}\alpha_2\mu_3^{(2)}\cdots\mu_m^{(2)}
\;+\;\dots\;+\;
\mu_1^{(1)}\cdots\mu_{m-1}^{(1)}\alpha_m
&=& 0.
\eqn

In the representation (\ref{f6.1}), the values of $\mu_k^{(j)}$ are uniquely
determined. In fact, $\mu_k^{(1)}$ is the eigenvalue of $M_k$ corresponding to
the non-trivial invariant subspace and $\mu_k^{(2)}$ is the other eigenvalue.
Let $\ep_k^{(1)}$ and $\ep_k^{(2)}$ be the eigenvalues of $E_k$ numbered in
such a way that $\mu_k^{(j)}=\exp(-2\pi i\ep_k^{(j)})$. Because of the
non-resonance of $E_k$, their values are again uniquely determined,
and so are the following integers,
\ben\label{f6.int}
n_1\;\;=\;\;\ep_1^{(1)}+\dots+\ep_m^{(1)},\qquad\quad
n_2\;\;=\;\;\ep_1^{(2)}+\dots+\ep_m^{(2)}.
\een

In case (C), the $m$-tuple $[M_1,\dots,M_m]$ is diagonalizable,
i.e., it is equivalent to (\ref{f6.1}) with
$\alpha_1=\dots=\alpha_m=0$. Here, in contrast, the values of
$\mu_k^{(1)}$ and $\mu_k^{(2)}$ can be interchanged. Consequently,
the integers $n_1$ and $n_2$ are only defined up to change of
order. In the description of the main results later on, we will
therefore assume without loss of generality that $n_1\ge n_2$.

The following proposition provides some information about the
equivalence classes of $m$-tuples of the form (\ref{f6.1})

\begin{proposition}\label{p6.1}
Let $[M_1,\dots,M_m]$ and $[\wt{M}_1,\dots,\wt{M}_m]$ be two $m$-tuples of
the form (\ref{f6.1}). Assume that both $m$-tuples have the same values
$\mu_1^{(1)},\dots,\mu_m^{(1)}$ and $\mu_1^{(2)},\dots,\mu_m^{(2)}$, but
possibly different values  $\alpha_1,\dots,\alpha_m$ and
$\wt{\alpha}_1,\dots,\wt{\alpha}_m$, respectively.
Then these two $m$-tuples are equivalent if and only if
there exist $\lambda\in\C\setminus\{0\}$ and $\rho\in\C$
such that
\bqn\label{f6.3}
\alpha_k &=&
\lambda\wt{\alpha}_k+\rho(\mu_k^{(1)}-\mu_k^{(2)})
\eqn
for each $k=1,\dots,m$.
\end{proposition}
\begin{proof}
The $m$-tuples $[M_1,\dots,M_m]$ and $[\wt{M}_1,\dots,\wt{M}_m]$ are equivalent
if and only if there exists a matrix $C\in G\C^{2\times 2}$ such that
$\wt{M}_k=CM_kC\iv$. If all matrices $M_1,\dots,M_m$ are scalar matrices, then
the assertion is trivial. Otherwise, if at least one $M_k$ is not scalar,
one rewrites $\wt{M}_kC=CM_k$, and now a straightforward computation
shows that $C$ is an upper triangular matrix. Using this, the relation
$\wt{M}_k=CM_kC\iv$ is easily seen to be equivalent to (\ref{f6.3}).
\end{proof}

As a conclusion, an $m$-tuple of the
form (\ref{f6.1}) is diagonalizable if and only if the vector
$\alpha=(\alpha_1,\dots,\alpha_m)$ is a multiple of the vector
$(\mu_1^{(1)}-\mu_1^{(2)},\dots,\mu_m^{(1)}-\mu_m^{(2)})$.

All $m$-tuples of the form (\ref{f6.1}) are parameterized by the
vector $\alpha$. Because of the linear condition (\ref{f6.2}), all
``admissible'' vectors $\alpha$ are taken from an
$(m-1)$-dimensional subspace of $\C^m$. The relation (\ref{f6.3})
establishes an equivalence relation between such vectors. The
statement of this proposition is that there is a one-to-one
correspondence between equivalence classes of such vectors and
equivalence classes of such $m$-tuples. Elaborating on
(\ref{f6.3}), it is easy to see that these equivalence classes can
be identified with $\sP^{m-2}\cup\{0\}$ in case
$\mu_k^{(1)}=\mu_k^{(2)}$ for all $k=1,\dots,m$ and with
$\sP^{m-3}\cup\{0\}$ in case $\mu_k^{(1)}\neq\mu_k^{(2)}$ for some
$k=1,\dots,m$. Here $\sP^{n}$ stands for the $n$-dimensional
complex projective space, and $\{0\}$ symbolizes the single
equivalence class containing the vector $\alpha=0$. We remark that
the singleton $\{0\}$ corresponds exactly to the class of matrices
$[M_1,\dots,M_m]$ which are diagonalizable.

In the following theorem, we consider all systems of standard form
$Y'(\tz)=A(z)Y(\tz)$ with indices
$\ka_1,\ka_2\in\Z$, $\ka_1\ge\ka_2$,
for which $A(z)$ can be written as
\bqn\label{f6.4}
A(z) &=& \frac{1}{p(z)}
\left(\ba{cc}
\ka_1z^{m-1}+p_1(z) & p_{12}(z) \\
0 & \ka_2z^{m-1}+p_2(z) \ea\right),
\eqn
where $p(z)=(z-a_1)\cdots(z-a_m)$, $p_1,p_2\in\cP_{m-2}$,
$p_{12}\in\cP_{m-2+\ka_1-\ka_2}$.
We will completely characterize the data which is associated
to such systems.

\begin{theorem}\label{t6.2} Let $\ka_1,\ka_2\in\Z$, $\ka_1\ge\ka_2$, and let
$[M_1,\dots,M_m]$ and $[E_1,\dots,E_m]$ be admissible data of $2\times 2$
matrices. Then the following two statements are equivalent:
\begin{itemize}
\item[(i)]
$[M_1,\dots,M_m]$ and $[E_1,\dots,E_m]$ is the data associated to
some system of the form (\ref{f6.4}).
\item[(ii)]
$[M_1,\dots,M_m]$ is
equivalent to some $m$-tuple of the form (\ref{f6.1}) with $n_1=\ka_1$,
$n_2=\ka_2$, where $n_1$ and $n_2$ are defined by (\ref{f6.int}).
\end{itemize}
\end{theorem}
\begin{proof}
(i)$\Rightarrow$(ii): Because $A(z)$ is of
triangular form we obtain from Lemma \ref{l5.1} that there exists a solution
$Y(\tz)$ which is also of triangular form. So let us write
$Y'(\tz)=A(z)Y(\tz)$ as
\bqn\label{f6.6}
\left(\ba{cc} y_1'(\tz) &
y_{12}'(\tz) \\ 0 & y_2'(\tz) \ea \right) &=& \left(\ba{cc} a_1(z) & a_{12}(z)
\\ 0 & a_2(z) \ea \right) \left(\ba{cc} y_1(\tz) & y_{12}(\tz) \\ 0 & y_2(\tz)
\ea \right).
\eqn
For $j=1,2$, we have
\bqn\label{f6.7} a_j(z) &=&
\frac{\ka_jz^{m-1}+p_j(z)}{p(z)} \;\;=\;\; \sum_{k=1}^m
\frac{\ep_k^{(j)}}{z-a_k},
\eqn
where the last equality defines the numbers
$\ep_1^{(j)},\dots,\ep_m^{(j)}$. They satisfy
$\ka_j=\ep_1^{(j)}+\dots+\ep_m^{(j)}$. Because $E_k$ is similar to the residue
of $A(z)$ at $z=a_k$, the numbers $\ep_k^{(1)}$ and $\ep_k^{(2)}$ are the
eigenvalues of $E_k$. Once we have shown that  $\mu_k^{(j)}=\exp(-2\pi
i\ep_k^{(j)})$, it follows that $\ka_1=n_1$ and $\ka_2=n_2$.

In fact, from (\ref{f6.7}) we obtain
that the solutions of $y_j'(\tz)=a_j(z)y_j(\tz)$ are
given by
$y_j(\tz) = \prod_{k=1}^m (\tz-a_k)^{\ep_k^{(j)}}$ up to a redundant constant.
Regardless of the precise expression for $y_{12}(\tz)$,
it follows that
the monodromy representation for $Y(\tz)$ is
given by (\ref{f6.1}) with $\mu_k^{(j)}=\exp(-2\pi i\ep_k^{(j)})$
and with certain $\alpha_k$. Hence $[M_1,\dots,M_m]$ is equivalent to
(\ref{f6.1}) and the $\mu_k^{(j)}$ and $\ep_k^{(j)}$
are indeed properly ordered.

(ii)$\Rightarrow$(i):
We consider all $m$-tuples
of the form (\ref{f6.1}) with fixed $\mu_k^{(j)}$, and
we fix also the numbers $\ep_k^{(j)}$ (i.e., the eigenvalues
of $E_k$).
Observing that $\ka_1=n_1$ and $\ka_2=n_2$, we
introduce the functions $a_1(z)$ and $a_2(z)$ by
(\ref{f6.7}).

In what follows we consider all systems of the form (\ref{f6.6})
with those $a_1(z)$ and $a_2(z)$ and with
$a_{12}(z)=p_{12}(z)/p(z)$ where $p_{12}(z)$ runs through
all polynomials in $\cP_{m-2+\ka_1-\ka_2}$. These systems are
indeed of standard form because $\ka_1\ge\ka_2$ and
because $E_1,\dots,E_m$ are assumed to be non-resonant.
{}From the first part of this proof we already know that
the monodromy representation of the corresponding solutions
are of the form (\ref{f6.1}) with certain
vectors $\alpha=(\alpha_1,\dots,\alpha_m)$.
It remains to show that if $p_{12}$ ranges through all of
$\cP_{m-2+\ka_1-\ka_2}$, then the corresponding vectors
$\alpha$ take values in all equivalence classes
(defined by the equivalence relation (\ref{f6.3})).

We first examine the question when two systems with
$A(z)$ and $\wt{A}(z)$ given by
$$
A(z) \;\;=\;\;\left(\ba{cc} a_1(z) & p_{12}(z)/p(z) \\
0 & a_2(z) \ea\right),\qquad
\wt{A}(z) \;\;=\;\;\left(\ba{cc} a_1(z) & \tilde{p}_{12}(z)/p(z) \\
0 & a_2(z) \ea\right)
$$
with $p_{12},\tilde{p}_{12}\in\cP_{m-2+\ka_1-\ka_2}$ are
equivalent. By definition this is the case if and only if
there exists a $V$ such that $\wt{A}=VAV\iv+V'V\iv$, where
$V$ is of the form
$$ V(z) \;\;=\;\;
\left(\ba{cc} \lambda_1 & v(z) \\ 0 & \lambda_2 \ea\right)
$$
with $\lambda_1,\lambda_2\in\C\setminus\{0\}$ and
$v\in\cP_{\ka_1-\ka_2}$ if $\ka_1>\ka_2$, and
$V(z)=V\in G\C^{2\times 2}$ if $\ka_1=\ka_2$.
If $\ka_1=\ka_2$, then $\wt{A}V=VA$ implies that
$V$ has to be an upper triangular matrix (except for the case
where $A(z)$ is given with $a_1(z)=a_2(z)$ and $p_{12}(z)=0$,
which can easily be dealt with separately).
Using this information about $V(z)$, it follows that the above
systems are equivalent if and only if
$\wt{p}_{12}=\lambda p_{12}+v(a_2-a_1)p+v'p$ for some
$\lambda\in\C\setminus\{0\}$ and $v\in\cP_{\ka_1-\ka_2}$.

This, in turn, can be rephrased by saying that
$\wt{p}_{12}-\lambda p_{12}\in\im \Xi$ for some
$\lambda\in\C\setminus\{0\}$, where $\im\Xi$ is the image
of the linear mapping
$\Xi:\cP_{\ka_1-\ka_2}\to\cP_{m-2+\ka_1-\ka_2},
v\mapsto v(a_2-a_1)p+v'p$.  We conclude that
$\dim\im\Xi\le \ka_1-\ka_2+1$.
If $a_1=a_2$ (hence $\ka_1=\ka_2$), then $\Xi$ is even
the zero mapping. Now decompose
$\cP_{m-2+\ka_1-\ka_2}=\im\Xi\oplus X$ as a direct sum,
and remark that $\dim X\ge m-2$ and even $\dim X\ge m-1$ if
$a_1=a_2$. Finally, we arrive at the following statement:
if $p_{12},\wt{p}_{12}\in X$, then the above systems
are equivalent if and only if $\wt{p}_{12}=\lambda p_{12}$
for some $\lambda\in\C\setminus\{0\}$.

Now choose any basis $p_{12}^{(1)},\dots,p_{12}^{(d)}$ in
$X$ where $d=\dim X$. Assume that a solution of the system
(\ref{f6.6}) with $a_{12}=p_{12}^{(r)}/p$ is given by
some $y_{12}=y_{12}^{(r)}$, $r=1,\dots,d$. The corresponding
monodromy representation is supposed to be given by
(\ref{f6.1}) with some vector $\alpha=\alpha^{(r)}\in\C^{m}$.
Now consider an arbitrary polynomial $p_{12}\in X$ and write
$p_{12}=\sum \gamma_r p_{12}^{(r)}$. A straightforward
computation shows that a solution of (\ref{f6.6}) with
$a_{12}=p_{12}/p$ is given with
$y_{12}=\sum \gamma_r y_{12}^{(r)}$
(notice that there are certainly further solutions).
Moreover, the monodromy representation
of this solution is given by (\ref{f6.1}) with
$\alpha=\sum \gamma_r \alpha^{(r)}$. In this way, we have defined
a linear mapping $\Lambda:X\to \C^{m}, p_{12}\mapsto \alpha$
from $X$ into the set of all admissible vectors $\alpha$,
which characterize the reducible monodromy.

Now we invoke the fact that there is a one-to-one correspondence
between equivalence classes of systems of standard form and
equivalence classes of given data (see Corollary \ref{c3.6}). This
fact specialized to the mapping $\Lambda$ means the following: for
$p_{12},\wt{p}_{12}\in X$, the vector $\Lambda(p_{12})$ is
equivalent to the vector $\Lambda(\wt{p}_{12})$ in the sense of
(\ref{f6.3}) if and only if $\wt{p}_{12}=\lambda p_{12}$ for some
$\lambda\in\C\setminus\{0\}$. This allows us easily to conclude
that the kernel of $\Lambda$ is trivial. Hence, in case $a_1=a_2$,
we have $\dim\im\Lambda=d\ge m-1$, and we are done because the set
of all admissible vectors $\alpha$ is an ($m-1$)-dimensional
subspace of $\C^m$. On the other hand, $a_1\neq a_2$ means that
$\ep_{k}^{(1)}\neq\ep_k^{(2)}$ for some $k=1,\dots,m$. Hence (by
the non-resonance assumption) the vector
$\mu=(\mu_1^{(1)}-\mu_1^{(2)},\dots,\mu_m^{(1)}-\mu_m^{(2)})$,
occurring in (\ref{f6.3}), is not the zero vector. Because
$\dim\im\Lambda=d\ge m-2$, the assertion is proved if we have
shown that $\mu\notin\im\Lambda$. Indeed, if
$\Lambda(p_{12})=\mu$, then the equivalence of $\mu$ to the zero
vector implies that $p_{12}=0$. Hence $\mu=0$, which is a
contradiction.
\end{proof}

The reader should realize the meaning of the condition
$n_1=\ka_1$ and $n_2=\ka_2$ in the previous theorem.
Because of the necessary assumption $\ka_1\ge\ka_2$, not all
reducible data is covered by (ii) of the theorem,
but only those for which $n_1\ge n_2$.
However, the case of diagonalizable
data is completely covered.

The next corollary resumes this, and gives in addition
some (but in general not a complete)
information about the indices in the remaining cases.

\begin{corollary}\label{c6.3}
Let $[M_1,\dots,M_m]$ and $[E_1,\dots,E_m]$ be admissible data
of $2\times 2$ matrices.
\begin{itemize}
\item[(i)]
If $[M_1,\dots,M_m]$ is diagonalizable, or if
$[M_1,\dots,M_m]$ possesses exactly one non-trivial
invariant subspace and $n_1\ge n_2$ holds,
then the indices are $\ka_1=n_1$ and $\ka_2=n_2$.
\item[(ii)]
If $[M_1,\dots,M_m]$ is irreducible, or if
$[M_1,\dots,M_m]$ possesses exactly one non-trivial
invariant subspace and $n_1< n_2$ holds,
then the indices fulfill the conditions:
\ben\label{f6.9}
\ka_1-\ka_2\;\;\le\;\; m-2\quad\mbox{ and }\quad
\ka_1+\ka_2\;\;=\;\;\ka
\een
\end{itemize}
\end{corollary}
\begin{proof}
Part (i) is a conclusion of the implication (ii)$\Rightarrow$(i)
of Theorem \ref{t6.2}. The second condition of part (ii) is just
relation (\ref{f4.1}).

Now assume that $\ka_1-\ka_2>m-2$. Then Theorem \ref{t4.3} implies
that the matrix $A(z)$ of the corresponding system is the form
(\ref{f6.4}). Using the implication (i)$\Rightarrow$(ii) of
Theorem \ref{t6.2}, it follows that $[M_1,\dots,M_m]$ is of the
form (\ref{f6.1}) with $n_1\ge n_2$, contradicting the assumption.
\end{proof}

Hence, for data satisfying the assumption in (i),
we have a complete knowledge of the indices.
Note that they do not depend on the location of the
singularities.

For data satisfying the assumption in (ii), we will give
precise information for the cases $m=3$ and $m=4$.
In case $m=3$, we need only invoke (\ref{f6.9}).

\begin{corollary}
Let $[M_1,M_2,M_3]$ and $[E_1,E_2,E_3]$ be admissible data
of $2\times 2$ matrices. Assume that
$[M_1,M_2,M_3]$ is irreducible, or that
$[M_1,M_2,M_3]$ possesses exactly one non-trivial
invariant subspace and $n_1< n_2$ holds.
\begin{itemize}
\item[(a)]
If $\ka$ is odd, then
$\ka_1=(\ka+1)/2$ and $\ka_2=(\ka-1)/2$.
\item[(b)]
If $\ka$ is even, then
$\ka_1=\ka_2=\ka/2$.
\end{itemize}
\end{corollary}

The results for $m=3$ presented here are in accordance with the
result obtained in \cite{SpTa} for the factorization problem.
Actually, they are a generalization
of the former result because the matrices
$E_1,E_2,E_3$ can be chosen less restrictive.
(We require only non-resonance instead that the real parts
of the eigenvalues are contained in an interval of length one.)

\begin{corollary}
Let $[M_1,\dots,M_4]$ and $[E_1,\dots,E_4]$ be admissible
data of $2\times 2$ matrices, and let
$a_1,\dots,a_4\in\C$ be distinct points.
Assume that $[M_1,\dots,M_4]$ is irreducible, or that
$[M_1,\dots,M_4]$ possesses exactly one non-trivial
invariant subspace and $n_1<n_2$ holds.
Let $\ep_k^{(1)},\ep_k^{(2)}$ be the
eigenvalues of $E_k$, $k=1,\dots,4$.
\begin{itemize}
\item[(a)]
If $\ka$ is odd, then
$\ka_1=(\ka+1)/2$ and $\ka_2=(\ka-1)/2$.
\item[(b)]
If $\ka$ is even, then either
$\ka_1=\ka_2=\ka/2$, or $\ka_1=\ka/2+1$ and $\ka_2=\ka/2-1$.
\item[(b*)]
For $\ka$ even, the indices are
$\ka_1=\ka/2+1$ and $\ka_2=\ka/2-1$ if and only if
$[M_1,\dots,M_4]$ is the monodromy of some second
order linear differential equation with Fuchsian singularities
$a_1,\dots,a_4$ and local exponents given by
$\{\ep_k^{(1)},\ep_k^{(2)}\}$ for $z=a_k$, $k=1,2,3$,
and $\{\ep_4^{(1)}+1-\ka/2,\ep_4^{(2)}+1-\ka/2\}$ for
$z=a_4$.
\end{itemize}
\end{corollary}
\begin{proof}
Again, assertions (a) and (b) are an immediate consequence of
(\ref{f6.9}). In order to prove (b*) we argue as follows.
The indices for the given data are
$\ka_1=\ka/2+1$ and $\ka_2=\ka/2-1$ if and only if
this data is associated to some system of standard form
for which $A(z)$ can be written as
\bqn
A(z) &=&
\frac{1}{p(z)}
\left(\ba{cc}
(\ka/2+1)z^3+p_1(z) & p_{12}(z) \\
\alpha & (\ka/2-1)z^3+p_2(z) \ea\right)\nn
\eqn
with $p(z)=(z-a_1)\cdots(z-a_4)$,
$\alpha\in\C$, $p_1,p_2\in\cP_2$ and
$p_{12}\in\cP_4$ (see Theorem \ref{t4.3}).
If $\alpha=0$, then we can apply Theorem
\ref{t6.2} and conclude that the monodromy is reducible and
$n_1=\ka/2+1>\ka/2-1=n_2$, which contradicts our assumption
on the data. Hence $\alpha\neq0$. Now we can use
Theorem \ref{t5.4}, which shows the desired equivalence.
\end{proof}

For the case where assumption (b) of the previous
corollary applies, we are led to the monodromy
problem for second order Fuchsian differential equations
with four singular points.  Note that the monodromy depends
on the location of the singularities. Hence so will the
corresponding indices. Unfortunately, no explicit
answer is known to this monodromy problem.
Only one particular case can be answered immediately
by means of
the remark made in the last paragraph of Section \ref{s5}.
Namely, if one of the matrices $M_1,\dots,M_4$ is a
scalar matrix (i.e., a multiple of $I$), then the
indices are necessarily $\ka_1=\ka_2=\ka/2$.

\reseteqn
\section{The {\boldmath $3\times 3$} matrix case with {\boldmath $3$} singularities}
\label{s7}

In this section we consider the $3\times 3$ matrix case with $3$ singularities.
In principle, the results are similar to the $2\times 2$ matrix case with $4$
singularities, although it requires some more effort to describe them.
In some cases of admissible data, the indices can be determined explicitly,
in other cases they depend on the description of the monodromy of third
order Fuchsian linear differential equations with three singularities.

\subsection{Classification of reducibility type}

We start again with some preparations related to the ``reducibility type''
of the data. Let $[M_1,M_2,M_3]$ and $[E_1,E_2,E_3]$ be admissible data
of $3\times 3$ matrices. Then the following (mutually exclusive) cases
can occur:
\begin{itemize}
\item[(A)]
$[M_1,M_2,M_3]$ is irreducible;
\item[(B-1)]
$[M_1,M_2,M_3]$ possesses exactly one non-trivial
invariant subspace $V_1$ such that $\dim V_1=1$;
\item[(B-2)]
$[M_1,M_2,M_3]$ possesses exactly one non-trivial
invariant subspace $V_2$ such that $\dim V_2=2$;
\item[(B-3)]
$[M_1,M_2,M_3]$ possesses exactly two non-trivial
invariant subspaces $V_1$ and $V_2$ such that $\dim V_1=1$, $\dim V_2=2$
and $V_1\cap V_2=\{0\}$;
\item[(C)]
$[M_1,M_2,M_3]$ possesses exactly two non-trivial
invariant subspaces $V_1$ and $V_2$ such that $\dim V_1=1$, $\dim V_2=2$
and $V_1\subset V_2$;
\item[(C-1)]
$[M_1,M_2,M_3]$ possesses (at least) two one-dimensional
invariant subspaces $V_1^{(a)}$ and $V_1^{(b)}$ and exactly one
two-dimensional invariant subspace $V_2=V_1^{(a)}\oplus V_1^{(b)}$;
\item[(C-2)]
$[M_1,M_2,M_3]$ possesses (at least) two two-dimensional
invariant subspaces $V_2^{(a)}$ and $V_2^{(b)}$ and exactly one
one-dimensional invariant subspace $V_1=V_2^{(a)}\cap V_2^{(b)}$;
\item[(C-3)]
$[M_1,M_2,M_3]$ possesses (at least) two one-dimensional invariant
subspaces and two two-dimensional invariant subspaces, where all
one-dimensional invariant subspaces are contained in a
two-dimensional subspace and the intersection of all
two-dimensional invariant subspaces is a one-dimensional subspace;
\item[(D)]
$[M_1,M_2,M_3]$ possesses (at least) three one-dimensional invariant subspaces
$V_1^{(a_)}$, $V_1^{(b)}$ and $V_1^{(c)}$ such that
$V_1^{(a_)}\oplus V_1^{(b)}\oplus V_1^{(c)}=\C^3$.
\end{itemize}
The readers may convince themselves that this classification is complete.

In all these cases we define
\bqn\label{f7.1}
\ka &=& {\rm trace\,}E_1+{\rm trace\,}E_2+{\rm trace\,}E_3.
\eqn

In case (B-1), the triple $[M_1,M_2,M_3]$ is equivalent to
\ben\label{f7.2}
\left[
\left(\ba{cc}\mu_1 & \alpha_1 \\ 0 & M_1^{\#} \ea \right),
\left(\ba{cc}\mu_2 & \alpha_2 \\ 0 & M_2^{\#} \ea \right),
\left(\ba{cc}\mu_3 & \alpha_3 \\ 0 & M_3^{\#} \ea \right)
\right]
\een
with numbers $\mu_k$, $1\times 2$ vectors $\alpha_k$ and $2\times 2$ matrices
$M_k^{\#}$, which satisfy $\mu_1\mu_2\mu_3=1$, $M_1^{\#}M_2^{\#}M_3^{\#}=I$
and
\bqn\label{f7.3}
\alpha_1M_2^{\#}M_3^{\#}+\mu_1\alpha_2M_3^{\#}+\mu_1\mu_2\alpha_3 &=& 0.
\eqn
The triple $[M_1^{\#},M_2^{\#},M_3^{\#}]$ is irreducible.
For each $k=1,2,3$, the value $\mu_k$ is uniquely determined and so is the
eigenvalue $\ep_k$ of $E_k$ that satisfies $\mu_k=\exp(-2\pi i\ep_k)$.
Hence the following integers are uniquely defined,
\bqn\label{f7.4}
\nu\;\;=\;\;\ep_1+\ep_2+\ep_3,\qquad
N\;\;=\;\;\ka-\nu.
\eqn

In case (B-2), the triple $[M_1,M_2,M_3]$ is equivalent to
\ben\label{f7.5}
\left[
\left(\ba{cc}M_1^{\#} & \alpha_1 \\ 0 & \mu_1 \ea \right),
\left(\ba{cc}M_2^{\#} & \alpha_2 \\ 0 & \mu_2 \ea \right),
\left(\ba{cc}M_3^{\#} & \alpha_3 \\ 0 & \mu_3 \ea \right)
\right]
\een
with numbers $\mu_k$, $2\times 1$ vectors $\alpha_k$ and $2\times 2$ matrices
$M_k^{\#}$ having the same properties as above, but with (\ref{f7.3}) replaced
by
\bqn\label{f7.6}
\alpha_1\mu_2\mu_3+M_1^{\#}\alpha_2\mu_3+M_1^{\#}M_2^{\#}\alpha_3 &=& 0.
\eqn
We also define the integers $\nu$ and $N$ by (\ref{f7.4}).

In case (B-3), the triple $[M_1,M_2,M_3]$ is equivalent both to
(\ref{f7.2}) and (\ref{f7.5}) with $\alpha_1=\alpha_2=\alpha_3=0$.
Again, the integers $\nu$ and $N$ are uniquely defined by (\ref{f7.4}).

In case (C), the triple $[M_1,M_2,M_3]$ is equivalent to
\ben\label{f7.7}
\left[
\left(\ba{ccc}
\mu_1^{(1)} & \alpha_1 &\gamma_1 \\
0 & \mu_1^{(2)} & \beta_1 \\ 0 & 0 & \mu_1^{(3)} \ea \right),
\left(\ba{ccc}
\mu_2^{(1)} & \alpha_2 &\gamma_2 \\
0 & \mu_2^{(2)} & \beta_2 \\ 0 & 0 & \mu_2^{(3)} \ea \right),
\left(\ba{ccc}
\mu_3^{(1)} & \alpha_3 &\gamma_3 \\
0 & \mu_3^{(2)} & \beta_3 \\ 0 & 0 & \mu_3^{(3)} \ea \right),
\right]
\een
where the numbers $\mu_k^{(j)}$ are uniquely determined and satisfy
$\mu_1^{(j)}\mu_2^{(j)}\mu_3^{(j)}=1$ for each $j=1,2,3$ and the
numbers $\alpha_k,\beta_k,\gamma_k$ satisfy the following conditions:
\bqn\label{f7.lin1}
\alpha_1\mu_2^{(2)}\mu_3^{(2)}+
\mu_1^{(1)}\alpha_2\mu_3^{(2)}+
\mu_1^{(1)}\mu_2^{(1)}\alpha_3
&=& 0,\\[.5ex]\label{f7.lin2}
\beta_1\mu_2^{(3)}\mu_3^{(3)}+
\mu_1^{(2)}\beta_2\mu_3^{(3)}+
\mu_1^{(2)}\mu_2^{(2)}\beta_3
&=& 0,\\[.5ex]\label{f7.lin3}
\mu_1^{(1)}\mu_2^{(1)}\gamma_3+
\mu_1^{(1)}\gamma_2\mu_3^{(3)}+
\gamma_1\mu_2^{(3)}\mu_3^{(3)} &=&
-\mu_1^{(1)}\alpha_2\beta_3-
\alpha_1\mu_2^{(2)}\beta_3-
\alpha_1\beta_2\mu_3^{(3)}.\qquad
\eqn
Note that $[M_1,M_2,M_3]$ is not equivalent to a
triple of the form (\ref{f7.7}) where all $\alpha_k$ are zero or all $\beta_k$
are zero. We introduce the integers $n_1,n_2,n_3$ by
\bqn\label{f7.8}
n_j &=& \ep_1^{(j)}+\ep_2^{(j)}+\ep_3^{(j)},\qquad j=1,2,3,
\eqn
where $\ep_k^{(j)}$ are the eigenvalues of
$E_k$ numbered in such a way that $\mu_k^{(j)}=\exp(-2\pi i\ep_k^{(j)})$.
The integers $n_1,n_2,n_3$ are uniquely defined in this case.

In case (C-1), the triple $[M_1,M_2,M_3]$ is equivalent to
(\ref{f7.7}) with $\alpha_k=0$ (but not all $\beta_k$ or all $\gamma_k$
can be made zero). The values $\mu_3^{(j)}$ are unique, but the values
$\mu_1^{(j)}$ and $\mu_2^{(j)}$ can be interchanged. Hence the integer $n_3$
is uniquely determined, and we will assume without loss of generality that
$n_1\ge n_2$.

In case (C-2), the triple $[M_1,M_2,M_3]$ is equivalent to
(\ref{f7.7}) with $\beta_k=0$ (but not all $\alpha_k$ or all $\gamma_k$
can be made zero). The values $\mu_1^{(j)}$ are unique, but the values
$\mu_2^{(j)}$ and $\mu_3^{(j)}$ can be interchanged. Hence the integer $n_1$
is uniquely determined, and we will assume without loss of generality that
$n_2\ge n_3$.

In case (C-3), the triple $[M_1,M_2,M_3]$ is equivalent to
(\ref{f7.7}) with $\beta_k=\gamma_k=0$ (but not all $\alpha_k$ can be made
zero). Here the values of $\mu_k^{(j)}$ are again unique for all $j,k=1,2,3$.
We define the integers $n_1,n_2,n_3$ by (\ref{f7.8}), but to remark the
different roles that these numbers play, we use the notation
\bqn\label{f7.9}
\nu_1\;=\;n_1,\quad \nu_2\;=\;n_2,\quad \nu_{\#}=n_3.
\eqn
Notice that (by permuting rows and columns) $[M_1,M_2,M_3]$ is also
equivalent to a triple of the form (\ref{f7.7}) with
$\alpha_k=\gamma_k=0$ or with $\alpha_k=\beta_k=0$. This observation will play
some role later on in treating this case. The numbers $\nu_1,\nu_2,\nu_{\#}$
would then have to be defined appropriately.

Finally, in case (D), the triple $[M_1,M_2,M_3]$ is diagonalizable,
i.e., equivalent to (\ref{f7.7}) with $\alpha_k=\beta_k=\gamma_k=0$.
The values of $\mu_k^{(j)}$ may be permuted and we will therefore assume
without loss of generality that the numbers $n_1,n_2,n_3$ defined as above
are such that $n_1\ge n_2\ge n_3$.

\subsection{Systems of triangular form}

We will first consider monodromy data which is triangularizable, i.e.,
which corresponds to the cases (C)--(D). The generic case is (C),
in which the following proposition provides some information.

\begin{proposition}\label{p7.1}
Let $[M_1,M_2,M_3]$ be a triple of the form (\ref{f7.7})
with $M_1M_2M_3 = I$, and assume that condition (C) is fulfilled.
Let $[\wt{M}_1,\wt{M}_2,\wt{M}_3]$ be another triple
of the same form and with the same entries except that
$\gamma_1,\gamma_2,\gamma_3$ are replaced by
$\wt{\gamma}_1,\wt{\gamma_2},\wt{\gamma}_3$. Then these triples are equivalent
if and only if there exist $\alpha,\beta,\gamma\in\C$ such that
$$
\ba{r@{\quad}l@{\quad\mbox{in case }\;\;}l}
1.) &
\wt{\gamma}_k \;=\; \gamma_k
+\gamma(\mu_k^{(1)}-\mu_k^{(3)}) &
(\mu_k^{(1)})_{k=1}^3\neq(\mu_k^{(2)})_{k=1}^3\neq(\mu_k^{(3)})_{k=1}^3,
\\[.5ex]
2.) &
\wt{\gamma}_k \;=\; \gamma_k
+\gamma(\mu_k^{(2)}-\mu_k^{(3)})+\alpha\beta_k &
(\mu_k^{(1)})_{k=1}^3=(\mu_k^{(2)})_{k=1}^3\neq(\mu_k^{(3)})_{k=1}^3,
\\[.5ex]
3.) &
\wt{\gamma}_k \;=\; \gamma_k
+\gamma(\mu_k^{(1)}-\mu_k^{(2)})+\beta\alpha_k &
(\mu_k^{(1)})_{k=1}^3\neq(\mu_k^{(2)})_{k=1}^3=(\mu_k^{(3)})_{k=1}^3,
\\[.5ex]
4.) &
\wt{\gamma}_k \;=\; \gamma_k
+\alpha\beta_k+\beta\alpha_k &
(\mu_k^{(1)})_{k=1}^3=(\mu_k^{(2)})_{k=1}^3=(\mu_k^{(3)})_{k=1}^3.
\ea
$$
\end{proposition}
\begin{proof}
These triples are equivalent if and only if there exists an invertible matrix
$S$ such that $\wt{M}_k=S\iv M_kS$. It follows that the matrix $S$ maps the
invariant subspaces of $[\wt{M}_1,\wt{M}_2,\wt{M}_3]$ into the invariant
subspaces of $[M_1,M_2,M_3]$. Because of the condition (C),
$V_1=\C\oplus\{0\}\oplus\{0\}$ and $V_2=\C\oplus\C\oplus\{0\}$ are the only
non-trivial invariant subspaces of both triples. Hence $V_1$ and $V_2$ are
also invariant subspaces of $S$, and this implies that $S$ is of triangular
form,
\bqn
S &=&
\left(
\ba{ccc} s_1 & \alpha & \gamma \\ 0 & s_2 & \beta \\
0 & 0 & s_3 \ea \right).
\eqn
Now the equality $\wt{M_k}=S\iv M_k S$ gives the following two
conditions,
\bqn\label{f7.14}
\alpha_k &=&
\frac{s_2}{s_1}\alpha_k+\frac{\alpha}{s_1}(\mu_k^{(1)}-\mu_k^{(2)}), \\[.5ex]
\beta_k &=& \frac{s_3}{s_2}\beta_k+\frac{\beta}{s_2}(\mu_k^{(2)}-\mu_k^{(3)}),
\label{f7.15}
\eqn
and a third condition, which will be stated later.
We first claim that $s_1=s_2$. Indeed, (\ref{f7.14}) can be rewritten
as $0=(s_2-s_1)\alpha_k+\alpha(\mu_k^{(1)}-\mu_k^{(2)})$. If $s_1\neq s_2$,
then Proposition \ref{p6.1} applied to the triple of $2\times 2$
matrices obtained from $[M_1,M_2,M_3]$ by removing the third rows and columns
implies that $[M_1,M_2,M_3]$ is equivalent to a triple of the form
(\ref{f7.7}) with $\alpha_k=0$. This contradicts condition (C). Similarly,
we obtain $s_2=s_3$. A simple thought shows that we may assume
$s_1=s_2=s_3=1$. Hence (\ref{f7.14}), (\ref{f7.15}) simplify to
\ben\label{f7.16}
\alpha(\mu_k^{(1)}-\mu_k^{(2)})\;\;=\;\;
\beta(\mu_k^{(2)}-\mu_k^{(3)})\;\;=\;\;0,
\een
and the third condition can be expressed as
\bqn\label{f7.17}
\wt{\gamma}_k &=& \gamma_k +
\gamma(\mu_k^{(1)}-\mu_k^{(3)})-\alpha\beta_k+\beta\alpha_k.
\eqn
The conditions (\ref{f7.16}) lead to four distinct cases.
If $\mu_k^{(1)}\neq\mu_k^{(2)}$ for some $k$, then $\alpha=0$,
and if $\mu_k^{(2)}\neq\mu_k^{(3)}$ for some $k$, then $\beta=0$.
Otherwise, $\alpha$ and $\beta$, resp., can be arbitrary.
Combining these statements with formula (\ref{f7.17}) proves the
assertion.
\end{proof}

We are now able to describe the sets of equivalence classes of triples
$[M_1,M_2,M_3]$, which satisfy condition (C) and for which all values
appearing in (\ref{f7.7}) are kept fixed except $\gamma_1,\gamma_2,\gamma_3$.
We first note that $\gamma_1,\gamma_2,\gamma_3$ have to satisfy the linear
condition (\ref{f7.lin3}).

In case 1, the set of equivalence classes can be identified with $\C$
if $\mu_k^{(1)}\neq\mu_k^{(3)}$ for some $k$ and with $\C^2$ if
$\mu_k^{(1)}=\mu_k^{(3)}$ for all $k$, respectively.

In cases 2 and 3,  the sets of equivalence classes are singletons.
In fact, as to case 2, the vectors $(\beta_k)$ and
$(\mu_k^{(2)}-\mu_k^{(3)})$ are linearly independent. If they were
not, i.e., if $(\beta_k)$ is a multiple of
$(\mu_k^{(2)}-\mu_k^{(3)})$, then Proposition \ref{p6.1} applied
to the triple of $2\times 2$ matrices obtained from
$[M_1,M_2,M_3]$ by removing the first rows and columns implies
that $[M_1,M_2,M_3]$ is equivalent to a triple (\ref{f7.7}) with
$\beta_k=0$. This conflicts with condition (C).

In case 4, the set of equivalence classes is a singleton if
the vectors $(\alpha_k)$ and $(\beta_k)$ are linearly independent.
Otherwise, this set can be identified with $\C$.

In the next theorem, we consider all system of standard form
$Y'(\tz)=A(z)Y(\tz)$ with indices $\ka_1\ge\ka_2\ge\ka_3$, for which
$A(z)$ can be written as
\bqn\label{f7.10}
A(z) &=& \frac{1}{p(z)}
\left(\ba{ccc}
\ka_1z^2+p_1(z) & p_{12}(z) & p_{13}(z) \\
0 & \ka_2z^2+p_2(z) & p_{23}(z) \\
0 & 0 & \ka_3z^2+p_3(z) \ea \right),
\eqn
where $p(z)=(z-a_1)(z-a_2)(z-a_3)$,
$p_j\in\cP_1$, $p_{jk}\in\cP_{1+\ka_j-\ka_k}$. We will completely characterize
the data associated to such systems.

\begin{theorem}\label{t7.2}
Let $\ka_1,\ka_2,\ka_3\in\Z$, $\ka_1\ge\ka_2\ge\ka_3$, and
$[M_1,M_2,M_3]$ and $[E_1,E_2,E_3]$ be admissible data of
$3\times 3$ matrices. Then the following two statements are equivalent:
\begin{itemize}
\item[(i)]
$[M_1,M_2,M_3]$ and $[E_1,E_2,E_3]$ is the data associated to some system
of the form  (\ref{f7.10}).
\item[(ii)]
$[M_1,M_2,M_3]$ is equivalent to some triple
(\ref{f7.7}) with $n_j=\ka_j$ for $j=1,2,3$, where the integers $n_j$ are
defined by (\ref{f7.8}).
\end{itemize}
\end{theorem}
\begin{proof}
(i)$\Rightarrow$(ii):
Because $A(z)$ is of triangular form, there exists a
solution which is also of triangular form
(see Lemma \ref{l5.1}).
Hence we can write $Y'(\tz)=A(z)Y(\tz)$ as
\bqn\label{f7.11}
\left(\ba{ccc} y_1' & y_{12}' & y_{13}'  \\
0 & y_2'  & y_{23}'  \\
0 & 0 & Y_3'  \ea \right) &=&
\left(\ba{ccc} a_1  & a_{12}  & a_{13}  \\
0 & a_2  & a_{23}  \\
0 & 0 & a_3  \ea \right)
\left(\ba{ccc} y_1  & y_{12}  & y_{13} \\
0 & y_2  & y_{23}  \\
0 & 0 & y_3  \ea \right).
\eqn
For $j=1,2,3$, we have
\bqn\label{f7.12}
a_j(z) &=& \frac{\ka_jz^{2}+p_j(z)}{p(z)}
\;\;=\;\;
\sum_{k=1}^3 \frac{\ep_k^{(j)}}{z-a_k}
\eqn
where the numbers $\ep_1^{(j)},\ep_2^{(j)},\ep_3^{(j)}$ are defined by the
last equality and satisfy $\ka_j=\ep_1^{(j)}+\ep_2^{(j)}+\ep_3^{(j)}$.
Because $E_k$ is similar to the residue of $A(z)$ at $z=a_k$, the numbers
$\ep_k^{(1)},\ep_k^{(2)},\ep_k^{(3)}$ are the eigenvalues of $E_k$.
Once it is shown that they are properly ordered, it follows
that $\ka_j=n_j$ for $j=1,2,3$.

The diagonal entries of the solution $Y(\tz)$ satisfy
$y_j'(\tz)=a_j(z)y_j(\tz)$. Therefore they are given by
$y_j(\tz) = \prod_{k=1}^3 (\tz-a_k)^{\ep_k^{(j)}}$
up to a redundant constant.
Regardless of the precise expressions for the off-diagonal
entries of $Y(\tz)$ it follows that
the corresponding monodromy representation is
given by (\ref{f7.7}) with $\mu_k^{(j)}=\exp(-2\pi i\ep_k^{(j)})$
and with certain $\alpha_k,\beta_k,\gamma_k$. Hence $[M_1,M_2,M_3]$ is
equivalent to (\ref{f7.7}) and the above relation between
$\mu_k^{(j)}$ and $\ep_k^{(j)}$ implies that the eigenvalues of $E_k$
are properly numbered.

(ii)$\Rightarrow$(i):
Suppose we are given $[M_1,M_2,M_3]$, which is equivalent to a triple of the
form (\ref{f7.7}), and assume that the numbers $\ep_k^{(j)}$ and $n_j=\ka_j$
are defined appropriately. We introduce $a_j(z)$ by (\ref{f7.12}) and consider
systems (\ref{f7.11}) with those $a_j$ and yet unspecified
$a_{12},a_{13},a_{23}$. We have to show that there exist
$a_{12},a_{13},a_{23}$ (being the ratio $p_{jk}/p$ of polynomials with
properties as described above) such that the monodromy of this system is given
by $[M_1,M_2,M_3]$.

We first consider some trivial cases. Assume that $[M_1,M_2,M_3]$
is equivalent to (\ref{f7.7}) with $\alpha_k=0$. The proof of Theorem
\ref{t6.2} tells us that there exists an $a_{13}$ such that the system
\bqn\label{f7.21}
\left(\ba{cc} y_1' & y_{13}' \\ 0 & y_3' \ea\right) &=&
\left(\ba{cc} a_1 & a_{13} \\ 0 & a_3 \ea\right)
\left(\ba{cc} y_1 & y_{13} \\ 0 & y_3 \ea\right)
\eqn
has monodromy given by
\ben\label{f7.22}
\left[
\left(\ba{cc} \mu_1^{(1)} & \gamma_1 \\ 0 & \mu_1^{(3)} \ea\right),
\left(\ba{cc} \mu_2^{(1)} & \gamma_2 \\ 0 & \mu_2^{(3)} \ea\right),
\left(\ba{cc} \mu_3^{(1)} & \gamma_3 \\ 0 & \mu_3^{(3)} \ea\right)
\right].
\een
Moreover, there exists an $a_{23}$ such that the system
\bqn\label{f7.23}
\left(\ba{cc} y_2' & y_{23}' \\ 0 & y_3' \ea\right) &=&
\left(\ba{cc} a_2 & a_{23} \\ 0 & a_3 \ea\right)
\left(\ba{cc} y_2 & y_{23} \\ 0 & y_3 \ea\right)
\eqn
has monodromy given by
\ben\label{f7.24}
\left[
\left(\ba{cc} \mu_1^{(2)} & \beta_1 \\ 0 & \mu_1^{(3)} \ea\right),
\left(\ba{cc} \mu_2^{(2)} & \beta_2 \\ 0 & \mu_2^{(3)} \ea\right),
\left(\ba{cc} \mu_3^{(2)} & \beta_3 \\ 0 & \mu_3^{(3)} \ea\right)
\right].
\een
In this connection we remark that (\ref{f7.22}) and (\ref{f7.24}) are indeed
``admissible'' triples in the sense that the product of the occurring matrices
is equal to $I$. Moreover, $a_{13}=p_{13}/p$ and $a_{23}=p_{23}/p$ where the
polynomials $p_{13}$ and $p_{23}$ have the proper degree. Now we consider the
system (\ref{f7.11}) with those $a_{13}$ and $a_{23}$ and with $a_{12}=0$. This
system has a solution of triangular form as given in (\ref{f7.11}) with
$y_{12}=0$. The monodromy representation with respect to this solution is of
the form (\ref{f7.7}) with desired values $\mu_k^{(j)}$ and $\alpha_k=0$.
However, we may have different values $\wt{\beta}_k$ and $\wt{\gamma}_k$
(instead of the desired values $\beta_k$ and $\gamma_k$). Analyzing the
relation between the monodromy of (\ref{f7.11}) and the monodromy of the
``subsystems'' (\ref{f7.21}) and (\ref{f7.23}), it follows that (\ref{f7.22})
is equivalent to (\ref{f7.22}) with $\gamma_k$ replaced by $\wt{\gamma}_k$ and
that (\ref{f7.24}) is equivalent to (\ref{f7.24}) with $\beta_k$ replaced by
$\wt{\beta}_k$. We apply Proposition \ref{p6.1} and obtain that
\ben
\wt{\gamma}_k=\lambda^{(1)}\gamma_k+\rho^{(1)}(\mu_k^{(1)}-\mu_k^{(3)}),
\qquad
\wt{\beta}_k=\lambda^{(2)}\beta_k+\rho^{(2)}(\mu_k^{(2)}-\mu_k^{(3)})
\een
with certain $\lambda^{(1)},\lambda^{(2)}\in\C\setminus\{0\}$ and
$\rho^{(1)},\rho^{(2)}\in\C$.
A simple thought shows that a representation (\ref{f7.7}) with $\alpha_k=0$
and having above $\beta_k$ and $\gamma_k$ as entries is
equivalent to the same representation but with above $\wt{\beta}_k$ and
$\wt{\gamma}_k$. In fact, the equivalence is established by a certain upper
triangular matrix having $(1,2)$-entry equal to zero.
Hence, the monodromy of the constructed system (\ref{f7.11}) is given by
$[M_1,M_2,M_3]$, and this concludes this case.

The other trivial case, where $[M_1,M_2,M_3]$ is equivalent to a triple of the
form (\ref{f7.7}) with $\beta_k=0$ can be treated in the same way.
Here we consider the system (\ref{f7.11}) with $a_{23}=0$ and $Y_{23}=0$ and
the appropriate ``subsystems''.

Now assume that $[M_1,M_2,M_3]$ does not fall into these trivial classes.
In other words, $[M_1,M_2,M_3]$ satisfies the condition (C) concerning the
reducibility type. Arguing similar as above, we find an appropriate function
$a_{12}$ such that the system
\bqn\label{f7.26}
\left(\ba{cc} y_1' & y_{12}' \\ 0 & y_2' \ea\right) &=&
\left(\ba{cc} a_1 & a_{12} \\ 0 & a_2 \ea\right)
\left(\ba{cc} y_1 & y_{12} \\ 0 & y_2 \ea\right)
\eqn
has monodromy given by
\ben\label{f7.27}
\left[
\left(\ba{cc} \mu_1^{(1)} & \alpha_1 \\ 0 & \mu_1^{(2)} \ea\right),
\left(\ba{cc} \mu_2^{(1)} & \alpha_2 \\ 0 & \mu_2^{(2)} \ea\right),
\left(\ba{cc} \mu_3^{(1)} & \alpha_3 \\ 0 & \mu_3^{(2)} \ea\right)
\right].
\een
and we also find a function $a_{23}$ such that the system (\ref{f7.23})
has the monodromy (\ref{f7.24}).
With those functions $a_{12}$ and $a_{23}$ and an arbitrary function
$a_{13}^{\#}$ we consider the system (\ref{f7.11}).
This system has a solution with certain functions $y_1,y_2,y_3$ and
$y_{12},y_{23},y_{13}^{\#}$. The monodromy representation of this solution
is given by (\ref{f7.7}) with the desired values $\mu_k^{(j)}$, but possibly
different values $\wt{\alpha}_k$, $\wt{\beta}_k$ and $\wt{\gamma}^{\#}_k$.
Analyzing the relation to the ``subsystems'', for which we know the monodromy,
it follows again from Proposition \ref{p6.1} that
\ben
\wt{\alpha}_k=\lambda^{(1)}\alpha_k+\rho^{(1)}(\mu_k^{(1)}-\mu_k^{(2)}),
\qquad
\wt{\beta}_k=\lambda^{(2)}\beta_k+\rho^{(2)}(\mu_k^{(2)}-\mu_k^{(3)})
\een
with certain $\lambda^{(1)},\lambda^{(2)}\in\C\setminus\{0\}$ and
$\rho^{(1)},\rho^{(2)}\in\C$.
By means of the construction of a suitable triangular matrix (as above), we
obtain that $[M_1,M_2,M_3]$ is equivalent to a triple (\ref{f7.7}) having
off-diagonal entries $\wt{\alpha}_k$, $\wt{\beta}_k$ and certain
$\wt{\gamma}_k$.

In other word, the monodromy of the system we have constructed so
far is ``almost'' as desired, only the values $\wt{\gamma}_k^{\#}$
might not yet be the desired ones $\wt{\gamma}_k$. In fact, in
cases 2 and 3 and also in the ``generic'' subcase of case 4, which
has been singled out in Proposition \ref{p7.1}, we are already
done because (as has been remarked in the paragraphs following
this proposition) there exists only a single equivalence class of
such triples. Hence the monodromy of the constructed system is
equivalent to $[M_1,M_2,M_3]$.

In the remaining cases we must still modify the system, and we proceed as
follows. Again the proof of Theorem \ref{t6.2}
tells us that there exists a suitable function $a_{13}^{\$}$ such that the
system
\bqn\label{f7.29}
\left(\ba{cc} y_1' & {y_{13}^{\$}}' \\ 0 & y_3' \ea\right) &=&
\left(\ba{cc} a_1 & a_{13}^{\$} \\ 0 & a_3 \ea\right)
\left(\ba{cc} y_1 & y_{13}^{\$} \\ 0 & y_3 \ea\right)
\eqn
has monodromy given by
\ben\label{f7.30}
\left[
\left(\ba{cc} \mu_1^{(1)} & \wt{\gamma}_1-\wt{\gamma}_1^{\#} \\ 0 &
\mu_1^{(3)} \ea\right),
\left(\ba{cc} \mu_2^{(1)} & \wt{\gamma}_2-\wt{\gamma}_2^{\#}\\ 0 &
\mu_2^{(3)} \ea\right),
\left(\ba{cc} \mu_3^{(1)} & \wt{\gamma}_3-\wt{\gamma}_3^{\#} \\ 0 &
\mu_3^{(3)} \ea\right) \right].
\een
In this connection we remark that both $\wt{\gamma}_k$ and
$\wt{\gamma}_k^{\#}$ satisfy a linear relation (\ref{f7.lin3})
with the same right hand side. The differences satisfy the
homogeneous linear relation, which corresponds to (\ref{f6.2}).
Hence the (\ref{f7.30}) is an admissible triple, i.e.,
the product of the matrices is equal to one.

The monodromy representation of the solution in (\ref{f7.29})
is also given by (\ref{f7.30}) except that the upper right entries are
replaced by certain values $\wt{\gamma}_k^{\$}$. The relation between these
and the former values is (see Proposition \ref{p6.1})
\bqn\label{f7.31}
\wt{\gamma}_k^{\$} &=& \lambda(\wt{\gamma}_k-\wt{\gamma}_k^{\#})+
\rho(\mu_k^{(1)}-\mu_k^{(3)})
\eqn
with certain $\lambda\in\C\setminus\{0\}$ and $\rho\in\C$.

Now we consider a new system (\ref{f7.11}) where $A(z)$ has the same entries
as the old one except that $a_{13}^{\#}$ is replaced by $a_{13}$
where we define $a_{13}:=a_{13}^{\#}+\lambda\iv a_{13}^{\$}$. It is easy to
see that the solution of this new system has the same entries as the old
solution except that $y_{13}^{\#}$ is replaced by a certain $y_{13}$.
Moreover, the monodromy representation with respect to this solution is of the
form (\ref{f7.7}) with off-diagonal entries $\wt{\alpha}_k$, $\wt{\beta}_k$
(as in the old one) and certain $\wt{\gamma}_k^{\&}$.

What remains to show is that this monodromy representation is
equivalent to $[M_1,M_2,M_3]$. Because we already pointed out that
$[M_1,M_2,M_3]$ is equivalent to (\ref{f7.7}) with off-diagonal
entries $\wt{\alpha}_k$, $\wt{\beta}_k$ and $\wt{\gamma}_k$, we
are done as soon as we have shown that $\wt{\gamma}_k$ and
$\wt{\gamma}_k^{\&}$ are related in the way described in
Proposition \ref{p7.1}.

In order to show this, we have to consider the defining relations
for the $(1,3)$-entries of the solutions and of the monodromy
representations explicitly. As to the old system, we have
\bqn\label{f7.32} {y_{13}^{\#}}' &=&
a_1y_{13}^{\#}+a_{12}y_{23}+a_{13}^{\#}y_3,\\[.5ex] \label{f7.33}
y_{13}^{\#}(\sigma_k\iv(\tz)) &=&
y_1(\tz)\wt{\gamma}_k^{\#}+y_{12}(\tz)\wt{\beta}_k+
y_{13}^{\#}(\tz)\mu_k^{(3)}, \eqn where $\sigma_k$, $k=1,2,3$, are
the deck transformations. The corresponding relations for the new
system are \bqn\label{f7.34} y_{13}' &=&
a_1y_{13}+a_{12}y_{23}+a_{13}y_3,\\[.5ex] \label{f7.35}
y_{13}(\sigma_k\iv(\tz)) &=&
y_1(\tz)\wt{\gamma}_k^{\&}+y_{12}(\tz)\wt{\beta}_k+
y_{13}(\tz)\mu_k^{(3)}. \eqn We also have to state the conditions
on the $(1,2)$-entries in (\ref{f7.29}) and (\ref{f7.30}),
\bqn\label{f7.36} {y_{13}^{\$}}' &=&
a_1y_{13}^{\$}+a_{13}^{\$}y_3,\\[.5ex] \label{f7.37}
y_{13}^{\$}(\sigma_k\iv(\tz)) &=& y_1(\tz)\wt{\gamma}_k^{\$}+
y_{13}^{\$}(\tz)\mu_k^{(3)}. \eqn Combining (\ref{f7.32}) and
(\ref{f7.36}) and using $a_{13}=a_{13}^{\#}+\lambda\iv
a_{13}^{\$}$, we obtain that \bqn (y_{13}^{\#}+\lambda\iv
y_{13}^{\$})' &=& a_1(y_{13}^{\#}+\lambda\iv
y_{13}^{\$})+a_{12}y_{23}+a_{13}y_3. \eqn Comparing this with
(\ref{f7.34}) we may conclude without loss of generality that
$y_{13}=y_{13}^{\#}+\lambda\iv y_{13}^{\$}$. (The solution of the
differential equation is not unique but we may take any of them.)
Using this identity and combining (\ref{f7.33}) with
(\ref{f7.37}), we arrive at \bqn y_{13}(\sigma_k\iv(\tz)) &=&
y_1(\tz)(\wt{\gamma}_k^{\#}+\lambda\iv \wt{\gamma}_k^{\$}) +
y_{12}(\tz)\wt{\beta}_k+ y_{13}(\tz)\mu_k^{(3)}. \eqn Comparing
this with (\ref{f7.35}) we conclude that \bqn \wt{\gamma}_k^{\&}
&=& \wt{\gamma}_k^{\#}+\lambda\iv \wt{\gamma}_k^{\$}\nn\\ &=&
\wt{\gamma}_k+\rho\lambda\iv (\mu_k^{(1)}-\mu_k^{(3)}), \eqn where
the last equality follows from (\ref{f7.31}). Now we conclude from
Proposition \ref{p7.1} that the triple (\ref{f7.7}) with
$\wt{\alpha}_k$, $\wt{\beta_k}$ and $\wt{\gamma}_k^{\&}$, which is
the monodromy representation for the new system, is equivalent to
the triple (\ref{f7.7}) with $\wt{\alpha}_k$, $\wt{\beta_k}$ and
$\wt{\gamma}_k$, which in turn is equivalent to $[M_1,M_2,M_3]$.
Hence the monodromy of the new system is given by $[M_1,M_2,M_3]$.
\end{proof}

\subsection{Systems of block-triangular form}

In what follows, we are going to consider $3\times 3$ systems of
certain block-triangular forms, which are not of triangular form. There
exists two types of them relating to two kinds of block-triangular partitions
of the $3\times 3$ matrix function $A(z)$.

In order to describe these systems and the data associated to
them, we have first to introduce some definitions and to state
some results related to $2\times 2$ systems and the corresponding
data.

\pagebreak[1]
Let $[M_1^{\#},M_2^{\#},M_3^{\#}]$ and $[E_1^{\#},E_2^{\#},E_3^{\#}]$ be admissible data of
$2\times 2$ matrices. We call this data {\em quasi-block data} if either
\begin{itemize}
\item[(I)]
$[M_1^{\#},M_2^{\#},M_3^{\#}]$ is irreducible, or,
\item[(II)]
$[M_1^{\#},M_2^{\#},M_3^{\#}]$ possesses exactly one non-trivial invariant subspace
and  we have $n_1^{\#}<n_2^{\#}$, where $n_1^{\#}$ and $n_2^{\#}$ are defined
by (\ref{f6.int}).
\end{itemize}

The next result is an immediate consequence of Theorem \ref{t6.2}.
It describes the data of $2\times 2$ systems with three singularities which are
not equivalent to a system of triangular form. We remark that such systems have
necessarily indices $\ka_1^{\#},\ka_2^{\#}\in\Z$ for which
$\ka_1^{\#}-\ka_2^{\#}\in\{0,1\}$.

\begin{lemma}\label{l7.3}
Let $[M_1^{\#},M_2^{\#},M_3^{\#}]$ and $[E_1^{\#},E_2^{\#},E_3^{\#}]$
be admissible data of $2\times 2$ matrices. Then the following two statements
are equivalent:
\begin{itemize}
\item[(i)] $[M_1^{\#},M_2^{\#},M_3^{\#}]$ and $[E_1^{\#},E_2^{\#},E_3^{\#}]$
is the data of a $2\times 2$ system of standard form with three singularities
which is not equivalent to a system of the form (\ref{f6.6}).
\item[(ii)]
$[M_1^{\#},M_2^{\#},M_3^{\#}]$ and $[E_1^{\#},E_2^{\#},E_3^{\#}]$
is quasi-block data.
\end{itemize}
\end{lemma}
\begin{proof}
Quasi-block data is exactly the data which is not covered by (ii)
of Theorem \ref{t6.2}. As was already remarked in the paragraph
after Theorem \ref{t6.2}, the data which is covered by this
theorem is exactly the data for which
$[M_1^{\#},M_2^{\#},M_3^{\#}]$ is reducible with the condition
$n_1^{\#}\geq n_2^{\#}$ or is diagonalizable without another
condition.
\end{proof}

Let $[M_1,M_2,M_3]$ and $[E_1,E_2,E_3]$ be admissible data of
$3\times 3$ matrices. We call this data {\em $(1,2)$-quasi-block
data} if $[M_1,M_2,M_3]$ is equivalent to a triple
\ben\label{f7.41} \left[ \left(\ba{cc}\mu_1 & \alpha_1 \\ 0 &
M_1^{\#} \ea \right), \left(\ba{cc}\mu_2 & \alpha_2 \\ 0 &
M_2^{\#} \ea \right), \left(\ba{cc}\mu_3 & \alpha_3 \\ 0 &
M_3^{\#} \ea \right) \right] \een and
$[M_1^{\#},M_2^{\#},M_3^{\#}]$ and $[E_1^{\#},E_2^{\#},E_3^{\#}]$
is quasi-block data. Here $\mu_k$ are numbers, $\alpha_k$ are
$1\times 2$ vectors, and $M_k^{\#}$ are $2\times 2$ matrices. The
matrices $E_k^{\#}$ and numbers $\ep_k$ are supposed to satisfy
$M_k^{\#}\sim\exp(-2\pi iE_k^{\#})$ and $\mu_k=\exp(-2\pi i
\ep_k)$ and are chosen in such a way that the eigenvalues of
$E_k^{\#}$ and the number $\ep_k$ are exactly the eigenvalues of
$E_k$. We note that for given $M_k^{\#}$ and $\mu_k$, the number
$\ep_k$ is uniquely determined and so is the matrix $E_k^{\#}$ up
to similarity (because of the non-resonance of $E_k$). We also
define the integers \bqn\label{f7.42} \nu \;\;=\;\;
\ep_1+\ep_2+\ep_3,\qquad N \;\;=\;\; {\rm
trace\,}(E_1^{\#}+E_2^{\#}+E_3^{\#}). \eqn We remark that if
$[M_1,M_2,M_3]$ possesses two different one-dimensional invariant
subspaces, then $[M_1,M_2,M_3]$ is equivalent to different triples
of the form (\ref{f7.41}) (with different $\mu_k$ and $M_k^{\#}$
in particular). In this case we obtain in general also different
values for $\ep_k, E_k^{\#}, \nu$ and $N$.

\begin{proposition}\label{p7.4}
Let $[M_1,M_2,M_3]$ be a triple of the form (\ref{f7.41}). Assume that
\begin{itemize}
\item[(a)]
$[M_1^{\#},M_2^{\#},M_3^{\#}]$ is not
diagonalizable;
\item[(b)]
If  $[M_1^{\#},M_2^{\#},M_3^{\#}]$ possess exactly one non-trivial invariant
subspace, then
\begin{itemize}
\item[(b1)]
for at least one $k$, $M_k^\#$ possesses two different eigenvalues;
\item[(b2)]
for at least one $k$, $\mu_k$ is not equal to the eigenvalue of $M_k^\#$
which corresponds to this invariant subspace.
\end{itemize}
\end{itemize}
Let $[\wt{M}_1,\wt{M}_2,\wt{M}_3]$ be another triple of the same form
(\ref{f7.41}) but with $\alpha_1,\alpha_2,\alpha_3$ replaced by
$\wt{\alpha}_1,\wt{\alpha}_2,\wt{\alpha}_3$. Then these triples are equivalent
if and only if there exist $\lambda\in\C\setminus\{0\}$ and
$c\in\C^{1\times 2}$ such that
\bqn\label{f7.43}
\alpha_k &=& \lambda\wt{\alpha}_k+cM_k^{\#}-\mu_k c
\qquad\quad\mbox{ for each }k=1,2,3.
\eqn
\end{proposition}
\begin{proof}
The triples $[M_1,M_2,M_3]$ and $[\wt{M}_1,\wt{M}_2,\wt{M}_3]$
are equivalent if and only if there exists a matrix $C\in G\C^{3\times 3}$
such that $\wt{M}_k=CM_kC\iv$. We write $C$ in the same block form as
the matrices in the triple (\ref{f7.41}),
\bqn\label{f7.44}
C &=& \left(\ba{cc} c_{11} & c_{12} \\
c_{21} & c_{22} \ea\right),
\eqn
and inspect the $(2,1)$-block entry in the resulting equality $\wt{M}_kC=CM_k$.
We obtain $M_k^{\#}c_{21}=c_{21}\mu_k$. Because of (a) and (b2), it follows
that $c_{21}=0$. Otherwise, ${\rm Span}\,\{c_{21}\}$ is an invariant subspace
for all $M_k^\#$ with the eigenvalue $\mu_k$.
Hence $C$ is also of block triangular form.

Now we look at the $(2,2)$-block entry of $\wt{M}_kC=CM_k$, and it
follows that $M_k^{\#}c_{22}=c_{22}M_k^{\#}$. We claim that
$c_{22}$ is a scalar matrix. Assume the contrary. Then $c_{22}$ is
either similar to a diagonal matrix with two different diagonal
entries or it is similar to a Jordan block. We may assume without
loss of generality that $c_{22}$ is actually equal to such a
diagonal matrix or to a Jordan block. (Otherwise one applies an
appropriate similarity transformation simultaneously to $c_{22}$
and $M_1^{\#},M_2^{\#},M_3^{\#}$, which will equally result in a
contradiction.)  By simple computations it follows in the first
case that all $M_1^{\#},M_2^{\#},M_3^{\#}$ are of diagonal form.
In the second case, it follows that all
$M_1^{\#},M_2^{\#},M_3^{\#}$ are either Jordan blocks or scalar
matrices. Both contradict the assumption on the triple
$[M_1^{\#},M_2^{\#},M_3^{\#}]$.

With this information about the entries of $C$, it is now easy to conclude
that $\wt{M}_k=CM_kC\iv$ is equivalent to (\ref{f7.43}).
\end{proof}

Now we are able to identify the set of equivalence classes for
triples of the form (\ref{f7.41}) with $\mu_k$ and $M_k^{\#}$
being fixed and satisfying the assumptions stated in the previous
proposition. All such triples are parameterized by $1\times 2$
vectors $\alpha_1,\alpha_2,\alpha_3$. It is convenient to
introduce the $1\times 6$ vector
$\alpha=(\alpha_1,\alpha_2,\alpha_3)$. First we have to observe
the linear condition (\ref{f7.3}) which can be restated as
\bqn\label{f7.45} \alpha \left(\ba{c} M_2^\# M_3^\# \\ \mu_1
M_3^\# \\ \mu_1\mu_2I_2 \ea \right) =0, \eqn where $I_2$ is the
$2\times 2$ identity matrix. Because the $6\times 2$ matrix
appearing in this linear condition has rank $2$, the vector
$\alpha$ is actually subject to two linear conditions. Hence all
``admissible'' vectors $\alpha$ are taken from a $4$-dimensional
linear subspace. The equivalence relation (\ref{f7.43}) can
shortly be written as \bqn\label{f7.46} \alpha &=&
\lambda\wt{\alpha} \;+\;c\,\Big(\ba{ccc} M_1^\#-\mu_1I_2, &
M_2^\#-\mu_2I_2, & M_3^\#-\mu_3I_2 \ea\Big) \eqn with
$\lambda\in\C\setminus\{0\}$ and $c\in\C^{1\times 2}$. Because of
the assumptions in Proposition \ref{p7.4}, the $2\times 6$ matrix
appearing there has either rank $2$ or rank $1$. It has rank $1$
if and only if $[M_1^\#,M_2^\#,M_3^\#]$ is equivalent to a triple
\ben\label{f7.46a} \left[ \left(\ba{cc} \hat{\mu}_1
&\hat{\alpha}_1 \\ 0 & \mu_1 \ea\right), \left(\ba{cc} \hat{\mu}_2
&\hat{\alpha}_1 \\ 0 & \mu_2 \ea\right), \left(\ba{cc} \hat{\mu}_3
&\hat{\alpha}_1 \\ 0 & \mu_3 \ea\right) \right] \een with certain
$\hat{\mu}_k$ and $\hat{\alpha}_k$. In other words, the above
matrix has rank $1$ if and only if $[M_1^\#,M_2^\#,M_3^\#]$
possesses exactly one non-trivial invariant subspace and, for each
$k=1,2,3$, the number $\mu_k$ is equal to the eigenvalue of
$M_k^\#$ which does not correspond to this invariant subspace. The
second part in (\ref{f7.46}) gives a ``linear'' factorization of
vector spaces and thus yields either a two-dimensional or a
three-dimensional quotient space. Finally, there is another
``multiplicative'' equivalence. So we arrive at the statement that
the set of equivalence classes of such triples $[M_1,M_2,M_3]$ can
be identified either with $\sP^1\cup\{0\}$ or with
$\sP^2\cup\{0\}$ (corresponding to the above distinction), where
$\sP^n$ is the $n$-dimensional complex projective space and
$\{0\}$ denotes the single equivalence class containing the vector
$\alpha=0$. This single equivalence class consists precisely of
the  triples which are ``block-diagonalizable''.

Next we consider systems $Y'(\tz)=A(z)Y(\tz)$ of standard form with indices
$\ka_1\ge\ka_2\ge\ka_3$ for which $A(z)$ can be written as
\bqn\label{f7.47}
A(z) &=& \frac{1}{p(z)}
\left(\ba{ccc}
\ka_1z^2+p_1(z) & p_{12}(z) & p_{13}(z) \\
0 & \ka_2z^2+p_2(z) & p_{23}(z) \\
0 & p_{32}(z) & \ka_3z^2 + p_3(z) \ea\right),
\eqn
where $p(z)=(z-a_1)(z-a_2)(z-a_3)$, $p_k\in\cP_1$,
$p_{jk}\in\cP_{1+\ka_j-\ka_k}$. If such a system is not equivalent to a system
of the form (\ref{f7.10}), then it will be called a {\em system of
$(1,2)$-block form}.
We remark that necessarily $\ka_2-\ka_3\leq 1$ because otherwise
$p_{32}(z)=0$ and the system is of the form (\ref{f7.10}).

The next theorem describes the data which is associated to such
systems.

\begin{theorem}\label{t7.5}
Let $\ka_1,\ka_2,\ka_3\in\Z$, $\ka_1\ge\ka_2\ge\ka_3$, and assume
$\ka_2-\ka_3\leq 1$. Let $[M_1,M_2,M_3]$ and $[E_1,E_2,E_3]$ be admissible
data of $3\times 3$ matrices. Then the following two statements are
equivalent:
\begin{itemize}
\item[(i)]
$[M_1,M_2,M_3]$ and $[E_1,E_2,E_3]$ is the data associated to a system of
$(1,2)$-block form with indices $\ka_1,\ka_2,\ka_3$.
\item[(ii)]
$[M_1,M_2,M_3]$ and $[E_1,E_2,E_3]$ is $(1,2)$-quasi-block data, and
$\nu=\ka_1$, $N=\ka_2+\ka_3$, where $\nu$ and $N$ are defined by (\ref{f7.42}).
\end{itemize}
\end{theorem}
\begin{proof}
(i)$\Rightarrow$(ii):
Because $A(z)$ is of block-triangular form, there exists a solution which is
also of block triangular form. Hence let us write
$Y'(\tz)=A(z)=Y(\tz)$ as
\bqn\label{f7.48}
\left(\ba{cc} y_1'(\tz) & y_{12}'(\tz) \\
0 & Y_{2}'(\tz) \ea\right) &=&
\left(\ba{cc} a_1(z) & a_{12}(z) \\
0 & A_{2}(z) \ea\right)
\left(\ba{cc} y_1(\tz) & y_{12}(\tz) \\
0 & Y_{2}(\tz) \ea\right).
\eqn

We first consider the scalar subsystem $y_1'(\tz)=a_1(z)y_1(\tz)$.
We know that $a_1(z)=(\ka_{1}z^2+p_1(z))/p(z)$.
If $\ep_1,\ep_2,\ep_3$ denote the residues of $a_1(z)$ at $z=a_1,a_2,a_3$,
then the monodromy is given by $[\mu_1,\mu_2,\mu_3]$, where
$\mu_k=\exp(-2\pi i\ep_k)$. Moreover, it follows that
$\ka_1=\ep_1+\ep_2+\ep_3$.

Now we consider the $2\times 2$ subsystem $Y_2'(\tz)=A_2(z)Y_2(\tz)$,
where $A_2(z)$ is given by the lower right block in (\ref{f7.47}).
We claim that this system is not equivalent to a system of the form
(\ref{f6.6}). Indeed, if it were equivalent, then a simple thought shows that
the system $Y'(\tz)=A(z)Y(\tz)$ is equivalent to a system of the form
(\ref{f7.10}), which contradicts the assumption that it is a system
of $(1,2)$-block form. Hence we are in a position to apply
Lemma \ref{l7.3}. We obtain that the data of the
system $Y_2'(\tz)=A_2(z)Y_2(\tz)$, $[M_1^\#,M_2^\#,M_3^\#]$ and
$[E_1^\#,E_2^\#,E_3^\#]$, is quasi-block data.
Note that $M_k^\#\sim \exp(-2\pi iE_k^\#)$ and that $E_k^\#$ is similar to the
residue of $A_2(z)$ at $z=a_k$. Hence, in particular, $\ka_2+\ka_3={\rm
trace\,}(E_1^\#+E_2^\#+E_3^\#)$.

Combining the information about the two subsystems, we obtain
immediately that the monodromy data $[M_1,M_2,M_3]$ of the
original system is equivalent to a triple (\ref{f7.41}) with
certain $\alpha_1,\alpha_2,\alpha_3 \in\C^{1\times 2}$. The
eigenvalues of the matrix $E_k$, which is similar to the residue
of $A(z)$ at $z=a_k$, are given by the number $\ep_k$ and the
eigenvalues of $E_k^\#$ because they are the residues of $a_1(z)$
and $A_2(z)$, respectively.  Hence, since the data of the $2\times
2$ subsystem is quasi-block data, the data of the original system
is $(1,2)$-quasi-block data. Moreover, from the values for $\ka_1$
and $\ka_2+\ka_3$ given above, it follows that $\ka_1=\nu$ and
$\ka_2+\ka_3=N$.

(ii)$\Rightarrow$(i): Assume that $[M_1,M_2,M_3]$ is a triple of
the form (\ref{f7.41}) such that (on defining $\ep_k$ and $E_k^\#$
as well as $\nu$ and $N$ by (\ref{f7.42})) the appropriate
conditions for $(1,2)$-quasi-block data are satisfied. In
particular, $[M_1^\#,M_2^\#,M_3^\#]$ and $[E_1^\#,E_2^\#,E_3^\#]$
is quasi-block data, and $\nu=\ka_1$, $N=\ka_2+\ka_3$. We have to
show that there exists a system of $(1,2)$-block form which has
monodromy $[M_1,M_2,M_3]$, the indices $\ka_1,\ka_2,\ka_3$, and
the residues of $A(z)$ at $z=a_k$ are similar to $E_k$.

We first define
\bqn
a_1(z) &=& \sum_{k=1}^3 \frac{\ep_k}{z-a_k} \;\;=\;\;
\frac{\ka_1z^2+p_1(z)}{p(z)},
\eqn
noting that $\ka_1=\nu=\ep_1+\ep_2+\ep_3$. Furthermore, because
$[M_1^\#,M_2^\#,M_3^\#]$ and $[E_1^\#,E_2^\#,E_3^\#]$ is quasi-block data,
Lemma \ref{l7.3} shows that there exists a $2\times 2$ system having this
data with singularities $a_1,a_2,a_3$ such that this system is not
equivalent to a system of the form (\ref{f6.6}). Denote this system by
$Y_2'(\tz)=A_2(z)Y_2(\tz)$, and thus define the matrix function $A_2(z)$. Let
$\ka_2^\#$ and $\ka_3^\#$ be the indices of this $2\times 2$ system. We claim
that $\ka_2=\ka_2^\#$ and $\ka_3=\ka_3^\#$. Indeed, because the $2\times 2$
system is not of triangular form, we have $0\le\ka_2^\#-\ka_3^\#\le1$, and by
assumption we also have $0\le\ka_2-\ka_3\le1$. Moreover,
$\ka_2+\ka_3=N={\rm trace\,}(E_1^\#+E_2^\#+E_3^\#)=\ka_2^\#+\ka_3^\#$,
where the last relation follows from formula (\ref{f4.1})
applied in the case of the $2\times 2$ system. These relations for the
integers prove the stated equality.

Now we consider systems $Y'(\tz)=A(z)Y(\tz)$ given by (\ref{f7.48}) with those
$a_1(z)$ and $A_2(z)$ and yet unspecified $a_{12}(z)$, the entries of which
have, of course, to be the ratio of appropriate polynomials. We note that
these systems are of the form (\ref{f7.47}). They are also systems of standard
form with indices $\ka_1,\ka_2,\ka_3$, because, in particular,
$\ka_1\ge\ka_2\ge\ka_3$ is assumed. The monodromy representation of a suitable
solution of these systems is given by (\ref{f7.41}) with certain
$\alpha_1,\alpha_2,\alpha_3$. The residue $\wh{E}_k$ of $A(z)$ at $z=a_k$ has
the eigenvalues given by $\ep_k$ and the eigenvalues of $E_k^\#$. Hence, as
the eigenvalues of $E_k$ are also given in this way, it follows that $E_k$ is
similar to $\wh{E}_k$. (In case of coinciding eigenvalues, one has also to
employ that $M_k\sim\exp(-2\pi iE_k)$ and $M_k\sim\exp(-2\pi i\wh{E}_k)$.)

What remains to show is, firstly, that the systems constructed in this way
are indeed of $(1,2)$-block form, i.e., that they are not equivalent to
systems (\ref{f7.10}). Secondly, one has to show that if $a_{12}(z)$ runs
through all admissible functions, then corresponding vectors
$\alpha_1,\alpha_2,\alpha_3$ appearing in the monodromy representation
(\ref{f7.41}) run through all possibilities modulo the equivalence
described in Proposition \ref{p7.4}.

Let us turn to the first problem. Assume the contrary, namely that
our system with $A(z)$ is equivalent to a system with $\wt{A}(z)$
being of the form (\ref{f7.10}). Theorem \ref{t3.5} says that both
systems have the same indices and that the equivalence is
established by a matrix function $V(z)$, which has a particular
structure. We also know that the $2\times 2$ subsystem $A_2(z)$ of
$A(z)$ is not equivalent to a triangular system, and we are going
to show a contradiction to this assertion. If $\ka_2>\ka_3$, then
the $(3,1)$- and $(3,2)$-entries of $V(z)$ vanish. Using
$A(z)=V\iv(z)\wt{A}(z)V(z)-V\iv(z)V'(z)$, it follows that $A(z)$
has $(3,2)$-entry equal to zero, which is a contradiction. If
$\ka_1>\ka_2$, then the $V(z)$ is of the same block structure as
$A(z)$. Now using $\wt{A}(z)=V(z)A(z)V\iv(z)+V'(z)V\iv(z)$, it is
easy to see that the lower right $2\times 2$ submatrix in $V(z)$
establishes an equivalence between the subsystem $A_2(z)$ and a
triangular $2\times 2$ system, namely the one which is given by
the lower right $2\times 2$ matrix in $\wt{A}(z)$. In the
remaining case, $\ka_1=\ka_2=\ka_3$, the function $V(z)$ is a
constant matrix $V$. As before, let $\wh{E}_k$ be the residues of
$A(z)$ at $z=a_k$, and let $\wt{E}_k$ be those of $\wt{A}(z)$. The
equivalence between $A(z)$ and $\wt{A}(z)$ can now be restated as
$\wt{E}_k=V\wh{E}_kV\iv$. We claim that each one-dimensional
invariant subspace of $[\wt{E}_1,\wt{E}_2,\wt{E}_3]$ is contained
in some two-dimensional invariant subspace of this triple. Indeed,
this is obvious if this triple has at least two one-dimensional
invariant subspaces. If it has exactly one, then it is necessarily
the space $\C\oplus\{0\}\oplus\{0\}$, which is contained in
$\C\oplus\C\oplus\{0\}$ (because of the triangular form of
$\wt{E}_k$). The above relation between $\wt{E}_k$ and $\wh{E}_k$
establishes a one-to-one correspondence between invariant
subspaces of the previous triple and
$[\wh{E}_1,\wh{E}_2,\wh{E}_3]$. Hence each one-dimensional
invariant subspace of the latter triple 
is contained in some two-dimensional
invariant subspace. A simple thought shows that the $2\times 2$
matrices obtained from $[\wh{E}_1,\wh{E}_2,\wh{E}_3]$ by removing
the first rows and columns are reducible. But those matrices are
the residues of the system $A_2(z)$, and thus we obtain a
contradiction.

Finally, we turn to the problem of showing that if $a_{12}$ runs
through all admissible vector functions, then the vector
$\alpha=(\alpha_1,\alpha_2,\alpha_3)\in\C^{1\times 6}$ runs
through all possibilities modulo the equivalence stated in
Proposition \ref{p7.4} (see also (\ref{f7.45}) and (\ref{f7.46})).
Before we have to show that the assumptions of Proposition
\ref{p7.4} are satisfied.

We know that $[M_1^\#,M_2^\#,M_3^\#]$ and $[E_1^\#,E_2^\#,E_3^\#]$
is quasi-block data. Hence $[M_1^\#,M_2^\#,M_3^\#]$ is not diagonalizable.
If this triple possesses one non-trivial invariant subspace and all
$M_k^\#$ have only a single eigenvalue, then so have $E_k^\#$. It follows
that $n_1^\#=n_2^\#$, which is a contradiction (see (II) in the definition of
quasi-block data). As to condition (b2), assume that (b2) is not
satisfied. Then it is easy to see (using the nonresonance of $E_k$)
that $\nu=n_1^\#$.
Because $N={\rm trace\,}(E_1^\#+E_2^\#+E_3^\#)=n_1^\#+n_2^\#$
and $n_1^\#<n_2^\#$ by (II) of
the definition of quasi-block data, it follows that
$N>2n_1^\#$. On the other hand, $\ka_1=\nu$ and $\ka_2+\ka_3=N$.
We obtain that $\ka_1=n_1^\#$ and $\ka_2+\ka_3>2n_1^\#$.
Hence $\ka_2+\ka_3>2\ka_1$, which conflicts with $\ka_1\ge\ka_2\ge\ka_3$.
So we have proved that the assumptions of Proposition \ref{p7.4} are
satisfied.

Next we examine the question when two system with $A(z)$ and
$\wt{A}(z)$ given by (\ref{f7.48}) with $a_1(z)$ and $A_2(z)$
being introduced above, but with possibly different functions
$a_{12}(z)$ and $\tilde{a}_{12}(z)$, respectively, are equivalent.
Assume that \ben\label{f7.50} a_{12}(z) \;=\;
\left(\frac{p_{12}(z)}{p(z)},\;\frac{p_{13}(z)}{p(z)}\right),
\qquad \tilde{a}_{12}(z) \;=\;
\left(\frac{\tilde{p}_{12}(z)}{p(z)},\;\frac{\tilde{p}_{13}(z)}{p(z)}\right),
\een where $p_{12},\tilde{p}_{12}\in\cP_{1+\ka_1-\ka_2}$ and
$p_{13},\tilde{p}_{13}\in\cP_{1+\ka_1-\ka_3}$. This is the case
(see Theorem \ref{t3.5}) if there exists a function $V(z)$ of a
particular form such that $\wt{A}=VAV\iv+V'V\iv$.

We first claim that the $(2,1)$- and $(3,1)$-entries of $V(z)$,
$V_{21}$ and $V_{31}$, are equal to zero. This is obvious in the
case $\ka_1>\ka_2$. In case $\ka_1=\ka_2>\ka_3$, we have
$V_{31}=V_{32}=0$ and the upper left $2\times 2$ matrix of $V(z)$
is constant. In order to show that $V_{21}=0$, we write
$\wt{A}V=VA+V'$, and inspect the $(3,1)$-entries. It is equal to
zero for $V'$ and also for $VA$. The $(3,1)$-entry of $\wt{A}V$
equals $V_{21}$ times the $(3,2)$-entry of $\wt{A}(z)$. We
conclude that $V_{21}=0$ because otherwise the $(3,2)$-entry of
$\wt{A}(z)$ is zero, which contradicts the assumption that
$\wt{A}(z)$ is not of triangular form. In the case,
$\ka_1=\ka_2=\ka_3$, we have $\wt{A}(z)=VA(z)V\iv$ with a constant
matrix $V$. This can be rephrased in terms of the residues
$\wh{E}_k$ and $\wt{E}_k$ of $A(z)$ and $\wt{A}(z)$, respectively,
namely $\wt{E}_k=V\wh{E}_kV\iv$. The triples
$[\wh{E}_1,\wh{E}_2,\wh{E}_3]$ and $[\wt{E}_1,\wt{E}_2,\wt{E}_3]$
possess exactly one one-dimensional invariant subspace
$\C\oplus\{0\}\oplus\{0\}$. (If there were another, then we could
argue similar to above and conclude that the $2\times 2$ matrices
obtained from these triples by removing the first rows and columns
are reducible.) The matrix $V$ must map this invariant subspace
into itself. Hence it is also an invariant subspace of $V$, which
implies that $V_{21}=V_{31}=0$.

So we have shown that \bqn V(z) &=& \left(\ba{cc} v_{1} & v(z) \\
0 & V_2(z) \ea\right), \eqn with the same block structure as
$A(z)$ and $\wt{A}(z)$ and with entries having certain properties.
The next claim is that $V_2(z)$ is a constant scalar matrix. From
$\wt{A}=VAV\iv+V'V\iv$, it follows that
$A_2=V_2A_2V_2\iv+V_2'V_2\iv$, and we know that the $2\times 2$
system with $A_2(z)$ is not equivalent to a triangular system. If
the indices $\ka_2,\ka_3$ coincide, then $V_2(z)$ is a constant
matrix, and the above can be rephrased in terms of the residues $
\wh{E}_k^\#$ of $A_2(z)$. Namely,
$\wh{E}_k^\#=V_2\wh{E}_k^\#V_2\iv$. If $V_2$ is not a scalar
matrix, then we can argue similar as in the second paragraph of
the proof of Proposition \ref{p7.4} and conclude that the triple
$[\wh{E}_1^\#,\wh{E}_2^\#,\wh{E}_3^\#]$ is reducible, which
contradicts the assumption on $A_2(z)$. In case $\ka_2-\ka_3=1$,
$V_2(z)$ is of triangular form with constant values on the
diagonal. Now a straightforward analysis of
$A_2=V_2A_2V_2\iv+V_2'V_2\iv$ (using that $A_2(z)$ is not
triangular) implies that $V_2$ is a scalar constant matrix.

Having this information about $V(z)$, we can now reformulate the
condition $\wt{A}=VAV\iv+V'V\iv$ as follows. The systems $A(z)$
and $\wt{A}(z)$ are equivalent if and only if \bqn
\tilde{a}_{12}(z) &=& \lambda a_{12}(z)+
v(z)A_2(z)-a_1(z)v(z)+v'(z), \eqn for some
$\lambda\in\C\setminus\{0\}$ and some $v(z)=(v_{12}(z),v_{13}(z))$
where $v_{12}\in\cP_{\ka_1-\ka_2}$ and
$v_{13}\in\cP_{\ka_1-\ka_3}$. Referring to (\ref{f7.50}) this can
be expressed as \bqn (\tilde{p}_{12},\tilde{p}_{13}) &=&
\lambda(p_{12},p_{13})+
p(v_{12},v_{13})A_2-a_1p(v_{12},v_{13})+p(v_{12}',v_{13}'). \eqn
Hence the two systems with $(p_{12},p_{13})$ and
$(\tilde{p}_{12},\tilde{p}_{13})$ are equivalent if and only if
there exists $\lambda\in\C\setminus\{0\}$ such that \bqn
(\tilde{p}_{12},\tilde{p}_{13})-\lambda(p_{12},p_{13})\;\in\;\im\Xi,
\eqn where $\im\Xi$ is the image of the linear mapping defined by
\bqn \Xi &:& \cP_{\ka_1-\ka_2}\oplus\cP_{\ka_1-\ka_3} \to
\cP_{1+\ka_1-\ka_2}\oplus\cP_{1+\ka_1-\ka_3} \nn\\[.5ex] &&
(v_{12},v_{13})\mapsto
p(v_{12},v_{13})A_2-a_1p(v_{12},v_{13})+p(v_{12}',v_{13}'). \eqn
Now we decompose
$\cP_{1+\ka_1-\ka_2}\oplus\cP_{1+\ka_1-\ka_3}=\im\Xi\oplus X$ as a
direct sum. We arrive at the following statement: If
$(p_{12},p_{13})\in X$ and $(\tilde{p}_{12},\tilde{p}_{13})\in X$,
then the corresponding systems are equivalent if and only if
$(\tilde{p}_{12},\tilde{p}_{13})=\lambda(p_{12},p_{13})$ for some
$\lambda\in\C\setminus\{0\}$.

Obviously, we have $\dim X\geq 2$. However, we claim that $\dim
X\geq 3$ in the case where the $2\times 6$ matrix appearing in
(\ref{f7.46}) has rank $1$, i.e., if $[M_1^\#,M_2^\#,M_3^\#]$ is
equivalent to a triple (\ref{f7.46a}). Indeed, this is shown as
soon as one has shown that the kernel of $\Xi$ is non-trivial. For
this we consider the $2\times 2$ system
$\wh{Y}_2'(\tz)=\wh{A}_2(z)\wh{Y}_2(\tz)$ with
$\wh{A}_2(z)=A_2(z)-a_1(z)I_2$. It is easy to see that this system
is also of standard form with indices $\ka_2-\ka_1$ and
$\ka_3-\ka_1$, and with data
$[M_1^\#\mu_1\iv,M_2^\#\mu_2\iv,M_3^\#\mu_3\iv]$ and
$[E_1^\#-\ep_1I_2,E_2^\#-\ep_2I_2,E_3^\#-\ep_3I_2]$. The solution
of this system is given by $\wh{Y}_2(\tz)=Y_2(\tz)y_1\iv(\tz)$.
{From} $\wh{Y}_2'=\wh{A}_2\wh{Y}_2$ it follows immediately that
$\wh{Y}_2\iv\wh{A}_2+(\wh{Y}_2\iv)'=0$ by multiplying on both
sides with $\wh{Y}_2\iv$. We may assume that the monodromy
representation of $\wh{Y}_2(\tz)$ is given by \bqn \left[
\left(\ba{cc} \hat{\mu}_1\mu_1\iv &\hat{\alpha}_1\mu_1\iv \\ 0 & 1
\ea\right), \left(\ba{cc} \hat{\mu}_2\mu_2\iv
&\hat{\alpha}_1\mu_2\iv \\ 0 & 1 \ea\right), \left(\ba{cc}
\hat{\mu}_3\mu_3\iv &\hat{\alpha}_1\mu_3\iv \\ 0 & 1 \ea\right)
\right]. \eqn Now we define $(v_{12},v_{13})$ to be the last row
in $\wh{Y}_2\iv$. Obviously, this vector function satisfies
$(v_{12},v_{13})(A_2(z)-a_1(z)I_2)+(v_{12}',v_{13}')=0$. What
remains to show is that $v_{12}$ and $v_{13}$ are polynomials of
appropriate degree. From the above monodromy representation for
$\wh{Y}_2$, it follows that $(v_{12},v_{13})$ is a single valued
function. This function is analytic on
$\C\setminus\{a_1,a_2,a_3\}$. Because of the non-resonance of
$[E_1,E_2,E_3]$ and because $\mu_k$ is an eigenvalue of $M_k^\#$,
it follows that $\ep_k$ is an eigenvalue of $E_k^\#$. Hence the
matrices $E_k^\#-\ep_kI_2$ have an eigenvalue equal to zero and
are also non-resonant. But these matrices describe the local
behavior of $\wh{Y}_2$ (and of its inverse) near the point $a_k$.
A straightforward computation (using the fact that
$(v_{12},v_{13})$ is single valued) shows that $(v_{12},v_{13})$
is actually analytic at $z=a_k$. Hence we have shown that
$(v_{12},v_{13})$ is entire analytic. That $v_{12}$ and $v_{13}$
are polynomials of degree not greater than $\ka_1-\ka_2$ and
$\ka_1-\ka_3$, respectively, can be obtained by employing the
behavior of $\wt{Y}_2$ at infinity, which is described by the
integers $\ka_2-\ka_1$ and $\ka_3-\ka_1$. In summary, we can
conclude that $(v_{12},v_{13})\neq0$ lies in the kernel of $\Xi$.

Similar as in the proof of Theorem \ref{t6.2}, we can define a linear mapping
$\Lambda:X\to \C^{1\times 6}$ from the space $X$ into the set of all
``admissible'' vectors $\alpha=(\alpha_1,\alpha_2,\alpha_3)$
(see also (\ref{f7.45})). This mapping is
characterized by the following property. For a system of the form (\ref{f7.48})
with above defined $a_1(z)$ and $A_2(z)$ and with $a_{12}$ given by
(\ref{f7.50}), where $(p_{12},p_{13})\in X$, the corresponding monodromy is
given by (\ref{f7.41}), where $\alpha_1,\alpha_2,\alpha_3$ are given by
$\alpha=(\alpha_1,\alpha_2,\alpha_3)=\Lambda((p_{12},p_{13}))$.

We have a one-to-one correspondence between equivalence classes of systems
of standard form and equivalence classes of data (see Corollary \ref{c6.3}).
In our situation we can apply Proposition \ref{p7.4} (see also
(\ref{f7.46})) and what has been stated above about the equivalence of the
systems under consideration. We obtain the following statement.
For $(p_{12},p_{13}),(\tilde{p}_{12},\tilde{p}_{13})\in X$, the vectors
$\Lambda((p_{12},p_{13}))$ and $\Lambda((\tilde{p}_{12},\tilde{p}_{13}))$ are
equivalent in the sense of (\ref{f7.46}) (or, (\ref{f7.43})) if and only if
there exists a $\lambda\in\C\setminus\{0\}$ such that
$(\tilde{p}_{12},\tilde{p}_{13})=\lambda(p_{12},p_{13})$.

The kernel of $\Lambda$ is trivial. Indeed,
$\Lambda((p_{12},p_{13}))=0=\Lambda(0)$ implies immediately
$(p_{12},p_{13})=0$. Hence $\dim\im\Lambda\geq 2$ or even $\dim\im\Lambda\geq
3$ depending on whether the $2\times 6$ matrix appearing in (\ref{f7.46})
has rank $2$ or $1$. So we are done as soon as we have shown that
the sum of $\im\Lambda$ and the following linear subspace of $\C^{1\times 6}$,
\bqn\label{f7.subsp}
\left\{\;
c\,\Big(\ba{ccc}
M_1^\#-\mu_1I_2, & M_2^\#-\mu_2I_2, & M_3^\#-\mu_3I_2 \ea\Big)\;:\;
c\in\C^{1\times 2}\;\right\},
\eqn
exhausts the whole four-dimensional subspace of $\C^{1\times 6}$
consisting of all admissible vectors
$\alpha$ (see again (\ref{f7.46})).
As we have estimates for the dimensions, this follows from the statement that
the intersection of (\ref{f7.subsp}) with $\im\Lambda$ contains only the
zero vector. Indeed, if $\alpha=\Lambda((p_{12},p_{13}))$ is
contained in the above set, then $\alpha$ is equivalent to $0=\Lambda(0)$,
whence again follows that $(p_{12},p_{13})=0$. Hence $\alpha=0$.
\end{proof}

Finally, we are going to consider data and systems which are also of block
form, but where the block partition of the matrices are $(2,1)$ rather than
$(1,2)$. The proofs of the corresponding statements are analogous, and will
therefore be omitted. However, we will state definitions and the results
completely.

Let $[M_1,M_2,M_3]$ and $[E_1,E_2,E_3]$ be admissible data of
$3\times 3$ matrices. We call this data {\em $(2,1)$-quasi-block
data} if $[M_1,M_2,M_3]$ is equivalent to a triple
\ben\label{f7.41+} \left[ \left(\ba{cc}M_1^{\#} & \alpha_1 \\ 0 &
\mu_1 \ea \right), \left(\ba{cc}M_2^{\#} & \alpha_2 \\ 0 & \mu_2
\ea \right), \left(\ba{cc}M_3^{\#} & \alpha_3 \\ 0 & \mu_3 \ea
\right) \right] \een and $[M_1^{\#},M_2^{\#},M_3^{\#}]$ and
$[E_1^{\#},E_2^{\#},E_3^{\#}]$ is quasi-block data. Here $\mu_k$
are numbers, $\alpha_k$ are $2\times 1$ vectors, and $M_k^{\#}$
are $2\times 2$ matrices. The matrices $E_k^{\#}$ and numbers
$\ep_k$ are supposed to satisfy $M_k^{\#}\sim\exp(-2\pi
iE_k^{\#})$ and $\mu_k=\exp(-2\pi i \ep_k)$ and are chosen in such
a way that the eigenvalues of $E_k^{\#}$ and the number $\ep_k$
are the exactly eigenvalues of $E_k$. For given $M_k^{\#}$ and
$\mu_k$, the number $\ep_k$ is uniquely determined and so is the
matrix $E_k^{\#}$ up to similarity (because of the non-resonance of
$E_k$). Again we define the integers \bqn\label{f7.42+} \nu
\;\;=\;\; \ep_1+\ep_2+\ep_3,\qquad N \;\;=\;\; {\rm
trace\,}(E_1^{\#}+E_2^{\#}+E_3^{\#}). \eqn As above, we remark
that if $[M_1,M_2,M_3]$ possesses two different two-dimensional
invariant subspaces, then $[M_1,M_2,M_3]$ is equivalent to
different triples of the form (\ref{f7.41+}) (with different
$\mu_k$ and $M_k^{\#}$ in particular). In general, we obtain also
different values for $\ep_k, E_k^{\#}, \nu$ and $N$.

The counterpart to Proposition \ref{p7.4} is the following result.

\begin{proposition}\label{p7.6}
Let $[M_1,M_2,M_3]$ be a triple of the form (\ref{f7.41+}). Assume that
\begin{itemize}
\item[(a)]
$[M_1^{\#},M_2^{\#},M_3^{\#}]$ is not
diagonalizable;
\item[(b)]
If  $[M_1^{\#},M_2^{\#},M_3^{\#}]$ possess exactly one non-trivial invariant
subspace, then
\begin{itemize}
\item[(b1)]
for at least one $k$, $M_k^\#$ possesses two different eigenvalues;
\item[(b2)]
for at least one $k$, $\mu_k$ is not equal to the eigenvalue of $M_k^\#$
which does not correspond to this invariant subspace.
\end{itemize}
\end{itemize}
Let $[\wt{M}_1,\wt{M}_2,\wt{M}_3]$ be another triple of the same form
(\ref{f7.41+}) but with $\alpha_1,\alpha_2,\alpha_3$ replaced by
$\wt{\alpha}_1,\wt{\alpha}_2,\wt{\alpha}_3$. Then these triples are equivalent
if and only if there exist $\lambda\in\C\setminus\{0\}$ and
$c\in\C^{2\times 1}$ such that
\bqn\label{f7.43+}
\alpha_k &=& \lambda\wt{\alpha}_k+M_k^{\#}c-c\mu_k
\qquad\quad\mbox{ for each }k=1,2,3.
\eqn
\end{proposition}

The set of equivalence classes for triples of the form
(\ref{f7.41+}) with $\mu_k$ and $M_k^{\#}$ being fixed can be
identified similar as before. All such triples are parameterized
by $2\times 1$ vectors $\alpha_1,\alpha_2,\alpha_3$, and we
introduce the $6\times 1$ vector
$\alpha=(\alpha_1^T,\alpha_2^T,\alpha_3^T)^T$. We have to take the
linear condition (\ref{f7.6}) into account. It can be restated as
\bqn\label{f7.45+} \Big(\ba{ccc} \mu_2\mu_3I_2, & M_1^\#\mu_3, &
M_1^\# M_2^\#\ea\Big) \alpha =0. \eqn The $2\times 6$ matrix
appearing here has rank $2$. Hence the vector $\alpha$ is again
subject to two linear conditions, and all ``admissible'' vectors
$\alpha$ are taken from a $4$-dimensional linear subspace of
$\C^{6\times 1}$. The equivalence relation (\ref{f7.43+}) can be
rewritten as \bqn\label{f7.46+} \alpha &=& \lambda\wt{\alpha}
\;+\;\left(\ba{c} M_1^\#-\mu_1I_2 \\ M_2^\#-\mu_2I_2 \\
M_3^\#-\mu_3I_2 \ea\right)c \eqn with $\lambda\in\C\setminus\{0\}$
and $c\in\C^{2\times 1}$. The $6\times 2$ matrix appearing here
has again either rank $2$ or rank $1$. It has rank one if and only
if $[M_1^\#,M_2^\#,M_3^\#]$ possesses exactly one non-trivial
invariant subspace and, for each $k=1,2,3$, the number $\mu_k$ is
equal to the eigenvalue of $M_k^\#$ which corresponds to this
invariant subspace. In conclusion, we obtain that
 the set of equivalence classes of the above triples can be
identified with either $\sP^1\cup\{0\}$ or with $\sP^2\cup\{0\}$.
The single equivalence class containing the vector $\alpha=0$
consists of the  triples which are ``block-diagonalizable''.

Next we consider systems $Y'(\tz)=A(z)Y(\tz)$ of standard form with indices
$\ka_1\ge\ka_2\ge\ka_3$ for which $A(z)$ can be written as
\bqn\label{f7.47+}
A(z) &=& \frac{1}{p(z)}
\left(\ba{ccc}
\ka_1z^2+p_1(z) & p_{12}(z) & p_{13}(z) \\
p_{21}(z) & \ka_2z^2+p_2(z) & p_{23}(z) \\
0 & 0 & \ka_3z^2 + p_3(z) \ea\right),
\eqn
where $p(z)=(z-a_1)(z-a_2)(z-a_3)$, $p_k\in\cP_1$,
$p_{jk}\in\cP_{1+\ka_j-\ka_k}$. If such a system is not equivalent to a system
of the form (\ref{f7.10}), then it will be called a {\em system of
$(2,1)$-block form}.
Note that necessarily $\ka_2-\ka_3\leq 1$ because otherwise
$p_{32}(z)=0$ and the system is of the form (\ref{f7.10}).

The data which is associated to such systems is described in the following
theorem.

\begin{theorem}\label{t7.7}
Let $\ka_1,\ka_2,\ka_3\in\Z$, $\ka_1\ge\ka_2\ge\ka_3$, and assume
$\ka_1-\ka_2\leq 1$. Let $[M_1,M_2,M_3]$ and $[E_1,E_2,E_3]$ be admissible
data of $3\times 3$ matrices. Then the following two statements are
equivalent:
\begin{itemize}
\item[(i)]
$[M_1,M_2,M_3]$ and $[E_1,E_2,E_3]$ is the data associated to a system of
$(2,1)$-block form with indices $\ka_1,\ka_2,\ka_3$.
\item[(ii)]
$[M_1,M_2,M_3]$ and $[E_1,E_2,E_3]$ is $(2,1)$-quasi-block data, and
$N=\ka_1+\ka_2$, $\nu=\ka_3$, where $\nu$ and $N$ are defined by
(\ref{f7.42+}).
\end{itemize}
\end{theorem}

We remark that the cases of $(1,2)$-block form systems and of
$(2,1)$-block form systems are not exclusive. They overlap exactly
if $\ka_1=\ka_2=\ka_3$ and if the system is equivalent to a system
of a ``block-diagonal form''. This can be deduced from the
argumentation in the proof of the main results given below.

\subsection{The main results in the case {\boldmath $n=m=3$}}

The main results for the case $m=n=3$, i.e., the identification of
the indices corresponding to given data is given in the following
theorem. We have to resort to the classification of the
reducibility type of $[M_1,M_2,M_3]$ and the definition of the
various integers $\ka$, $\nu$, $N$, $n_1$, $n_2$, $n_3$, $\nu_1$,
$\nu_2$, $\nu_{\#}$ stated at the beginning of this section.
Recall that we have assumed $n_1\geq n_2$ in case (C-1), $n_2\geq
n_3$ in case (C-2), and $n_1\geq n_2\geq n_3 $ in case (D) without
lost of generality. In the other cases, the integers of interest
were defined uniquely.

\begin{theorem}
Let $[M_1,M_2,M_3]$ be admissible data of $3\times 3$ matrices.
Then in the following cases the indices $\ka_1,\ka_2,\ka_3\in\Z$,
$\ka_1\geq\ka_2\geq\ka_3$, are as follows.
\begin{itemize}
\item[(1)]
In cases (B-1) or (B-3) with $2\nu\geq N$, the indices are\\[.5ex]
\hspace*{6ex}
$\ka_1=\nu$, $\ka_2=N/2$, $\ka_3=N/2$ if $N$ is even, or,\\
\hspace*{6ex}
$\ka_1=\nu$, $\ka_2=(N+1)/2$, $\ka_3=(N-1)/2$ if $N$ is odd.
\item[(2)]
In cases (B-2) or (B-3) with $2\nu\leq N$, the indices are\\[.5ex]
\hspace*{6ex}
$\ka_1=N/2$, $\ka_2=N/2$, $\ka_3=\nu$ if $N$ is even, or,\\
\hspace*{6ex}
$\ka_1=(N+1)/2$, $\ka_2=(N-1)/2$, $\ka_3=\nu$ if $N$ is odd.
\item[(3)]
In case (C) with $n_1\geq n_2 \geq n_3$ or in case (C-1) with $n_2\geq n_3$
or in case (C-2) with $n_1\geq n_2$ or in case (D), the indices are \\[.5ex]
\hspace*{6ex}
$\ka_1=n_1$, $\ka_2=n_2$, $\ka_3=n_3$.
\item[(4)]
In case (C-3) with $\nu_1\geq\nu_2$, the indices are\\[.5ex]
\hspace*{6ex}
$\ka_1=\nu_1$, $\ka_2=\nu_2$, $\ka_3=\nu_{\#}$ if $\nu_2\geq\nu_{\#}$,
or,\\
\hspace*{6ex}
$\ka_1=\nu_1$, $\ka_2=\nu_{\#}$, $\ka_3=\nu_2$ if $\nu_1\geq\nu_{\#}\geq\nu_2$,
or,\\
\hspace*{6ex}
$\ka_1=\nu_{\#}$, $\ka_2=\nu_1$, $\ka_3=\nu_2$ if $\nu_{\#}\geq\nu_1$.
\item[(5)]
In cases (C) or (C-1) with $n_2<n_3$ and $2n_1\geq n_2+n_3$,
the indices are\\[.5ex]
\hspace*{6ex}
$\ka_1=n_1$, $\ka_2=(n_2+n_3)/2$, $\ka_3=(n_2+n_3)/2$ if $n_2+n_3$ is even,
or,\\ \hspace*{6ex}
$\ka_1=n_1$, $\ka_2=(n_2+n_3+1)/2$, $\ka_3=(n_2+n_3-1)/2$ if $n_2+n_3$ is odd.
\item[(6)]
In cases (C) or (C-2) with $n_1<n_2$ and $n_1+n_2\geq 2n_3$,
the indices are\\[.5ex]
\hspace*{6ex}
$\ka_1=(n_1+n_2)/2$, $\ka_2=(n_1+n_2)/2$, $\ka_3=n_3$ if $n_1+n_2$ is even,
or,\\ \hspace*{6ex}
$\ka_1=(n_1+n_2+1)/2$, $\ka_2=(n_1+n_2-1)/2$, $\ka_3=n_3$ if $n_1+n_2$ is odd.
\item[(7)]
In case (C-3) with $\nu_1<\nu_2$ and $2\nu_{\#}\geq \nu_1+\nu_2$, the
indices are \\[.5ex]
\hspace*{6ex}
$\ka_1=\nu_{\#}$, $\ka_2=(\nu_1+\nu_2)/2$, $\ka_3=(\nu_1+\nu_2)/2$,
if $\nu_1+\nu_2$ is even, or,\\
\hspace*{6ex}
$\ka_1=\nu_{\#}$, $\ka_2=(\nu_1+\nu_2+1)/2$, $\ka_3=(\nu_1+\nu_2-1)/2$,
if $\nu_1+\nu_2$ is odd.
\item[(8)]
In case (C-3) with $\nu_1<\nu_2$ and $\nu_1+\nu_2\geq 2\nu_{\#}$, the indices are
\\[.5ex]
\hspace*{6ex}
$\ka_1=(\nu_1+\nu_2)/2$, $\ka_2=(\nu_1+\nu_2)/2$, $\ka_3=\nu_{\#}$
if $\nu_1+\nu_2$ is even, or,\\
\hspace*{6ex}
$\ka_1=(\nu_1+\nu_2+1)/2$, $\ka_2=(\nu_1+\nu_2-1)/2$, $\ka_3=\nu_{\#}$
if $\nu_1+\nu_2$ is odd.
\end{itemize}
In the remaining cases, i.e.,
\begin{itemize}
\item[(9)] in case (A),
\item[(10)] in case (B-1) with $2\nu<N$,
\item[(11)] in case (B-2) with $2\nu>N$,
\item[(12)] in case (C) with $n_1+n_2<2n_3$ and $2n_1<n_2+n_3$,
\item[(13)] in case (C-1) with $2n_1<n_2+n_3$,
\item[(14)] in case (C-2) with $n_1+n_2<2n_3$,
\end{itemize}
the indices are determined by the following statements:
\begin{itemize}
\item[(a)]
If $\ka\equiv 1\;{\rm mod}\, 3$, then $\ka_1=(\ka+2)/3$,
$\ka_2=\ka_3=(\ka-1)/3$.
\item[(b)]
If $\ka\equiv -1\;{\rm mod}\, 3$, then $\ka_1=\ka_2=(\ka+1)/3$,
$\ka_3=(\ka-2)/3$.
\item[(c)]
If $\ka\equiv 0\;{\rm mod}\, 3$, then either $\ka_1=\ka_2=\ka_3=\ka/3$ or
$\ka_1=\ka/3+1$, $\ka_2=\ka/3$, $\ka_3=\ka/3-1$.
\item[(c*)]
If $\ka\equiv 0\;{\rm mod}\, 3$, then the indices are
$\ka_1=\ka/3+1$, $\ka_2=\ka/3$, $\ka_3=\ka/3-1$ if and only if
$[M_1,M_2,M_3]$ is the monodromy of some third order linear differential
equation with Fuchsian singularities $a_1,a_2,a_3$ and local
exponents given by $\{\ep_k^{(1)},\ep_k^{(2)},\ep_k^{(3)}\}$
for $z=a_k$, $k=1,2$, and
$\{\ep_3^{(1)}+1-\ka/3,\ep_k^{(2)}+1-\ka/3,\ep_k^{(3)}+1-\ka/3\}$
for $z=a_3$, where $\ep_k^{(1)},\ep_k^{(2)},\ep_k^{(3)}$ are the eigenvalues of
$E_k$, $k=1,2,3$.
\end{itemize}
\end{theorem}
\begin{proof}
First of all, the reader can easily verify that the classification into the
cases (1)--(14) is complete. Some of these cases do overlap
(namely, case (B-3) with $2\nu=N$ is stated in (1) and (2), case (C-3) with
$\nu_1\geq\nu_2$ and $\nu_1=\nu_\#$ or $\nu_2=\nu_\#$ is stated multiply
in (4), and case (C-3) with
$\nu_1<\nu_2$ and $\nu_1+\nu_2=2\nu_\#$ is given in both (7) and (8)),
but in these cases, the description of the partial indices is the same.
We remark in particular that the cases given in (9)--(14) are exactly the
cases which are not contained in (1)-- (8).

In Theorem \ref{t7.2} we have completely characterized the data which is
associated to systems of the form (\ref{f7.10}).
We obtain that given data is associated to such a system if and only if
$[M_1,M_2,M_3]$ is equivalent to a certain triple (\ref{f7.7}) such that for
the corresponding integers $n_1,n_2,n_3$ (which are defined by resorting to
the properly numbered eigenvalues of $E_1,E_2,E_3$) the inequality $n_1\geq
n_2\geq n_3$ holds. This inequality follows from the statement that
$\ka_k=n_k$ for $k=1,2,3$, and from the ordering $\ka_1\geq\ka_2\geq\ka_3$.
We note that in some cases, there exist essentially different ways to
establish an equivalence between $[M_1,M_2,M_3]$ and a triangular triple
(\ref{f7.7}), which in turn results into different definitions of
$n_1,n_2,n_3$. The statement (i) of Theorem \ref{t7.2} is true if and only
if such a relationship exists.

Examining the different reducibility types for a triple
(\ref{f7.7}) under the condition $n_1\geq n_2\geq n_3$ we arrive
at the following cases. First of all, only the reducibility types
(C)--(D) are covered. Case (C) is covered exactly for $n_1\geq
n_2\geq n_3$ as the enumeration of these integers is unique. Case
(C-1) is covered exactly for $n_2 \geq n_3$ as $n_3$ is defined
uniquely and $n_1\geq n_2$ without loss of generality. Similarly,
case (C-2) is covered exactly for $n_1\geq n_2$ as $n_1$ is unique
and $n_2\geq n_3$ without loss of generality. The case (D) is
obviously covered completely. So we arrive at the statement (3).

However, there are some other possibilities leading to case (C-3)
and statement (4). In case (C-3), the integers $\nu_1,\nu_2,\nu_3$
are uniquely defined, but an equivalence with a triple
(\ref{f7.7}) can be established in three different ways. The first
way is if the triple has entries $\beta_k=\gamma_k=0$. This gives
$n_1=\nu_1$, $n_2=\nu_2$, $n_3=\nu_\#$, and hence the first
subcase in (4). The second possibility is if the triple has
entries $\alpha_k=\beta_k=0$. Here we obtain $n_1=\nu_1$,
$n_2=\nu_\#$, $n_3=\nu_2$, which is the second subcase in (4).
Finally, we may have $\alpha_k=\gamma_k=0$, which leads to
$n_1=\nu_\#$, $n_2=\nu_1$, $n_3=\nu_2$, and the third subcase in
(4).

We have now identified all the data associated to system (\ref{f7.10})
in terms of the classification established at the beginning of this section.

Next we do the same for the systems of $(1,2)$-block form (see
(\ref{f7.47})). The data of such systems has been described in
Theorem \ref{t7.5}. With the numbers $\nu$ and $N$ defined by
(\ref{f7.42}) we obtain that the indices are $\ka_1=\nu$,
$\ka_2=[(N+1)/2]$, and $\ka_3=[N/2]$. The latter follows from the
fact that $\ka_2+\ka_3=N$ and $0\leq\ka_2-\ka_3\leq 1$. The
ordering $\ka_1\geq\ka_2\geq\ka_3$ can be rephrased by
$\nu\geq[(N+1)/2]$, which is equivalent to $2\nu\geq N$. We arrive
at conclusion that the data of such systems is precisely the data
which is of $(1,2)$-quasi-block form and satisfies $2\nu\geq N$.

Now we have to examine what does this mean in terms of the
classification given at the beginning of this section. In the
triple (\ref{f7.41}), there appears a triple
$[M_1^\#,M_2^\#,M_3^\#]$ which along with $[E_1^\#,E_2^\#,E_3^\#]$
has to be quasi-block data. We have to distinguish the cases (I)
and (II) stated in the definition of quasi-block data.

In case (I), the triple (\ref{f7.41}) possesses exactly one one-dimensional
invariant subspace $\C\oplus\{0\}\oplus\{0\}$ and no two-dimensional invariant
subspace containing it. Hence the data $[M_1,M_2,M_3]$ and $[E_1,E_2,E_3]$ is
precisely the data of reducibility type (B-1) or (B-3) with the condition
$2\nu\geq N$, where $\nu$ and $N$ are now defined by (\ref{f7.4}). In fact,
these numbers are the same as the former ones defined in (\ref{f7.47}).
So we arrive at the statement (1).

In case (II), we may assume $[M_1^\#,M_2^\#,M_3^\#]$ to be of triangular form.
We have $N=n_1^\#+n_2^\#$ and $n_1^\#<n_2^\#$ with the integers $n_1^\#$ and
$n_2^\#$ defined in terms of the eigenvalues of $[E_1^\#,E_2^\#,E_3^\#]$.
Hence the data is given by triples (\ref{f7.47}) with $2\nu\geq n_1^\#+n_2^\#$
and $n_1^\#<n_2^\#$. These triples are of reducibility type (C), (C-1) or
(C-3). We do not obtain (C-2) or (D) because of the
non-diagonalizability of $[M_1^\#,M_2^\#,M_3^\#]$. In case of reducibility
type (C), we have $n_1=\nu$, $n_2=n_1^\#$ and $n_3=n_2^\#$, which gives the
condition $2n_1\geq n_2+n_3$ and $n_2<n_3$. In case (C-1) we have the same
identification, but in addition we may also have (interchanging $n_1$ and
$n_2$) $n_2=\nu$, $n_1=n_1^\#$ and $n_3=n_2^\#$. This would give the condition
$2n_2\geq n_1+n_3$ and $n_1<n_3$, which, however, conflicts with the assumption
$n_1\geq n_2$. So we arrive at the statement (5).

In case (C-3), the only possibility not contradicting the
non-diagonalizability of $[M_1^\#,M_2^\#,M_3^\#]$ is if the triple
$[M_1,M_2,M_3]$ is equivalent to a triple (\ref{f7.7}) with
$\alpha_k=\gamma_k=0$. This results in the identification $\nu_\#=\nu$,
$\nu_1=n_1^\#$ and $\nu_2=n_2^\#$. Hence we obtain the conditions $2\nu_\#\geq
\nu_1+\nu_2$ and $\nu_1<\nu_2$, and this gives the statement (7).

So we have also settled the cases of data corresponding to systems
of $(1,2)$-block form.

The analysis of the data corresponding to systems of $(2,1)$-block
form (\ref{f7.47+}) is completely analogous and will therefore be
omitted. Here we arrive at the statements (2), (6) and (8) instead
of (1), (5) and (7).

Finally, we have to examine the indices for the data in the remaining cases.
We know that this data is associated to certain systems which are neither of
$(1,2)$-block form, $(2,1)$-block form nor of triangular form.

It follows that $\ka_1-\ka_2\le1$ because otherwise (by Theorem \ref{t4.3})
the $(2,1)$- and $(3,1)$-entries of the matrix $A(z)$ of the system vanish.
Similarly, it follows that $\ka_2-\ka_3\le 1$. From the fact that
$\ka_1\geq\ka_2\geq\ka_3$ and $\ka_1+\ka_2+\ka_3=\ka$, we immediately
obtain the assertions (a), (b) and (c).

In order to prove (c*), we have to examine when the indices are
such that $\ka_1-\ka_2=\ka_2-\ka_3=1$. Again by Theorem
\ref{t4.3}, this is exactly the case if the system is of the form
\bqn A(z) &=& \frac{1}{p(z)} \left(\ba{ccc} (\ka/3+1)z^2+p_1(z) &
p_{12}(z) & p_{13}(z) \\ \alpha & \ka/3z^2+p_2(z) & p_{23}(z) \\ 0
& \beta & (\ka/3-1)z^2+p_3(z) \ea\right),\nn \eqn where
$p(z)=(z-a_1)(z-a_2)(z-a_3)$, $p_k\in\cP_1$,
$p_{jk}\in\cP_{1+\ka_j-\ka_k}$ and $\alpha,\beta\in\C$. If
$\alpha=0$ or $\beta=0$, then $A(z)$ is of block-triangular or
even triangular form, which has been excluded. Hence
$\alpha\beta\neq0$. (In fact, if $\alpha\beta\neq0$, then $A(z)$
is not equivalent to a system of triangular or block-triangular
form.) Applying Theorem \ref{t5.4}, we obtain the desired
assertion (c*).
\end{proof}

In the cases where assumption (c) of the previous theorem applies,
we are led to the monodromy problem for third order Fuchsian
differential equations with three singular points.  Here the
monodromy does not depend on the location of the singularities,
and neither do the corresponding indices. Again, as far as the
authors know, no explicit answer exists to this monodromy problem.
However, some particular cases can be tackled directly by help of
the remark made in the last paragraph of Section \ref{s5}. Namely,
if one of the matrices $M_1,M_2$ or $M_3$ has the property that
$\mdg(M_k)\leq 2$, then the indices are necessarily
$\ka_1=\ka_2=\ka_3=\ka/3$. We note that in the cases (9)--(14),
matrices $[M_1,M_2,M_3]$ with this property really exist.

\vspace*{3ex} \noindent {\bf Acknowledgments.} The authors want to
thank  A.~A.~Bolibruch and A.~R.~Its for fruitful discussions with
one of them (T.E.) about various aspects of the monodromy and the
Riemann-Hilbert problem during their stay at MSRI.

\end{document}